\documentclass[12pt]{amsart}
\usepackage{amssymb}
\usepackage{amsmath}
\usepackage{amsthm}
\usepackage{amstext}
\usepackage{amsopn}
\usepackage{mathrsfs}
\usepackage{latexsym}
\usepackage{mathabx}
\allowdisplaybreaks
\newtheorem{Definition}{Definition}[section]
\newtheorem{Proposition}{Proposition}[section]
\newtheorem{Lemma}{Lemma}[section]
\newtheorem{Theorem}{Theorem}[section]
\newtheorem{Corollary}{Corollary}[section]
\newtheorem{Remark}{Remark}[section]
\newtheorem{Example}{Example}[section]

\begin{document}
\bibliographystyle{plain}
\footnotetext{
This work is partially supported by MNiSW research grant No 1 P03A 013 30}
\title[Operators related to subordination]
{Operators related to subordination 
for free multiplicative convolutions}
\author[R. Lenczewski]{Romuald Lenczewski}
\address{Romuald Lenczewski, \newline
Instytut Matematyki i Informatyki, Politechnika Wroc\l{}awska, \newline
Wybrze\.{z}e Wyspia\'{n}skiego 27, 50-370 Wroc{\l}aw, Poland  \vspace{10pt}}
\email{Romuald.Lenczewski@pwr.wroc.pl}
\maketitle
\begin{abstract}
It has been shown by Voiculescu and Biane that 
the analytic subordination property holds for free additive
and multiplicative convolutions.
In this paper, we present an operatorial approach to
subordination for free multiplicative convolutions. 
This study is based on the concepts of `freeness with subordination', or 
`s-free independence', and `orthogonal independence', introduced recently in the context of
free additive convolutions. In particular, we introduce and study
the associated multiplicative convolutions and construct 
related operators, called `subordination operators' and 
`subordination branches'. Using orthogonal independence,
we derive decompositions of subordination branches
and related decompositions of s-free and free multiplicative
convolutions. The operatorial methods lead to several new 
types of graph products, called `loop products',
associated with different notions of independence (monotone, boolean, orthogonal, 
s-free). We also prove that the enumeration of rooted 
`alternating double return walks' on the loop products of graphs 
and on the free product of graphs gives the moments of 
the corresponding multiplicative convolutions.
\end{abstract}
\section{Introduction}
Multiplication of free random variables $X_1$ and $X_2$ with distributions $\mu_1$ and $\mu_2$, 
respectively, leads to the multiplicative convolution $\mu_1 \boxtimes \mu_2$, introduced by 
Voiculescu [23] in the $C^{*}$--algebra framework, which gives the distribution of the product of $X_{1}$ and $X_{2}$
(the general case of measures with unbounded support was studied by Bercovici and Voiculescu [7]).

Let ${\mathcal  M}_{{\mathbb R}_{+}}$ denote the set
of probability measures on ${\mathbb R}_{+}=[0, \infty)$.
If $\mu\in{\mathcal M}_{{\mathbb R}_{+}}$, we can define
\begin{equation}\tag{1.1}
\psi_{\mu}(z)=\int_{{\mathbb R}_{+}}\frac{zt}{1-zt}d\mu(t), \;\;\; z\in {\mathbb C}\setminus 
{\mathbb R}_{+},
\end{equation}
which, in the case when $\mu$ has finite moments of all orders, becomes the 
moment generating function $\psi_{\mu}(z)=\sum_{n=1}^{\infty}\mu(X^n)z^{n}$, 
where $\mu(X^n)$ are the moments
of the unique functional $\mu:{\mathbb C}[X]\rightarrow {\mathbb C}$ defined by $\mu$.
In order to study $\mu_1 \boxtimes \mu_2$, Voiculescu introduced the S-transform of $\mu$ 
defined by $S_{\mu}(z)=(1+z)\psi_{\mu}^{-1}(z)/z$, where $\psi_{\mu}^{-1}(z)$ denotes the inverse
of $\psi_{\mu}(z)$ with respect to composition. The key multiplicative formula for the S-transforms
is given by $S_{\mu_1\,\boxtimes \,\mu_2}(z)=S_{\mu_1}(z)S_{\mu_2}(z)$.

In our approach, a central role is played by the transform related to $\psi_{\mu}(z)$, namely 
\begin{equation}\tag{1.2}
\eta_{\mu}(z)=\frac{\psi_{\mu}(z)}{1+\psi_{\mu}(z)}
\end{equation}
where $z\in {\mathbb C}\setminus {\mathbb R}_{+}$. Using this transform, 
Biane [8] has proved the subordination property for free multiplicative convolutions
of probability measures on ${\mathbb R}_{+}$ (and also for probability 
measures on the unit circle ${\mathbb T}$).
Earlier, the subordination property for free additive convolutions [22]
was discovered by Voiculescu [24] for compactly supported measures
on ${\mathbb R}$, generalized by Biane [8] to arbitrary measures 
on ${\mathbb R}$ (see also [17] for a related approach).

For instance, subordination for free multiplicative convolutions of probability
measures on ${\mathbb R}_{+}$ says that for given 
$\mu_1, \mu_2\in {\mathcal M}_{{\mathbb R}_{+}}\setminus \{\delta_{0}\}$, there exist analytic self-maps
$\eta_{1},\eta_{2}$ of ${\mathbb C}\setminus {\mathbb R}_{+}$, such that
\begin{equation}\tag{1.3}
\eta_{\mu_1  \boxtimes \, \mu_2}(z)=\eta_{\mu_1}(\eta_{1}(z))=\eta_{\mu_2}(\eta_{2}(z)).
\end{equation}
One says that $\eta_{\mu_1\,\boxtimes\,\mu_2}$
is {\it subordinate} to both $\eta_{\mu_1}$ and $\eta_{\mu_2}$, with $\eta_1$ and $\eta_2$ 
being the so-called {\it subordination functions}. 
These functions play a key role in the
analytical study of free convolutions [3,4,10].

The functions $\eta_{1}$ and $\eta_{2}$ are unique and can be viewed as $\eta$-transforms 
of certain probability measures on ${\mathbb R}_{+}$ which are not 
concentrated at zero. This defines a binary operation $\boxslash$ on 
${\mathcal M}_{{\mathbb R}_{+}}\setminus \{\delta_{0}\}$, namely
\begin{equation}\tag{1.4}
\eta_{1}(z)=\eta_{\mu_2 \boxslash\, \mu_1}(z),\;\;\;  {\rm and} \;\;\;
\eta_{2}(z)=\eta_{\mu_1 \boxslash\, \mu_2}(z).
\end{equation}
The associated convolution $\mu_1 \boxslash \mu_2$, introduced in this paper, is 
called the {\it s-free multiplicative convolution} and it plays the role of the 
multiplicative analog of the s-free additive convolution 
$\mu_1 \boxright \mu_2$ studied in [16].

The subordination formulas (1.3) are related to the so-called 
{\it monotone multiplicative convolution} of 
probability measures on ${\mathbb R}_{+}$, introduced and 
studied by Bercovici [6].
This convolution can be defined by the equation
\begin{equation}\tag{1.5}
\eta_{\mu_1\circlearrowright \mu_2}(z)=\eta_{\mu_1}(\eta_{\mu_2}(z))
\end{equation}
for $z\in {\mathbb C}\setminus {\mathbb R}_{+}$,
where  $\mu_1, \mu_2\in {\mathcal M}_{{\mathbb R}_{+}}$.
Using this convolution and the s-free multiplicative convolution, 
we can write subordination equations (1.3) in the convolution form
\begin{equation}\tag{1.6}
\mu_1 \boxtimes \mu_2 = \mu_1 \circlearrowright (\mu_2 \boxslash \mu_1)=  
\mu_2 \circlearrowright (\mu_1 \boxslash \mu_2)
\end{equation}  
for $\mu_1, \mu_2\in {\mathcal M}_{{\mathbb R}_{+}}\setminus \{\delta_{0}\}$.

Further, we show that one can decompose both s-free and free multiplicative
convolutions of compactly supported measures on ${\mathbb R}_{+}$ which are
not concentrated at zero in terms of simpler convolutions.
For that purpose we introduce another convolution 
of $\mu_1,\mu_2\in {\mathcal M}_{{\mathbb R}_{+}}\setminus \{\delta_{0}\}$
called the {\it orthogonal multiplicative convolution},
denoted $\mu_1 \angle \mu_2$.
It plays the role of a multiplicative analog of the orthogonal additive 
convolution on ${\mathbb R}$, introduced and studied in [16].
Using transforms, one can define this convolution by 
\begin{equation}\tag{1.7}
\eta_{\mu_1 \angle \mu_2}(z)=\frac{z\eta_{\mu_1}(\eta_{\mu_2}(z))}{\eta_{\mu_2}(z)}
\end{equation}
for $z\in {\mathbb C}\setminus {\mathbb R}_{+}$.
Using this convolution, we obtain a decomposition of the s-free multiplicative
convolution of $\mu_1 \boxslash \mu_2$ in terms of an infinite sequence of 
alternating $\mu_1,\mu_2$ if these are compactly supported.
In view of (1.6), this leads to a decomposition of the free multiplicative convolution.

In a similar way one can treat probability measures on the unit circle
${\mathbb T}$, denoted ${\mathcal M}_{{\mathbb T}}$. In that case, $\eta_{\mu}$ 
is the integral over ${\mathbb T}$ and $z$ lies inside the open unit disc ${\mathbb D}$.
One defines subordination functions for $\mu_1,\mu_2\in  {\mathcal M}_{*}$,
where ${\mathcal M}_{*}={\mathcal M}_{{\mathbb T}}\setminus \{\mu:\mu(X)=0\}$.
Let us also observe that a number of results hold for
distributions $\mu_1, \mu_2\in \Sigma$, where $\Sigma$ denotes the set 
functionals $\mu: {\mathbb C}[X]\rightarrow {\mathbb C}$, which send $1$ into $1$.
In that case, the functions $\eta_{\mu}$ are understood as formal power series.

One of the main points of this paper is that our study of the relations between 
subordination functions and their decompositions uses operatorial techniques.
Namely, we construct bounded operators on Hilbert spaces, which correspond
to all compactly supported convolutions which appear in the subordination equations 
and in the decompositions of s-free convolutions. Our approach 
is based on the decomposition of the free product of Hilbert spaces 
$({\mathcal H},\xi)=({\mathcal H}_{1},\xi_{1})*({\mathcal H}_{2},\xi_{2})$, 
as (two different) orthogonal direct sums
\begin{equation}\tag{1.8}
{\mathcal H}=\bigoplus_{n=1}^{\infty}{\mathcal H}^{(n-1)}(\iota)\oplus {\mathcal H}^{(n)}(\bar{\iota})
\end{equation}
for each $\iota \in I=\{1,2\}$, where $\overline{1}=2 $ and $\overline{2}=1$,
with ${\mathcal H}^{(n)}(\iota)$ denoting the subspace spanned by 
alter\-nating tensor products of order $n$ which do not end with a vector 
from ${\mathcal H}_{\iota}^{0}:={\mathcal H}_{\iota}\ominus {\mathbb C}\xi_{\iota}$, following the 
original notation of Voiculescu [21] (see also [25]), and we set 
${\mathcal H}^{(0)}(\iota)={\mathbb C}\xi$ for each $\iota \in I$. By 
$P_{\iota}(n)$ we denote the canonical projection from
${\mathcal H}$ onto ${\mathcal H}^{(n-1)}(\iota)\oplus {\mathcal H}^{(n)}(\overline{\iota})$ and we set
$P_{\xi}$ to be the projection onto ${\mathbb C}\xi$.
Finally, we define the so-called {\it vacuum state} $\varphi(\cdot)=\langle \cdot\xi,\xi\rangle$
on $\mathcal{B}(\mathcal{H})$.

According to the above decomposition, we can represent bounded free random variables in 
the form of `orthogonal series'
\begin{equation}\tag{1.9}
X_{\iota}= \sum_{j=1}^{\infty}X_{\iota}(j),
\end{equation}
where $\iota\in I$ and $X_{\iota}(j)=P_{\iota}(j)X_{\iota}P_{\iota}(j)$.
These series play a crucial role in our study of subordination for 
both free additive convolutions [16] and free multiplicative 
convolutions, studied in this paper. 
Let us remark that the `orthogonal series' were introduced in a more general 
context of `monotone closed *-algebras' of operators `affiliated' with unital *-algebras [14]. 
In that approach, we used the `monotone tensor product' $\overline{\otimes}$ 
(reminding the von Neumann tensor product) to represent free random variables. 

Let  $\mu_1$, $\mu_2$ be $\varphi$--distributions of bounded 
random variables $X_1$ and $X_2$, respectively, which are free with respect to $\varphi$.
In the case of free additive convolutions, crucial was 
the decomposition of their sum in terms of `additive subordination branches', 
namely $X_1+X_2=B_{1}+B_2$, with
\begin{equation}\tag{1.10}
B_{\iota}=S_{\iota}^{\scriptscriptstyle {\rm odd}} + 
S_{\overline{\iota}}^{\scriptscriptstyle {\rm even}},
\end{equation}
where
\begin{equation}\tag{1.11}
S_{\iota}^{\scriptscriptstyle {\rm odd}}=
\sum_{\scriptscriptstyle j \;{\rm odd}}X_{\iota}(j)\;\;\; {\rm and}\;\;\;
S_{\iota}^{\scriptscriptstyle {\rm even}} =
\sum_{\scriptscriptstyle j\; {\rm even}} X_{\iota}(j),
\end{equation}
for each $\iota \in I$. Although the subordination branch $B_{\iota}$ is a bounded
operator on the `free Fock space' ${\mathcal H}$ for each $\iota \in I$, 
it acts non-trivially only on its subspace, the `s-free Fock space' ${\mathcal K}_{\iota}$, defined as
\begin{equation}\tag{1.12}
{\mathcal K}_{\iota}=
\bigoplus_{\scriptscriptstyle n\;{\rm odd}}{\mathcal H}^{(n-1)}(\iota)\oplus {\mathcal H}^{(n)}(\overline{\iota}),
\end{equation}
and therefore, it is appropriate to view $B_{\iota}$ as an element of $B({\mathcal K}_{\iota})$.

For operatorial subordination in the mutliplicative case, we need
to modify the operators $S_{\iota}^{ \scriptscriptstyle {\rm odd}}$, 
$S_{\iota}^{ \scriptscriptstyle {\rm even}}$ and take 
\begin{equation}\tag{1.13}
R_{\iota}^{ \scriptscriptstyle {\rm odd}}=S_{\iota}^{ \scriptscriptstyle {\rm odd}}\;\;\;
{\rm and}
\;\;\;
R_{\iota}^{ \scriptscriptstyle {\rm even}} =P_{\xi}+S_{\iota}^{ \scriptscriptstyle {\rm even}}.
\end{equation}
Moreover, in order to include all bounded positive random variables, we shall 
need the operation $\boxslash$ to be defined for all compactly supported 
$\mu_1, \mu_2\in {\mathcal M}_{{\mathbb R}_{+}}$. For that purpose we set 
\begin{equation}\tag{1.14}
\mu \boxslash \delta_{0}= \delta_{\mu(X)}\;\;\; {\rm and}\;\;\;\delta_{0}\boxslash \mu =\delta_{0}
\end{equation}
for compactly supported $\mu \in {\mathcal M}_{{\mathbb R}_{+}}$, which 
turns out natural in our operatorial setting.

Using these notations, and choosing, for simiplicity, 
an existence-type formulation, we can summarize our subordination result for
positive operators as follows. \\
\indent{\par}
{\sc Subordination for positive operators.}
{\it If $X_1$ and $X_2$ are positive, then there
exists a decomposition $\sqrt{X_1}X_2\sqrt{X_1}=\sqrt{y}\,\,Y \!\sqrt{y}$,
such that the $\varphi$--distributions of $y$ and $Y$ are $\mu_1$ and 
$\mu_2\boxslash \mu_1$, respectively, and the pair $(y-1,Y-1)$ is monotone 
independent w.r.t. $\varphi$. Moreover, the $\varphi$--distribution of 
$Y$ agrees with that of $\sqrt{R_1^{\scriptscriptstyle {\rm even}}}
R_2^{\scriptscriptstyle {\rm odd}}
\sqrt{R_1^{\scriptscriptstyle {\rm even}}}$, 
and the pair $(R_{2}^{\scriptscriptstyle {\rm odd}}-1_{{\mathcal K}_{2}},
R_{1}^{\scriptscriptstyle {\rm even}}-1_{{\mathcal K}_{2}})$
is s-free independent w.r.t. $(\varphi, \psi)$, where $\psi$ is 
the state associated with any unit vector $\zeta\in {\mathcal H}_{2}^{0}$.}\\
\indent{\par}
Moreover, we further decompose the `positive subordination branches'
\begin{equation}\tag{1.15}
\sqrt{R_1^{\scriptscriptstyle {\rm even}}}R_2^{\scriptscriptstyle {\rm odd}}
\sqrt{R_1^{\scriptscriptstyle {\rm even}}}
\end{equation}
using the notion of `orthogonal independence'.
These decompositions corresponds to decompositions of s-free convolutions 
in terms of orthogonal convolutions.
In a similar way we obtain operatorial subordination results
for unitary operators, or even more generally, for bounded operators.

The main examples of operators related to subordination for the free 
additive convolution, or, more generally, to various notions of independence 
${\mathcal I}$, are the adjacency matrices of subgraphs of the corresponding ${\mathcal I}$--product 
graphs ${\mathcal G}_{1}{\mathcal I}{\mathcal G}_{2}$ (in our study, ${\mathcal I}$ stands for orthogonal, comb, star, s-free or free). In fact, to each independence ${\mathcal I}$, 
we can associate an additive and a multiplicative convolution,
\begin{equation}\tag{1.16}
\mu_1 +_{\mathcal I} \mu_2 \;\;\;{\rm and}\;\;\;\mu_1 \times_{\mathcal I} \mu_2
\end{equation}
respectively. Recall that the additive ${\mathcal I}$--convolution
of spectral distributions of rooted graphs corresponds to the addition 
of `monochromatic' ${\mathcal I}$-independent adjacency matrices
$S_{\iota}$, $\iota \in I$, and that, in turn, is related to 
the enumeration of rooted (i.e. root-to-root) walks on 
${\mathcal G}_{1}{{\mathcal I}}{\mathcal G}_{2}$ (for details, see [2] and its references). 

In order to find a `universal' multiplicative analog of this theorem,
one needs to introduce a new concept of a product of 
${\mathcal G}_{1}$ and ${\mathcal G}_{2}$ for each ${\mathcal I}$--independence, 
which we call the ${\mathcal I}$ {\it loop product} and denote 
${\mathcal G}_{1}{\mathcal I}_{\ell} {\mathcal G}_{2}$. 
Roughly speaking, ${\mathcal G}_{1}{\mathcal I}_{\ell}{\mathcal G}_{2}$
is obtained by adding colored loops to 
${\mathcal G}_{1}{{\mathcal I}}{\mathcal G}_{2}$ in a suitable way. 
We assume that the product graph has a `natural coloring', by which 
we mean that each copy of ${\mathcal G}_{\iota}$
is colored with color $\iota$, where $\iota \in I$.
The procedure of adding loops is equivalent to replacing each $S_{\iota}$ 
by its `unitization' $R_{\iota}$ obtained from 
$S_{\iota}$ by adding some projection, in such a way that
makes the pair $(R_{1}-1, R_{2}-1)$ ${\mathcal I}$--independent. 

Quite naturally, in order to formulate our result, 
we shall use the formal power series $\eta_{Z}$ corresponding to
the variable $Z=R_{2}R_{1}$, namely
\begin{equation}\tag{1.17}
\eta_{Z}(z)=\sum_{n=1}^{\infty}N_{Z}(n)z^{n},
\end{equation}
and interpret $N_{Z}(n)$ as the `first return moments' in the state $\varphi_{e}$
associated with the root $e$. These moments are related to the enumeration of 
walks of the same class, which we find to be 
{\it rooted alternating d-walks} (`d-walk' is our abbreviation of 
`double return walks originating with color 1') 
counted on different products. Finally, in view of 
the relation to independence ${\mathcal I}$ just mentioned, these moments
agree with the $\eta$--moments of the corresponding multiplicative 
convolutions $\mu_1\times_{{\mathcal I}}\mu_2$. It is not hard to see that
this result can be generalized to the framework of random walks [26].\\
\indent{\par}
{\sc Multiplication theorem.}
{\it Let ${\mathcal G}_{1}{\mathcal I}_{\ell}{\mathcal G}_{2}$ be naturally colored
and let $A({\mathcal G}_{1}{\mathcal I}_{\ell}{\mathcal G}_{2})=R_1+R_2$ be the
decomposition of its adjacency matrix induced by the coloring. Then
\begin{equation}\tag{1.17}
N_{Z}(n)=N_{\mu_1\times_{{\mathcal I}}\mu_2}(n)=|D_{2n}(e)|
\end{equation}
where $Z=R_2R_1$ and $D_{2n}(e)$ denotes the set of rooted alternating d-walks 
on ${\mathcal G}_{1}{\mathcal I}_{\ell}{\mathcal G}_{2}$ 
of lenght $2n$, where $n\in {\mathbb N}$. Moreover, the pair $(R_{1}-1, R_{2}-1)$ is 
${\mathcal I}$-independent.}\\
\indent{\par}
The paper is organized as follows. 
In Section 2, we introduce basic notions, including
the s-free multiplicative convolution. 
In Section 3, we introduce and study the concepts
of comb-- and star loop products of graphs and find 
relations between rooted alternating d-walks on these graphs and
the monotone and boolean multiplicative convolutions, respectively.
In Section 4, we introduce and study the 
orthogonal multiplicative convolution.
In Section 5 we define ans study the corresponding notion of the orthogonal loop product of
rooted graphs.
In Section 6 we show, by means of analytical methods, 
that the definition of the orthogonal multiplicative convolution
can be extended to arbitrary measures on ${\mathbb R}_{+}$ and ${\mathbb T}$.
The main operatorial results of the paper are contained in Sections 7 and 8, 
where we introduce and study operators on the free and s-free Fock spaces 
which are related to subordination for multiplicative 
free convolutions and their decompositions.  Finally, in Section 9 we 
find a relation between free and s-free multiplicative convolutions and the
enumeration of rooted alternating d-walks on the free product of graphs 
and on the s-free loop product of graphs, respectively.

Throughout the whole paper we understand that $I=\{1,2\}$ and we adopt the
notation $\overline{1}=2$ and $\overline{2}=1$. 
\\[10pt]
\section{Preliminaries}
By a {\it non-commutative probability space} we understand a pair $({\mathcal A}, \varphi)$,
where ${\mathcal A}$ is a unital algebra over ${\mathbb C}$ and $\varphi$ is
a linear functional $\varphi: {\mathcal A}\rightarrow {\mathbb C}$ such that
$\varphi(1)=1$. If ${\mathcal A}$ is a unital *-algebra and $\varphi$ is positive (called a {\it state}),
then $({\mathcal A},\varphi)$ is called a *-{\it probability space}. If, in addition,
${\mathcal A}$ is a $C^{*}$-algebra, then $({\mathcal A}, \varphi)$ is called a
$C^{*}$-{\it probability space}. By the Gelfand-Naimark-Segal theorem, a $C^{*}$-probability
space can always be realized as a subalgebra of bounded operators on a Hilbert space
${\mathcal H}$ with a distinguished unit vector $\xi$, for which $\varphi(a)=\langle a\xi, \xi \rangle$
for $a\in {\mathcal A}$. 

By a {\it random variable} we will understand any element $a$ of the considered algebra ${\mathcal A}$.
If ${\mathcal A}$ is equipped with an involution, then a random variable $a$ will be called
{\it self-adjoint} if $a=a^{*}$. The $\varphi$-{\it distribution} of a random variable $a$ is
the functional $\mu_{a}:{\mathbb C}[X]\rightarrow {\mathbb C}$
given by $\mu_{a}(1)=1$, $\mu_{a}(X^{n})=\varphi(a^{n})$. In particular,
if $({\mathcal A}, \varphi)$ is a $C^{*}$-probability space, then
the distribution $\mu_{a}$ of a self-adjoint random variable $a\in {\mathcal A}$ extends to
a compactly supported probability measure $\mu$ on the real line.
In that case we will often use the same notation $\mu$ for both the
distribution of $a$ and the associated compactly supported probability measure.

Various notions of `independence' ${\mathcal I}$ lead to several types of convolutions of 
distributions (probability measures). 
If we have two random variables, $X_{1}\in {\mathcal A}_{1}$ and  $X_{2}\in {\mathcal A}_{2}$
with $\varphi$-distributions $\mu_1$ and $\mu_2$, respectively, , 
where ${\mathcal A}_{1}$ and ${\mathcal A}_{2}$ are ${\mathcal I}$-independent subalgebras of a noncommutative
probability space $({\mathcal A}, \varphi)$, then the additive convolution of $\mu$ and $\nu$ 
associated with ${\mathcal I}$-independence is the $\varphi$-distribution of the sum $X_{1}+X_{2}$. In turn, the 
multiplicative convolution of $\mu$ and $\nu$ is the $\varphi$-distribution of $X_{2}X_{1}$. 
In this paper we are interested in the free and s-free multiplicative convolutions
associated with free and s-free independence, respectively.
In order to decompose them, we shall use the monotone multiplicative convolution introduced by 
Bercovici [6] and we will introduce the orthogonal multiplicative convolution associated with the notion
of `orthogonal independence' [16].
\begin{Definition}
{\rm Let $\mu_{1},\mu_{2}\in {\mathcal M}_{{\mathbb R}_{+}}\setminus \{\delta_{0}\}$ 
and let $\eta_{1}$ and $\eta_{2}$ be the associated subordination functions.
The unique probability measures $\sigma_{1},\sigma_{2}\in {\mathcal M}_{{\mathbb R}_{+}}\setminus
\{\delta_{0}\}$ such that $\eta_{1}=\eta_{\sigma_{2}}$ and $\eta_{2}=\eta_{\sigma_{1}}$ will be called 
the {\it s-free multiplicative convolution}
of $\mu_{1},\mu_{2}$ and $\mu_2, \mu_1$, respectively.
We set $\sigma_{1}=\mu_{1}\boxslash \mu_{2}$ and $\sigma_{2}=\mu_{2}\boxslash \mu_1$.
In a similar way we define the s-free multiplicative 
convolution of $\mu_1,\mu_2\in {\mathcal M}_{*}$.}
\end{Definition}

The s-free multiplicative convolution defines a binary operation $\,\boxslash\,$ on 
both ${\mathcal M}_{{\mathbb R}_{+}}\setminus \{\delta_{0}\}$ and ${\mathcal M}_{*}$.
It can be seen that it is neither commutative nor associative. 
We will show later that it is related to the concept of `freeness with subordination', or `s-free independence',
introduced in [16] (a relation between freeness and monotone independence
was also studied in [15]).
Moreover, it allows us to rewrite (1.3) in terms of convolutions. 
\begin{Proposition}
If $\mu_1,\mu_2\in {\mathcal M}_{{\mathbb R}_{+}}\setminus \{\delta_{0}\}$, then
\begin{equation}\tag{2.1}
\mu_1 \boxtimes \mu_2 = \mu_1 \circlearrowright (\mu_2 \boxslash \mu_1) =
\mu_2 \circlearrowright (\mu_1 \boxslash \mu_2)
\end{equation}
is the decomposition of the free multiplicative convolution 
corresponding to the subordination equations (1.3).
A similar result holds for $\mu_1, \mu_2\in {\mathcal M}_{*}$.
\end{Proposition}
{\it Proof.}
This fact is an immediate consequence of (1.3), in view of (1.5).
\hfill $\blacksquare$\\
\indent{\par}
By a {\it rooted set} we understand a pair $(X,e)$, where $X$ is a countable
set and $e$ is a distinguished element of $X$ called {\it root}.
By a {\it rooted graph} we understand a pair $({\mathcal G},e)$, where
${\mathcal G}=(V,E)$ is a non-oriented graph with the set of {\it vertices} $V=V({\mathcal G})$,
and the set of edges $E=E({\mathcal G})\subseteq \{(x,x'): \; x,x'\in V\}$
and $e\in V$ is a distinguished vertex called the {\it root}. 
We identify $(x,x')$ with $(x',x)$ since we consider non-oriented graphs, but 
when speaking of walks we find it convenient to say that $(x,x')$ {\it begins} with $x$ and {\it terminates} with $x'$.
Another reason for using this notation for edges is that 
of main interest to us are graphs which have {\it loops}, i.e. edges of the form $(x,x)$, where $x\in V$.
The notion of a rooted graph can be easily generalized to allow multiple edges.
Formally, we then obtain a {\it rooted multigraph}, but we will still use the term `rooted graph', 
or simply `graph' since all graphs will be considered to have a root.
Moreover, we will very often omit the root in our notation and denote by ${\mathcal G}$ the rooted graph 
$({\mathcal G},e)$ if no confusion arises.
Thus, in a graph ${\mathcal G}$ we denote by $n(x,x')$ the number of edges connecting $x$ and $x'$ (it may
be zero).

For (rooted) graphs we will also use the notation
\begin{equation}\tag{2.2}
V^{0}=V\setminus \{e\}.
\end{equation}
Two vertices $x,x'\in V$ are called {\it adjacent} if $x$ and $x'$ are connected by an edge,
which we denote $x \sim x'$. The {\it degree} of $x\in V$ is defined by $\kappa(x)=\sum_{x'\sim x}n(x,x')$.
A graph is called {\it locally finite} if $\kappa(x)<\infty$ for every $x\in V$.
It is called {\it uniformly locally finite} if ${\rm sup}\{\kappa(x): x\in V\}<\infty$.
All graphs considered in this paper will be connected and uniformly locally finite.

For $x\in V$, let $\delta(x)$ be the indicator function of the one-element set $\{x\}$.
Then $\{\delta(x),\,x\in V\}$ is an orthonormal basis of the Hilbert space $l^{2}(V)$ of
square integrable functions on the set $V$, with the usual inner product.
The {\it adjacency matrix} $A=A({\mathcal G})$ of a graph ${\mathcal G}$ is a 
matrix defined by $A({\mathcal G})_{x,x'}=n(x,x')$
We identify $A$ with the densely defined symmetric operator on $l^{2}(V)$ defined by
\begin{equation}\tag{2.3}
A({\mathcal G})\delta(x)=\sum_{x\sim x'}n(x,x')\delta(x')
\end{equation}
for $x\in V$. Notice that the sum on the right-hand-side is finite since our graph is assumed
to be locally finite. It is known that $A({\mathcal G})$ is bounded iff ${\mathcal G}$ is uniformly
locally finite. 

By the {\it spectral distribution of $A({\mathcal G})$ in a state} $\psi$ on $l_{2}(V)$
we understand the measure $\mu$ for which
\begin{equation}\tag{2.4}
\psi (A^{n})= \int_{{\mathbb R}}x^{n}\mu(dx), \;\; n\in {\mathbb N} \cup \{0\},
\end{equation}
and by the {\it spectral distribution of the rooted graph $({\mathcal G},e)$}
we understand the spectral distribution of $A({\mathcal G})$ in the state
$\varphi _{e}(.)=\langle . \delta(e),\delta(e)\rangle$.

A {\it walk} from $v_{0}$ to $v_{n}$ on a graph ${\mathcal G}$ is an alternating sequence 
of vertices and edges of the form
\begin{equation}\tag{2.5}
w=(v_{0}, \beta_{1}, v_{1}, \beta_{2}, v_{2}, \ldots , v_{n-1}, \beta_{n}, v_{n}),
\end{equation}
with vertices $v_{0}, v_{1}, \ldots , v_{n}\in V({\mathcal G})$ 
and edges $\beta_{1}, \beta_{2}, \ldots , \beta_{n}\in E({\mathcal G})$, such that
$\beta_{i}$ is an edge connecting $v_{i-1}$ and $v_{i}$. The lenght of $w$ is given by $l(w)=n$. 
We allow $v_{i-1}=v_{i}$, in particular this happens if $\beta_{i}$ is a loop. 

A {\it subwalk} of $w$ is a subsequence of $w$ of the form
\begin{equation}\tag{2.6}
w'=(v_{i}, \beta_{i}, v_{i+1}, \beta_{i+1}, \ldots , \beta_{j-1}, v_{j}),
\end{equation}
where $0\leq i <j \leq n$. 
In the case when ${\mathcal G}$ does not have 
multiple edges, we can identify $w$ with the sequence $(v_{0},v_{1},\ldots ,v_{n})$ with the understanding that
$v_{i-1}\sim v_{i}$ (in this case there is no confusion which edge connecting $v_{i-1}$ and $v_{i}$ is chosen)
for $1\leq i\leq n$. 

The set of walks from $v$ to $v'$ will be denoted $W(v,v')$ and we set $W(v,v)=W(v)$. 
A walk from $v=v_{0}$ to $v=v_{n}$, in which $v_{i}\neq v$ for all $1 \leq i \leq n-1$ 
will be called an {\it f-walk}. Note that we do not require all vertices 
$v_{0},v_{1}, \ldots , v_{n-1}$ in an f-walk to be distinct.
The set of f-walks from $v$ to $v$ will be denoted $F(v)$.
Subsets of $W(v)$ and $F(v)$ consisting of walks and f-walks of length $n$, respectively, 
will be denoted $W_{n}(v)$ and $F_{n}(v)$, respectively.  

A walk $w$ can also be represented as a sequence of subwalks; in particular, 
an alternating sequence of edges and subwalks of the form 
\begin{equation}\tag{2.7}
w=(w_{0}, \beta_{1}, w_{1}, \beta_{2}, w_{2}, \ldots , w_{r-1}, \beta_{r}, w_{r}),
\end{equation}
where $w_{0}$ is a walk from $v_{0}$ to some $v_{i(1)}$, $\beta_{1}$ is an edge connecting $v_{i(1)}$ with $v_{i(1)+1}$,
$w_{1}$ is a walk from $v_{i(1)+1}$ to some $v_{i(2)}$, $\beta_{2}$ is an edge connecting $v_{i(2)}$ with $v_{i(2)+1}$, 
and so on, finally, 
$w_{r}$ is a walk from $v_{i(r)+1}$ to $v_{i(r+1)}=v_{n}$. We shall also use similar representations
of walks which begin or end with an edge. Representing a walk in terms of edges and
subwalks is very convenient when the subwalks $w_{1},w_{2}, \ldots, w_{r}$ are of special type 
(for instance, are f-walks). 

A graph ${\mathcal G}$ whose edges are colored with colors 
from the set $I$ will be called an 
$I$-{\it edge-colored} graph.
In our case, we will always assume that $I=\{1,2\}$, so 
we can also use the term `$2$-edge-colored graph'. A walk in 
such a graph will be called ${\it alternating}$ if the colors of
its edges alternate.
\begin{Definition}
{\rm Let ${\mathcal G}$ be $I$-colored. A walk $w\in W(v)$ on ${\mathcal G}$ 
will be called an {\it alternating double return walk},
or simply an {\it alternating d-walk}, if it is alternating, 
begins with an edge of color $1$, and can be represented as a pair $(u_1,u_2)$ of 
subwalks, such that $u_1,u_2\in F(v)$. The set of alternating d-walks from 
$v$ to $v$ will be denoted $D(v)$, and we set $D_{m}(v)=W_{m}(v)\cap D(v)$.}
\end{Definition}

Edge-coloring by a two-element set is natural in the case of many 
products of rooted graphs. In particular, the products of rooted graphs of type
${\mathcal G}_{1}{\mathcal I}{\mathcal G}_{2}$ have the property that every edge of the product graph 
belongs either to a copy of ${\mathcal G}_{1}$, or to a copy of ${\mathcal G}_{2}$.
Thus we can color the edges of all copies of ${\mathcal G}_{\iota}$
with color $\iota$, where $\iota \in I$, in which case we will say that 
${\mathcal G}_{1}{\mathcal I}{\mathcal G}_{2}$ is {\it naturally colored}.
In turn, we will say that the product 
${\mathcal G}_{1}{\mathcal I}_{\ell}{\mathcal G}_{2}$ is {\it naturally colored} if 
its coloring is inherited from ${\mathcal G}_{1}{\mathcal I}{\mathcal G}_{2}$.
Let us also remark that our condition that an alternating walk should begin with 
an edge of color $1$ (and thus end with an adge of color $2$) 
is caused by the fact that we want to identify the given walk $w$ and 
its `reverse' obtained by reversing the order in (2.12).

We end this Section with elementary formulas for the `first return moments' of a random variable $X$, 
by which we understand the coefficients of the formal power series $\eta_{X}$ associated with
the distribution of $X$. If $X$ is a random variable in a non-commutative probability space
$({\mathcal A}, \varphi)$, then the `moment generating function' and the `first return moment generating function', respectively, associated with the $\varphi$-distribution $\mu$ of $X$, 
are given by formal power series
\begin{equation}\tag{2.8}
\psi_{X}(z)=\sum_{n=1}^{\infty}\mu(X^{n})(n)z^{n}, \;\;\; \eta_{X}(z)=\frac{\psi_{X}(z)}{1+\psi_{X}(z)}=\sum_{n=1}^{\infty}N_{X}(n)z^{n},
\end{equation}
where the numbers $N_{X}(n)$, $n\in {\mathbb N}$, 
will sometimes be called `first return moments of $X$'.

Below we give a convenient algebraic formula for these `first return moments'.
For that purpose, let us take the extension of $({\mathcal A},\varphi)$ to a larger 
noncommutative probability space, namely, the free product with identified units 
${\mathcal A}*{\mathbb C}[P]$, where $P^2=P$, with the state 
given by the linear extension of
\begin{equation}\tag{2.9}
\varphi(P)=1, \;\;\;\varphi(w_{1}Pw_{2})=\varphi(w_{1})\varphi(w_{2})
\end{equation}
for any $w_1, w_2\in {\mathcal A}*{\mathbb C}[P]$, where, slightly abusing notation,
we also denote the new state by $\varphi$ [13]. In particular, in the $C^{*}$-algebra context,
$P$ can be interpreted as the projection onto ${\mathbb C}\xi$, where
$\xi$ is the cyclic unit vector of the GNS triple associated with 
$({\mathcal A}, \varphi)$.
\begin{Proposition}
Let $X$ be a random variable in a noncommutative probability space
$({\mathcal A},\varphi)$. Then
\begin{equation}\tag{2.10}
N_{X}(n)  =\varphi (X(P^{\perp}X)^{n-1})
\end{equation}
{\it for $n\in {\mathbb N}$, where $P^{\perp}=1-P$.}
\end{Proposition}
{\it Proof.}
This is a straightforward consequence of (1.2). \hfill $\blacksquare$\\
\indent{\par}
Finally, let us prove an elementary fact about a relation
between the `first return moments' of products of `monochromatic' 
adjacency matrices in certain $I$-colored graphs and the
cardinalities of the sets $D_{2n}(e)$.
\begin{Proposition}
{\it Let $({\mathcal G},e)$ be an $I$-colored uniformly locally finite graph, 
and let $A({\mathcal G})=A_1+A_2$ be the decomposition of $A({\mathcal G})$ 
induced by the coloring. If the set of rooted alternating f-walks of even lenghts
is empty, then
\begin{equation}\tag{2.11}
N_{Z}(n)=|D_{2n}(e)|,
\end{equation}
where $Z=A_2A_1$ and the numbers $N_{Z}(n)$ correspond to the state 
$\varphi_{e}(.)=\langle .\delta(e), \delta(e)\rangle$.}
\end{Proposition}
{\it Proof.}
If we take $Z=A_2A_1$ in Proposition 2.2, then
$N_{Z}(n)$ is equal to the number of rooted 
alternating walks of lenght $2n$, which begin 
with an edge of color $1$ and such that intermediate returns to the root 
can occur only after odd numbers of steps. 
If there are no rooted alternating f-walks of even lenght, then
the first return to the root must occur after an odd number of steps
and thus the second return to the root must occur after an even number
of steps. Therefore, if we consider rooted alternating 
walks which have an even number of steps and such that intermediate
returns occur only after an odd number of steps, these
are d-walks. This proves our assertion.
\hfill $\blacksquare$\\
\indent{\par}
Let us note that all graph products considered in this paper 
satisfy the assumptions of Proposition 2.3, and therefore, 
walk-counting is always reduced to rooted alternating d-walks.\\[10pt]

\section{Comb and star loop products of graphs}
We begin with recalling the notions of the additive [18] and multiplicative [6] 
monotone convolutions and show that
the $\eta$-moments of the latter are related to the enumeration of alternating d-walks
on a new version of the comb product of graphs called the `comb loop product'.
\begin{Definition}
{\rm Two subalgebras ${\mathcal A}_{1}, {\mathcal A}_{2}$ of a unital algebra ${\mathcal A}$
are {\it monotone independent} with respect to a normalized linear functional $\varphi$ 
on ${\mathcal A}$ if \\
\indent{\par}
$\varphi(w_{1}a_{1}b)=\varphi(b)\varphi(w_{1}a_{1})$\;\;{\rm and}\;\;
$\varphi(ba_{2}w_{2})=\varphi(b)\varphi(a_{2}w_{2})$,
\indent{\par}
$\varphi(w_{1}a_{1}ba_{2}w_{2})=
\varphi(b)\varphi(w_{1}a_{1}a_{2}w_{2})$,\\[5pt]
whenever $a_{1},a_{2}\in {\mathcal A}_{1}$, $b\in {\mathcal A}_{2}$ and
$w_{1},w_{2}$ are arbitrary elements of the unital algebra ${\rm alg}({\mathcal A}_{1},{\mathcal A}_{2})$ 
generated by ${\mathcal A}_{1}$ and ${\mathcal A}_{2}$. In particular, we will say that the pair $(a,b)$ of elements of
${\mathcal A}$ is monotone independent w.r.t. $\varphi$ if the algebras generated by 
these elements are monotone independent.} 
\end{Definition}

Suppose $\varphi$-distributions of random variables $a_1$ and $a_2$ are $\mu_1$ and $\mu_2$, respectively.
If the pair $(a_1,a_2)$ is monotone independent w.r.t. $\varphi$, then
the $\varphi$-distribution of the sum $a_1+a_2$ is the monotone additive convolution 
$\mu_1\vartriangleright\mu_2$. In turn, if the pair $(a_1-1,a_2-1)$ is monotone independent w.r.t. $\varphi$, 
then the $\varphi$-distribution of the product $a_2a_1$ is the {\it monotone multiplicative convolution} 
$\mu_1\circlearrowright\mu_2$. The corresponding formal power series
satisfy (1.5). If $\mu_1,\mu_2\in {\mathcal M}_{{\mathbb R}_{+}}$ or $\mu_1,\mu_2\in {\mathcal M}_{{\mathbb T}}$, then (1.5) is understood in terms of analytic self-maps
of ${\mathbb C}\setminus {\mathbb R}_{+}$ or ${\mathbb T}$,
respectively. For details, see [19] (additive case) and [6] (multiplicative case).

It is known that the monotone independence can be associated with the comb product of graphs [1]. Let us recall the definition of the comb product and follow up
with the corresponding loop product. 
\begin{Definition}
{\rm The {\it comb product} of rooted graphs $({\mathcal G}_{1},e_1)$ and $({\mathcal G}_{2},e_2)$
is the rooted graph $({\mathcal G}_{1}\vartriangleright {\mathcal G}_{2},e)$
obtained by attaching a copy of ${\mathcal G}_{2}$ by its root $e_{2}$
to each vertex of ${\mathcal G}_{1}$, where we denote by $e$ the vertex obtained
by identifying $e_1$ and $e_2$. If no confusion arises, we denote the comb product by
${\mathcal G}_{1}\vartriangleright {\mathcal G}_{2}$. If we identify its set of vertices 
with $V_{1}\times V_{2}$, then its root is identified with $e_{1}\times e_{2}$.} 
\end{Definition}
\begin{Definition}
{\rm Suppose the edges of the comb product $({\mathcal G}_{1}\vartriangleright {\mathcal G}_{2},e)$
are naturally colored. The {\it comb loop product} of rooted graphs 
$({\mathcal G}_{1},e_1)$ and $({\mathcal G}_{2},e_2)$
is the rooted graph $({\mathcal G}_{1}\vartriangleright_{\ell}{\mathcal G}_{2},e)$
obtained from $({\mathcal G}_{1}\vartriangleright {\mathcal G}_{2},e)$
by attaching a loop of color $1$ to each vertex but the root of each
copy of ${\mathcal G}_{2}$.}
\end{Definition}

Let us justify the above definition. We follow the observation made by Bercovici [6] that
in order to introduce a multiplicative convolution associated with monotone independence, 
one needs to `unitize' the usual monotone independent variables. It has been shown in [5]
that one of the possible choices is to take variables of the form
\begin{equation}\tag{3.1}
R_1=a_1\otimes P_2+ 1_1\otimes P_2^{\perp} \;\;\;{\rm and}\;\;\; R_2=1_1\otimes a_2,
\end{equation}
where $a_{\iota},1_{\iota}\in B({\mathcal H}_{\iota})$, with $({\mathcal H}_{\iota}, \xi_{\iota})$,
being Hilbert spaces with distinguished unit vectors and identity operators $1_{\iota}$, 
and $\iota \in I$. Then the moments
of the product $R_{2}R_{1}$ in the state associated with the vector $\xi_1\otimes \xi_2$ agree with
the moments of $\mu_1\circlearrowright \mu_2$. 

If $a_{\iota}$ is taken to be
the adjacency matrix of a (uniformaly locally finite) graph 
${\mathcal G}_{\iota}$, $\iota\in I$, then the term $1_1\otimes P_{2}^{\perp}$ corresponds to
glueing a loop of color $1$ to each vertex but the root of each copy of ${\mathcal G}_{2}$. 
Thus, $R_1$ is obtained from the usual comb product adjacency matrix $S_{1}$ by adding a projection 
$L_{1}$, wheras $R_{2}=S_{2}$, where
\begin{equation}\tag{3.2}
S_1=a_1\otimes P_2, \;\; L_1=1_1\otimes P_2^{\perp}
\end{equation}
and $R_1$, $R_2$ are the adjacency matrices of `monochromatic', $1$ and $2$-colored subgraphs, 
respectively. This realization can be generalized to arbitrary random variables in noncommutative
probability spaces, using the extensions of functionals given by (2.9). 
\begin{Remark}
{\rm It is not difficult to show that the comb loop product of graphs ${\mathcal G}_{1}$ and
${\mathcal G}_{2}$ is isomorphic to the usual comb product of larger graphs 
$\widetilde{{\mathcal G}}_{1}$ and $\widetilde{{\mathcal G}}_{2}$.
However, 
in order to keep the same coloring of the product graph 
(wich is needed for counting alternating walks) 
one has to use two colors for $\widetilde{\mathcal G}_{2}$. 
Namely, $\widetilde{{\mathcal G}}_{2}$ is obtained from ${\mathcal G}_{2}$ by glueing a 
loop of color $1$ to each vertex but the root of ${\mathcal G}_{2}$. 
Then we indeed have 
${\mathcal G}_{1}\vartriangleright_{\ell} {\mathcal G}_{2}=
\widetilde{\mathcal G}_{1}\vartriangleright \widetilde{\mathcal G}_{2}$,
but usefulness of this relation seems to be limited (roughly speaking, 
we have a simpler product but a more complicated coloring).}
\end{Remark}
\begin{Example}
{\rm Let us consider the example of the comb loop product of graphs given in Figure 1. 
We choose both ${\mathcal G}_{1}$ and ${\mathcal G}_{2}$ to have a loop at the root
(to distinguish the loops, we draw loops of color 1 smaller than loops of color 2,
which is helpful in the enumeration of alternating walks).
Then ${\mathcal G}_{1}\vartriangleright_{\ell} {\mathcal G}_{2}$
has loops of two types: loops whose origin can be traced back to graphs ${\mathcal G}_{1}$ 
and ${\mathcal G}_{2}$ (these are all loops which are at the glueing points), 
or are added to the usual comb product in the process of forming the comb loop product.
\begin{figure}
\unitlength=1mm
\special{em.linewidth 2pt}
\linethickness{0.5pt}
\begin{picture}(60.00,30.00)(20.00,0.00)
\put(0.00,10.00){\line(1,0){20.00}}
\put(35.00,10.00){\line(0,1){6.00}}

\put(0.00,10.00){\circle*{1.60}}
\put(0.00,8.00){\circle{4.00}}
\put(10.00,10.00){\circle*{1.00}}
\put(20.00,10.00){\circle*{1.00}}

\put(35.00,10.00){\circle*{1.60}}
\put(35.00,6.50){\circle{7.00}}
\put(35.00,16.00){\circle*{1.00}}
\put(30.00,20.00){\circle*{1.00}}
\put(40.00,20.00){\circle*{1.00}}

\put(65.00,22.00){\circle{4.00}}
\put(70.00,18.00){\circle{4.00}}
\put(75.00,22.00){\circle{4.00}}

\put(80.00,22.00){\circle{4.00}}
\put(85.00,18.00){\circle{4.00}}
\put(90.00,22.00){\circle{4.00}}

\put(95.00,22.00){\circle{4.00}}
\put(100.00,18.00){\circle{4.00}}
\put(105.00,22.00){\circle{4.00}}

\put(1.00,11.00){\scriptsize $e_{1}$}
\put(71.00,11.00){\scriptsize $e$}
\put(10.00,11.00){\scriptsize $x$}
\put(20.00,11.00){\scriptsize $x'$}
\put(37.00,10.00){\scriptsize $e_{2}$}
\put(36.00,13.50){\scriptsize $y$}
\put(39.00,22.00){\scriptsize $y''$}
\put(29.00,22.00){\scriptsize $y'$}
\put(9.00,15.00){\scriptsize ${\mathcal G}_{1}$}
\put(26.00,15.00){\scriptsize ${\mathcal G}_{2}$}
\put(52.00,15.00){\scriptsize ${\mathcal G}_{1}\vartriangleright_{\ell} {\mathcal G}_{2}$}





\put(70.00,10.00){\line(1,0){30.00}}
\put(70.00,10.00){\line(0,1){6.00}}
\put(85.00,10.00){\line(0,1){6.00}}
\put(100.00,10.00){\line(0,1){6.00}}

\put(70.00,8.00){\circle{4.00}}
\put(70.00,6.50){\circle{7.00}}
\put(85.00,6.50){\circle{7.00}}
\put(100.00,6.50){\circle{7.00}}

\put(70.00,16.00){\circle*{1.00}}
\put(75.00,20.00){\circle*{1.00}}
\put(65.00,20.00){\circle*{1.00}}

\qbezier(65,20)(65,16)(70,16)
\qbezier(70,16)(75,16)(75,20)

\qbezier(80,20)(80,16)(85,16)
\qbezier(85,16)(90,16)(90,20)
\qbezier(95,20)(95,16)(100,16)
\qbezier(100,16)(105,16)(105,20)

\qbezier(30,20)(30,16)(35,16)
\qbezier(35,16)(40,16)(40,20)

\put(70.00,10.00){\circle*{1.60}}
\put(85.00,10.00){\circle*{1.00}}
\put(100.00,10.00){\circle*{1.00}}

\put(85.00,16.00){\circle*{1.00}}
\put(100.00,16.00){\circle*{1.00}}
\put(90.00,20.00){\circle*{1.00}}
\put(80.00,20.00){\circle*{1.00}}
\put(105.00,20.00){\circle*{1.00}}
\put(95.00,20.00){\circle*{1.00}}

\end{picture}
\caption{An example of the comb loop product}
\end{figure}
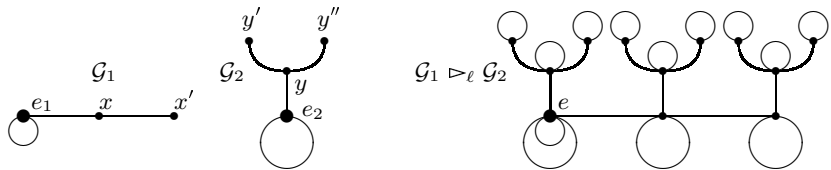
Denote by $(v,v)_{\iota}$ the loop from $v$ to $v$ of color $\iota$. The following 
walk is a rooted alternating d-walk:
$$
(\beta_{1},w_{1},\beta_{2}, w_{2},\beta_{3}, w_{3},\beta_{4}, w_{4}) \in D_{8}(e),
$$
where the sequence of edges of color 1, 
$$
(\beta_{1}, \beta_{2}, \beta_{3}, \beta_{4})=((e,x),(x,x'),(x',x),(x,e))
$$
forms a rooted f-walk of color $1$, which is interlaced with the sequence $(w_{1},w_{2},w_{3},w_{4})$
of loops of color 2:
$w_{1}=w_{3}=(x,x)_{2}\in F_{1}(x)$, $w_{2}=(x',x')_{2}\in F_{1}(x')$ and 
$w_{4}=(e,e)_{2}\in F_{1}(e)$.
Note that the first return to $e$ occurs at the end of $\beta_4$, whereas
the second return occurs as the end of $w_{4}$.
When any of the loops, $w_{i}$, $1\leq i \leq 4$, 
is replaced by an arbitrary alternating f-walk from the corresponding vertex to itself
which begins and ends with an edge of color $2$,
we still get a rooted alternating d-walk.
The simplest examples of a rooted alternating d-walk are of course given by:  
$(\beta_{0},w_{0})\in D_{2}(e)$,
where $\beta_{0}=(e,e)_{1}$, $w_{0}=(e,e)_{2}$ and $(\beta_{1},w_{1},\beta_{4}, w_{4}) \in D_{4}(e)$.}
\end{Example}
\begin{Theorem}
The Multiplication Theorem holds for monotone independence, the associated
loop product ${\mathcal G}_{1}\vartriangleright_{\ell} {\mathcal G}_{2}$, and the multiplicative convolution 
$\mu_1\circlearrowright \mu_2$.
\end{Theorem}
{\it Proof.}
From Definition 3.3 and (3.1) it follows that 
$$
Z=A({\mathcal G}_{1}\vartriangleright_{\ell} {\mathcal G}_{2})=R_1+R_2,
$$
where the pair $(R_{1}-1,R_{2}-1)$ is monotone independent w.r.t. $\varphi_{e}$, the state
associated with $e$ and $R_{\iota}$, $\iota \in I$, is the adjacency matrix
of the $\iota$--colored subgraph. Moroever, by Proposition 2.2, 
$N_{Z}(n)$ is the number of alternating d-walks on ${\mathcal G}_{1}\vartriangleright_{\ell} {\mathcal G}_{2}$
of lenght $2n$. This proves the first equation in (1.13) for monotone independence.
Now, from equation (1.5), we obtain the combinatorial formula 
$$
N_{\mu_1 \circlearrowright \mu_2}(n)=
\sum_{r=1}^{n}N_{\mu_1}(r)\sum _{k_{1}+k_{2}+\ldots +k_{r}=n}
N_{\mu_2}(k_{1})N_{\mu_2}(k_{2})\ldots
N_{\mu_2}(k_{r}),
$$
where it is tacitly assumed that the indices $k_{1}, k_{2}, \ldots , k_{r}$ are positive integers
and numbers $N_{\mu}(n)$ are coefficients of $\eta_{\mu}$ treated as formal power series 
$$
\eta_{\mu}(z)=\sum_{n=1}^{\infty}N_{\mu}(n)z^{n}.
$$
The above combinatorial formula allows us to find a correspondence between these $\eta$--moments and 
rooted alternating d-walks on ${\mathcal G}_{1}\vartriangleright_{\ell}{\mathcal G}_{2}$. 
Let us observe that $N_{\mu_1}(r)$ is the number of rooted 
f-walks of lenght $r$ on ${\mathcal G}_{1}$ and $N_{\mu_2}(k)$ is the number of
rooted f-walks of lenght $k$ on ${\mathcal G}_{2}$. 
On the other hand, recall that in the comb loop product
${\mathcal G}_{1}\vartriangleright_{\ell} {\mathcal G}_{2}$ there is only 
one copy of ${\mathcal G}_{1}$ (with $e_{1}$ identified with the root $e$ of the product graph), 
with a copy of ${\mathcal G}_{2}$ attached by its root to every vertex $x$ of ${\mathcal G}_{1}$
(the vertex $x$ becomes then the only common vertex of this copy of ${\mathcal G}_{2}$
and the original copy of ${\mathcal G}_{1}$). In addition, loops of color $1$ are glued to all
vertices but the roots of these copies of ${\mathcal G}_{2}$.
Therefore, each rooted alternating d-walk on ${\mathcal G}_{1}\vartriangleright_{\ell} {\mathcal G}_{2}$ 
consists of a sequence of edges of ${\mathcal G}_{1}$ which forms an f-walk of color $1$, 
namely $c=(v_{0},v_{1}, v_{2}, \ldots , v_{r-1},v_{r})\in F(e)$, and
alternating f-walks $c_{i}\in F(v_{i})$, $i=1, \ldots , r$, with the first edge
of color $2$, attached to every vertex of $c$.
Note that the number of alternating f-walks from $v_{i}$ to $v_{i}$ with the first edge
of color $2$ is equal to the number of walks from $e_2$ to $e_2$ in graph ${\mathcal G}_{2}$ --
the only difference is that if the latter has lenght $k_{i}$, the fomer has
lenght $2k_{i}-1$ since each edge must be followed by a loop of color $2$.
Thus, the contribution from each product of the form
$$
N_{\mu_1}(r)N_{\mu_2}(k_{1})N_{\mu_2}(k_{2})\ldots
N_{\mu_2}(k_{r})
$$ 
to the RHS of the formula for $N_{\mu_1 \circlearrowright \mu_2}(n)$
is equal to the number of all such f-walks $c$ of color 1 which have $r$ edges,
interlaced with $r$ alternating f-walks of lenghts $2k_{1}-1, 2k_{2}-1, \ldots , 2k_{r}-1$. 
The summation over $k_{1}+k_{2}+ \ldots + k_{r}=n$ indicates that
all alternating f-walks involved must have $2n-r$ edges together. 
The summation over $1\leq r \leq n$  gives exactly the cardinality of 
$D_{2n}(e)$, which finishes the proof.\hfill $\blacksquare$\\
\indent{\par}
Let us observe here that (3.1) can be generalized to arbitrary random variables 
in noncommutative probability spaces with distributions $\mu_1, \mu_2\in \Sigma$
if one treats $P_1$ and $P_2$ as idempotents and uses extensions of states 
$\varphi_1, \varphi_2$ given by (2.9). Then, the combinatorial formula
in the proof of Theorem 3.1 remains valid for $\mu_1, \mu_2\in \Sigma$.
\begin{Example}
{\rm If $\mu_1, \mu_2\in\Sigma$, the lowest order moments of 
$N_{\mu_1 \circlearrowright \mu_2}$ are given by
\begin{eqnarray*}
N_{\mu_1 \circlearrowright \mu_2}(1) &=& N_{\mu_1}(1)N_{\mu_2}(1),\\
N_{\mu_1 \circlearrowright \mu_2}(2) &=& N_{\mu_1}(2)N_{\mu_2}^{2}(1)+ N_{\mu_1}(1)N_{\mu_2}(2),\\
N_{\mu_1 \circlearrowright \mu_2}(3) &=& N_{\mu_1}(3)N_{\mu_2}^{3}(1) +2N_{\mu_1}(2)N_{\mu_2}(2)N_{\mu_2}(1)+
N_{\mu_1}(1)N_{\mu_2}(3),\\
N_{\mu_1 \circlearrowright \mu_2}(4) &=& N_{\mu_1}(4)N_{\mu_2}^{4}(1) +
3N_{\mu_1}(3)N_{\mu_2}(2)N_{\mu_2}^{2}(1)+ N_{\mu_1}(2)N_{\mu_2}^{2}(2)\\
&& +2N_{\mu_1}(2)N_{\mu_2}(3)N_{\mu_2}(1)+N_{\mu_1}(1)N_{\mu_2}(4).
\end{eqnarray*}
Let us apply these formulas to the graph ${\mathcal G}_{1}\vartriangleright_{\ell} {\mathcal G}_{2}$ in Fig.1.
Let $\mu_1$ and $\mu_2$ be the spectral distributions associated with ${\mathcal G}_{1}$ and ${\mathcal G}_{2}$,
respectively. By counting f-walks on ${\mathcal G}_{1}$ and ${\mathcal G}_{2}$, we easily get
$N_{\mu_1}(1)=N_{\mu_1}(2)=N_{\mu_1}(4)=1$, $N_{\mu_1}(3)=0$ and $N_{\mu_2}(1)=N_{\mu_2}(2)=1$,
$N_{\mu_2}(3)=0$, $N_{\mu_2}(4)=2$. In turn, using the above formulas and Theorem 3.1, we get
$D_{2}(e)=1$, $D_{4}(e)=D_{6}(e)=2$ and $D_{8}(e)=4$, which can also be verified directly by 
counting comb walks on ${\mathcal G}_{1}\vartriangleright_{\ell} {\mathcal G}_{2}$.}
\end{Example}
For completeness, let us also briefly discuss the case of the star product of graphs 
${\mathcal G}_{1}\star {\mathcal G}_{2}$ and find its relation to the boolean mutliplicative
convolution. 
\begin{Definition}
{\rm The {\it star product} of rooted graphs $({\mathcal G}_{1},e_1)$ and $({\mathcal G}_{2},e_2)$
is the rooted graph $({\mathcal G}_{1}\star {\mathcal G}_{2},e)$
obtained by glueing ${\mathcal G}_{1}$ and ${\mathcal G}_{2}$
at their roots and taking this common root to be the root of the product.
The {\it star loop product} $({\mathcal G}_{1}\star_{\ell} {\mathcal G}_{2},e)$
is obtained from naturally colored $({\mathcal G}_{1}\star {\mathcal G}_{2},e)$
by adding loops of color $\iota$ to every vertex but the root of ${\mathcal G}_{\overline{\iota}}$,
where $\iota \in I$.}
\end{Definition}
\indent{\par}
Recall here the definition of the boolean multiplicative convolution
of distributions $\mu_1, \mu_2\in \Sigma$ [5]. Treating $\eta_{\mu}(z)$
as a formal power series, we define 
\begin{equation}\tag{3.3}
\rho_{\mu}(z)=\frac{\eta_{\mu}(z)}{z},
\end{equation}
which allows us to define the boolean multiplicative 
convolution of $\mu_1,\mu_2\in \Sigma$ by the formula
\begin{equation}\tag{3.4}
\rho_{\mu_1\boxasterisk \mu_2}(z)=\rho_{\mu_1}(z)\rho_{\mu_2}(z).
\end{equation}
This formula defines $\boxasterisk$ as a binary operation
on ${\mathcal M}_{{\mathbb T}}$ [11].
However, if $\mu_1, \mu_2\in {\mathcal M}_{{\mathbb R}_{+}}$, 
then the RHS does not always give
the $\rho$-transform of a probability measure on ${\mathbb R}_{+}$ [5].

Nevertheless, the well-known relation [8,24] 
\begin{equation}\tag{3.5}
\rho_{\mu_1\, \boxtimes \,\mu_2}(z)=\rho_1(z)\rho_2(z)
\end{equation}
where $\rho_{\iota}(z)=\eta_{\iota}(z)/z$, $\iota \in I$, 
can still be written in terms of convolutions as
\begin{equation}\tag{3.6}
\mu_1 \boxtimes \mu_2 = (\mu_1 \boxslash \mu_2) \boxasterisk (\mu_2 \boxslash \mu_1).
\end{equation}
which was our motivation to use the symbol $\boxslash$ for the s-free
convolution (an analogous formula holds for $\mu_1,\mu_2\in {\mathcal M}_{*}$).
\begin{figure}
\unitlength=1mm
\special{em.linewidth 2pt}
\linethickness{0.5pt}
\begin{picture}(60.00,30.00)(20.00,0.00)
\put(0.00,10.00){\line(1,0){20.00}}
\put(35.00,10.00){\line(0,1){6.00}}

\put(0.00,10.00){\circle*{1.60}}
\put(0.00,8.00){\circle{4.00}}
\put(10.00,10.00){\circle*{1.00}}
\put(20.00,10.00){\circle*{1.00}}

\put(35.00,10.00){\circle*{1.60}}
\put(35.00,6.50){\circle{7.00}}
\put(35.00,16.00){\circle*{1.00}}
\put(30.00,20.00){\circle*{1.00}}
\put(40.00,20.00){\circle*{1.00}}

\put(65.00,22.00){\circle{4.00}}
\put(70.00,18.00){\circle{4.00}}
\put(75.00,22.00){\circle{4.00}}

\put(1.00,11.00){\scriptsize $e_{1}$}
\put(71.00,11.00){\scriptsize $e$}
\put(10.00,11.00){\scriptsize $x$}
\put(20.00,11.00){\scriptsize $x'$}
\put(37.00,10.00){\scriptsize $e_{2}$}
\put(36.00,13.50){\scriptsize $y$}
\put(39.00,22.00){\scriptsize $y''$}
\put(29.00,22.00){\scriptsize $y'$}
\put(9.00,15.00){\scriptsize ${\mathcal G}_{1}$}
\put(26.00,15.00){\scriptsize ${\mathcal G}_{2}$}
\put(52.00,15.00){\scriptsize ${\mathcal G}_{1}\star_{\ell} {\mathcal G}_{2}$}





\put(70.00,10.00){\line(1,0){30.00}}
\put(70.00,10.00){\line(0,1){6.00}}

\put(70.00,8.00){\circle{4.00}}
\put(70.00,6.50){\circle{7.00}}
\put(85.00,6.50){\circle{7.00}}
\put(100.00,6.50){\circle{7.00}}

\put(70.00,16.00){\circle*{1.00}}
\put(75.00,20.00){\circle*{1.00}}
\put(65.00,20.00){\circle*{1.00}}

\qbezier(65,20)(65,16)(70,16)
\qbezier(70,16)(75,16)(75,20)

\qbezier(30,20)(30,16)(35,16)
\qbezier(35,16)(40,16)(40,20)

\put(70.00,10.00){\circle*{1.60}}
\put(85.00,10.00){\circle*{1.00}}
\put(100.00,10.00){\circle*{1.00}}

\end{picture}
\caption{An example of the star loop product}
\end{figure}
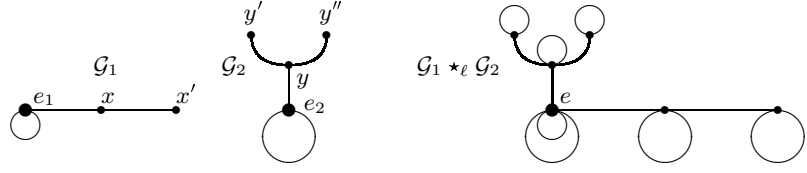
\indent{\par}
If $a_{\iota}$ is the adjacency matrix of a graph ${\mathcal G}_{\iota}$, where 
$\iota \in I$, then, using a similar reasoning as in the case of the comb loop product, 
we introduce operators of the form
\begin{eqnarray*}
R_{1}&=&a_{1}\otimes P_{2}+1 \otimes P_{2}^{\perp}, \\
R_{2}&=&P_{1}\otimes a_2+P_1^{\perp} \otimes 1_{2},
\end{eqnarray*}
which can be viewed as the adjacency matrices of the $1$- and $2$-colored subgraphs 
of the product graph $({\mathcal G}_{1}\star_{\ell} {\mathcal G}_{2},e)$. Moreover, 
these are exactly the operators which give the moments of the boolean multiplicative convolution
of spectral distributions of ${\mathcal G}_{1}$ and ${\mathcal G}_{2}$ since the pair
$(R_1-1,R_2-1)$ is boolean independent w.r.t. $\varphi_{e}$. 
An example of a star loop product is given in Fig.2. 

Let us also note that formulas (3.3)-(3.4) can be used if $a_1,a_2$
are random variables in arbitrary noncommutative probability spaces --
this is obtained if one uses extended states (2.9). This gives a realization 
of the boolean multiplicative convolution of arbitrary distributions 
$\mu_1, \mu_2\in \Sigma$.
\begin{Proposition}
The Multiplication Theorem holds for boolean independence, the associated
loop product ${\mathcal G}_{1}\star_{\ell}{\mathcal G}_{2}$, and the multiplicative convolution 
$\mu_1\boxasterisk \mu_2$.
\end{Proposition}
{\it Proof.}
The proof is similar to that of Theorem 3.1 and is based on the easy 
combinatorial formula 
$$
N_{\mu_1\boxasterisk\, \mu_2}(n)=\sum_{j+k=n+1}N_{\mu_1}(j)N_{\mu_2}(k),
$$
where the sum runs over positive indices $j,k$, obtained from (1.2) and (3.4).
Details are left to the reader.
\hfill $\blacksquare$
\begin{Example}
{\rm If $\mu_1, \mu_2\in\Sigma$, the lowest order `first return moments' 
$N_{\mu_1 \boxasterisk \,\mu_2}(n)$ are given by
\begin{eqnarray*}
N_{\mu_1 \boxasterisk\, \mu_2}(1) &=& N_{\mu_1}(1)N_{\mu_2}(1),\\
N_{\mu_1 \boxasterisk\, \mu_2}(2) &=& N_{\mu_1}(1)N_{\mu_2}(2)+ N_{\mu_1}(2)N_{\mu_2}(1),\\
N_{\mu_1 \boxasterisk\, \mu_2}(3) &=& N_{\mu_1}(1)N_{\mu_2}(3) +N_{\mu_1}(2)N_{\mu_2}(2)+
N_{\mu_1}(3)N_{\mu_2}(1),\\
N_{\mu_1 \boxasterisk\, \mu_2}(4) &=& N_{\mu_1}(1)N_{\mu_2}(4) +
N_{\mu_1}(2)N_{\mu_2}(3) + N_{\mu_1}(3)N_{\mu_2}(2)+N_{\mu_1}(4)N_{\mu_2}(1).
\end{eqnarray*}
These formulas can be used to count rooted alternating d-walks on the product
graph given in Fig. 2. Using Proposition 3.2 and 
the values of $N_{\mu_1}(j)$ and $N_{\mu_2}(k)$, given in Example 3.1,
we get $D_{2}(e)=D_{4}(e)=D_{6}(e)=1$
and $D_{8}(e)=3$, which can be also verified by direct computations.}
\end{Example}

\section{Orthogonal multiplicative convolution}

For given positive bounded random variables $X_{1}$ and $X_{2}$ 
with distributions $\mu_1$ and $\mu_2$,
we want to find a positive random variable $Z$ on some Hilbert space whose distribution is given by
the multiplicative analog of the orthogonal convolution, denoted $\mu_1 \angle \mu_2$.

Let us first recall the notion of orthogonal independence [16].
\begin{Definition}
{\rm Let $({\mathcal A},\varphi, \psi)$ be a unital algebra with a pair
of linear normalized functionals and let ${\mathcal A}_{1}$ and ${\mathcal A}_{2}$
be non-unital subalgebras of ${\mathcal A}$.
We say that ${\mathcal A}_{2}$ is {\it orthogonal} to ${\mathcal A}_{1}$
with respect to $(\varphi, \psi)$ if\\
\indent{\par}
(i)
$\;\varphi(bw_{2})=\varphi(w_{1}b)=0$,
\indent{\par}
(ii)
$\varphi(w_{1}a_{1}ba_{2}w_{2})=\psi(b)
\left(\varphi(w_{1}a_{1}a_{2}w_{2})- \varphi(w_{1}a_{1})\varphi(a_{2}w_{2})\right)$,\\[5pt]
for any $a_{1},a_{2}\in {\mathcal A}_{1}$, $b\in {\mathcal A}_{2}$,
and any elements $w_1,w_2$ of the unital algebra ${\rm alg}({\mathcal A}_{1},{\mathcal A}_{2})$
generated by ${\mathcal A}_{1}$ and ${\mathcal A}_{2}$.
We say that the pair $(a,b)$ of elements of ${\mathcal A}$ is {\it orthogonal} with respect to $(\varphi, \psi)$
if the algebra generated by $b\in {\mathcal A}$ is orthogonal to the algebra generated by
$a\in {\mathcal A}$.}
\end{Definition}
\begin{Definition}
{\rm Let $({\mathcal A}, \varphi, \psi)$ be a unital algebra with a pair of normalized linear functionals and let
${\mathcal A}_{1}, {\mathcal A}_{2}$ be a pair of (in general, non-unital) subalgebras of ${\mathcal A}$, 
such that ${\mathcal A}_{2}$ is orthogonal to ${\mathcal A}_{1}$ w.r.t. $(\varphi, \psi)$. 
Let $a_1,a_2\in {\mathcal A}$ be random variables with
$\varphi$-distribution $\mu_1$ and $\psi$-distribution $\mu_2$, respectively. The {\it orthogonal multiplicative
convolution} of $\mu_1$ and $\mu_2$, denoted $\mu_1 \angle \mu_2$, is the distribution of $a_2 a_1$, where
$(a_1-1, a_2-1)$ is orthogonal w.r.t. $(\varphi, \psi)$.}
\end{Definition}
In order to find a Hilbert space realization of the 
orthogonal multiplicative convolution of compactly supported probability measures,
recall the definition of the orthogonal product of two Hilbert spaces with
distinguished unit vectors $({\mathcal H}_{1}, \xi_{1})$ and $({\mathcal H}_{2}, \xi_{2})$.
Thus, the {\it orthogonal product} of $({\mathcal H}_{1}, \xi_{1})$ and $({\mathcal H}_{2}, \xi_{2})$
is the pair $({\mathcal H}, \xi)$, where
\begin{equation}\tag{4.1}
{\mathcal H}={\mathbb C}\xi \oplus {\mathcal H}_{1}^{0}\oplus ({\mathcal H}_{2}^{0}\otimes {\mathcal H}_{1}^{0}),
\end{equation}
with ${\mathcal H}_{i}^{0}={\mathcal H}_{i}\ominus {\mathbb C}\xi_{i}$ and where $\xi$ is a unit vector.
We denote it by $({\mathcal H}_{1}, \xi_{1})\vdash ({\mathcal H}_{2}, \xi_{2})$,
or simply ${\mathcal H}_{1}\vdash {\mathcal H}_{2}$, if no confusion arises.
It has been shown in [16] that ${\mathcal H}$ gives a Hilbert space realization of orthogonal random variables.
For that purpose we defined an isometry
$U: {\mathcal H}\rightarrow {\mathcal H}_{1}\otimes {\mathcal H}_{2}$ given by
\begin{equation}\tag{4.2}
U(\xi)=\xi_{1}\otimes \xi_{2},\;\; U(h_{1})=h_{1}\otimes \xi_{2},\;\; U(h_{2}\otimes h_{1})=h_{1}\otimes h_{2}
\end{equation}
for any $h_{1}\in {\mathcal H}_{1}^{0}$ and $h_{2}\in {\mathcal H}_{2}^{0}$. 
Using $U$, we define faithful non-unital *-homomorphisms
$\tau_{i}:{\mathcal B}({\mathcal H}_{i})\rightarrow {\mathcal B}({\mathcal H})$
by
\begin{equation}\tag{4.3}
\tau_{1}(a_1)= U^{*}(a_1\otimes P_{2}) U, \;\;\;
\tau_{2}(a_2)= U^{*}(P_{1}^{\perp} \otimes a_2)U,
\end{equation}
where $P_{1}$, $P_{2}$ are the projections
onto ${\mathbb C}\xi_{1}$ and ${\mathbb C}\xi_{2}$, respectively.
Then the pair $(\tau_{1}(a_1),\tau_{2}(a_2))$ is orthogonal w.r.t. $(\varphi, \psi)$, where
$\varphi$ and $\psi$ are states on ${\mathcal B}({\mathcal H})$ associated with the unit vector $\xi$
and any unit vector $\zeta\in {\mathcal H}_{1}^{0}$, respectively. For details, see [16].

For simplicity, we will identify vectors from ${\mathcal H}$ with their images 
in ${\mathcal H}_{1}\otimes {\mathcal H}_{2}$ under the isometry $U$ as well as bounded operators
$a, b$ with $\tau_{1}(a_1), \tau_{2}(a_2)$,
respectively, and we will carry out computations in ${\mathcal H}_{1}\otimes {\mathcal H}_{2}$. 
\begin{Proposition}
Let $a_1\in B({\mathcal H}_{1})$ have $\varphi_{1}$-distribution
$\mu_1$ and $a_2\in B({\mathcal H}_{2})$ -- $\varphi_{2}$-distribution $\mu_2$.
Then there exist a Hilbert space ${\mathcal H}$ with unit vectors $\xi$ and $\zeta$,
and random variables $R_1, R_2\in B({\mathcal H})$ such that $R_1$ has 
$\varphi$-distribution $\mu_1$, $R_2$ has $\psi$-distribution $\mu_2$ and
$R_2R_1$ has $\varphi$-distribution $\mu_1\angle \mu_2$, where $\varphi$ and $\psi$
are states associated with $\xi$ and $\zeta$, respectively.
\end{Proposition} 
{\it Proof.}
Let ${\mathcal H}={\mathcal H}_{1}\vdash {\mathcal H}_{2}$ with the canonical unit 
vector $\xi=\xi_{1}\otimes \xi_{2}$ and let $\zeta=\zeta_{1}\otimes \xi_{2}$, 
where $\zeta_{1}$ is an arbitrary unit vector from ${\mathcal H}_{1}^{0}$.
Define 
\begin{eqnarray*}
R_1&=&a_1\otimes P_{2}+1_{1}\otimes P_{2}^{\perp}, \\
R_2&=&P_{1}^{\perp}\otimes a_2+P_{1}\otimes 1_{2},
\end{eqnarray*}
where $P_{i}^{\perp}=1_{i}-P_{i}$ and $P_{i}$ is the canonical projection onto ${\mathbb C}\xi_{i}$
and $1_{i}$ denotes the identity in $B({\mathcal H}_{i})$, $i=1,2$.  
Therefore, denoting $1=1_{1}\otimes 1_{2}$, we get
\begin{eqnarray*}
R_1-1&=&(a_1-1_{1})\otimes P_{2},\\
R_2-1&=&P_{1}^{\perp}\otimes (a_2-1_{2}),
\end{eqnarray*} 
which, due to the appropriate tensor product form, implies that the pair $(R_1-1, R_2-1)$
is orthogonal w.r.t. $(\varphi, \psi)$, where $\varphi$ is the state associated with 
$\xi$ and $\psi$ is the state associated with $\zeta$. 
Moreover, elementary calculations show that the $\varphi$-distribution of 
$R_1$ is $\mu_1$ and that the $\psi$-distribution of $R_2$ is $\mu_2$,
which completes the proof.\hfill $\blacksquare$
\begin{Corollary}
Let $\mu_1, \mu_2$ be compactly supported probability measures on ${\mathbb R}_{+}$.
Then there exist bounded positive random variables $R_1,R_2$ on some Hilbert space
${\mathcal H}$ and unit vectors $\xi,\zeta\in {\mathcal H}$, such that $R_1$ has $\varphi$-distribution
$\mu_1$, $R_2$ has $\psi$-distribution $\mu_2$ and $\sqrt{R_2}R_1\sqrt{R_2}$
has $\varphi$-distribution $\mu_1\angle \mu_2$.
\end{Corollary}
{\it Proof.}
There exist positive operators $x_1\in B(L^{2}({\mathbb R}_{+},\mu_1))$ and 
$x_2\in B(L^{2}({\mathbb R}_{+}, \mu_2))$
(standard multiplication operators) which have distributions $\mu_1$ and $\mu_2$, respectively, with respect
to the states $\varphi_{1}$ and $\varphi_{2}$ 
associated with the constant functions equal to one, $h_{1}$ and $h_{2}$, respectively.
We can then use the framework of Proposition 4.1 to get 
self-adjoint bounded random variables on 
${\mathcal H}=L^{2}({\mathbb R}, \mu_1)\otimes L^{2}({\mathbb R},\mu_2)$  of the form
\begin{eqnarray*}
R'_{1}&=&(x_1-1_{1})\otimes P_{2},\\
R'_{2}&=&P_{1}^{\perp}\otimes (x_2-1_{2}),
\end{eqnarray*}
which are orthogonal w.r.t. $(\varphi, \psi)$, where $\varphi$ is associated with
vector $\xi=h_{1}\otimes h_{2}$ and $\psi$ is associated with $\zeta= f\otimes h_{2}$, where
$f$ is any function such that $\int_{{\mathbb R}}f(x)d\mu_1(x)=0$ and $\int_{{\mathbb R}}f^{2}(x)d\mu_1(x)=1$.
Now, let $R_1=R'_{1}+1$ and $R_2=R'_{2}+1$. It follows from the proof of Proposition 4.1 that
$R_1$ and $R_2$ have distributions $\mu_1$ and $\mu_2$, respectively. Moreover, both
are positive operators and thus, in particular, $\sqrt{R_2}$ exists. In fact,
$$
\sqrt{R_2}=P_{1}^{\perp}\otimes \sqrt{a_2}+P_{1}\otimes 1_{2},
$$
and $\sqrt{R_2}R_1\sqrt{R_2}$ is a bounded positive operator on ${\mathcal H}$. 
Moreover, $\sqrt{R_2}R_1\sqrt{R_2}$ has the same $\varphi$-distribution as $R_1R_2$,
which is obtained from the following calculation:
\begin{eqnarray*}
\varphi((\sqrt{R_2}R_1\sqrt{R_2})^{n})&=&
\varphi(\sqrt{R_2}(R_1R_2)^{n-1}R_{1}\sqrt{R_2})\\
&=&
(\varphi_{1}\otimes \varphi_{2})((\sqrt{P_1}\otimes 1_{2})(R_1R_2)^{n-1}R_{1}(P_{1}\otimes 1_{2}))\\
&=&
(\varphi_{1}\otimes \varphi_{2})((R_1R_2)^{n-1}R_{1}(P_1\otimes 1_{2}))\\
&=&
\varphi  ((R_1R_2)^{n}).
\end{eqnarray*}
By using the adjoints, we can interchange the order of $R_{1}$ and $R_{2}$ under 
the symbol $\varphi$. Therefore, the probability distribution 
of $\sqrt{R_2}R_1\sqrt{R_2}$ extends to the compactly supported measure $\mu_1 \angle \mu_2$ 
on ${\mathbb R}_{+}$.\hfill $\blacksquare$\\
\indent{\par}
We will now compute the $\varphi$--distribution of the product $Z:=R_2R_1$, where
$R_1$ and $R_2$ are given by Proposition 4.1, and derive
a formula for $\eta_{Z}$. In this context, let us observe that one can 
generalize Proposition 4.1 to arbitrary random variables in noncommutative
probability spaces since formulas for $R_{1}$ and $R_{2}$ remain valid provided
one takes extensions of states given by (2.9) and uses Proposition 2.2.
Therefore, further developments in this Section will hold 
for arbitrary random variables in noncommutative probability spaces.
This corresponds to taking distributions $\mu\in \Sigma$ and the 
corresponding formal power series $\eta_{\mu}$ and $\rho_{\mu}$.

In order to compute the $\varphi$--distribution of $Z$, let us decompose $Z$ as
$$
Z=Z_{1}+Z_{2}+Z_{3}+Z_{4},
$$
where
$$
Z_{1}=P_{1}^{\perp}a_1\otimes a_2P_{2}, \;\;\;
Z_{2}=P_{1}^{\perp}\otimes a_2P_{2}^{\perp},\;\;\;
Z_{3}=P_{1}a_1\otimes P_{2}, \;\;\;
Z_{4}=P_{1}\otimes P_{2}^{\perp}.
$$
In order to compute the `first-return moments' $N_{Z}(n)$ in the state 
$\varphi$, we need to compute the mixed $\varphi$--moments of $Z_{i}$'s. 
However, it turns out that there are only two types of them which give a non-vanishing 
contribution. We compute them in Propositions 4.3-4.4.
\begin{Proposition}
For $r$ odd and $j_{1},j_{2}, \ldots , j_{r}\geq 1$, we have
$$
\varphi(Z_{3}Z_{1}^{j_{1}}Z_{2}^{j_{2}}Z_{1}^{j_{3}}\ldots Z_{1}^{j_{r}})=
N_{a_1}(m)N_{a_2}(k_{1})\ldots N_{a_2}(k_{m-1}),
$$
where 
$m=1+\sum_{k\,{\rm odd}}j_{k}$
and all $k_{i}$'s are equal to $1$, except
$$
k_{j_{1}+1}=j_{2}+1,\;\; k_{j_{1}+j_{3}+1}=j_{4}+1, 
\ldots, k_{j_{1}+j_{3}+\ldots +j_{r-2}+1}=j_{r-1}+1,
$$
and where the `first return moments' $N_{a_{\iota}}(s)$ refer to the state $\varphi_{\iota}$, $\iota \in I$.
\end{Proposition}
{\it Proof.}
We have 
\begin{eqnarray*}
\varphi(Z_{3}Z_{1}^{j_{1}}Z_{2}^{j_{2}}Z_{1}^{j_{3}}\ldots Z_{1}^{j_{r}})
&=&
\varphi_{1}(P_{1}a_1(P_{1}^{\perp}a_1)^{j_{1}}P_{1}^{\perp}(P_{1}^{\perp}a_1)^{j_{3}}\ldots P_{1}^{\perp}(P_{1}^{\perp}a_1)^{j_{r}})\\
&&
\varphi_{2} (P_{2}(a_2P_{2})^{j_{1}}(a_2P_{2}^{\perp})^{j_{2}}(a_2P_{2})^{j_3-1}
\ldots (a_2P_{2})^{j_{r}})\\
&=&
\varphi_{1}(a_1(P_{1}^{\perp}a_1)^{j_{1}+j_{3}+\ldots +j_{r}})\\
&&
\varphi_{2}((a_2P_{2})^{j_{1}}(a_2P_{2}^{\perp})^{j_{2}}(a_2P_{2})^{j_{3}}\ldots (a_2P_{2})^{j_{r}})\\
&=&
N_{a_1}(1+j_{1}+j_{3}+\ldots +j_{r})(N_{a_2}(1))^{j_{1}+(j_{3}-1)+\ldots +(j_{r}-1)}\\
&&
N_{a_2}(j_{2}+1)N_{a_2}(j_{4}+1)\ldots N_{a_2}(j_{r-1}+1),
\end{eqnarray*}
using the following properties: 
$$
P_{1}a_1P_{1}=\varphi_{1}(a_1)P_{1}=N_{a_1}(1)P_{1}, \;\;\;
P_{2}a_2P_{2}=\varphi_{2}(a_2)P_{2}=N_{a_2}(1)P_{2}.
$$ 
Setting $k_{1}=k_{2}=\ldots =k_{j_{1}}=1$, $k_{j_{1}+1}=j_{2}+1$, $k_{j_{1}+2}=\ldots =k_{j_{1}+j_{3}}=1$,
$k_{j_{1}+j_{3}}=1$, $k_{j_{1}+j_{3}+1}=j_{4}+1$, $\ldots$, $k_{j_{1}+j_{3}+\ldots +j_{r-2}+1}=j_{r-1}+1$,
$k_{j_{1}+j_{3}+\ldots +j_{r-2}+2}=\ldots = k_{m-1}=1$, we get the desired result.
\hfill $\blacksquare$\\
\begin{Proposition}
For $r$ even and $j_{1},j_{2}, \ldots , j_{r}\geq 1$, we have
$$
\varphi(Z_{3}Z_{2}^{j_{1}}Z_{1}^{j_{2}}Z_{2}^{j_{3}}\ldots Z_{1}^{j_{r}})=
N_{a_1}(m)N_{a_2}(k_{1})\ldots N_{a_2}(k_{m-1}),
$$
where 
$m=1+\sum_{k\,{\rm even}}j_{k}$
and all $k_{i}$'s are equal to $1$, except
$$
k_{j_{2}+1}=j_{3}+1,\;\; k_{j_{2}+j_{4}+1}=j_{5}+1, 
\ldots, k_{j_{2}+j_{4}+\ldots +j_{r-2}+1}=j_{r-1}+1,
$$
and where the `first return moments' $N_{a_{\iota}}(s)$ refer to the state $\varphi_{\iota}$, $\iota \in I$.
\end{Proposition}
{\it Proof.}
Using similar arguments as in the proof of Proposition 4.3, we get
\begin{eqnarray*}
\varphi(Z_{3}Z_{2}^{j_{1}}Z_{1}^{j_{2}}Z_{2}^{j_{3}}\ldots Z_{1}^{j_{r}})
&=&
\varphi_{1}(P_{1}a_1P_{1}^{\perp}(P_{1}^{\perp}a_1)^{j_{2}}P_{1}^{\perp}(P_{1}^{\perp}a_1)^{j_{4}}\ldots (P_{1}^{\perp}a_1)^{j_{r}})\\
&&
\varphi_{2} (P_{2}(a_2P_{2}^{\perp})^{j_{1}}(a_2P_{2})^{j_{2}}(a_2P_{2}^{\perp})^{j_3}
\ldots (a_2P_{2})^{j_{r}})\\
&=&
\varphi_{1}(a_1(P_{1}^{\perp}a_1)^{j_{2}+j_{4}+\ldots +j_{r}})\\
&&
\varphi_{2}((a_2P_{2}^{\perp})^{j_{1}}(a_2P_{2})^{j_{2}}(a_2P_{2}^{\perp})^{j_{3}}\ldots (a_2P_{2})^{j_{r}})\\
&=&
N_{a_1}(1+j_{2}+j_{4}+\ldots +j_{r})(N_{a_2}(1))^{(j_{2}-1)+(j_{4}-1)+\ldots +(j_{r}-1)}\\
&&
N_{a_2}(j_{1}+1)N_{a_2}(j_{3}+1)\ldots N_{a_2}(j_{r-1}+1).
\end{eqnarray*} 
Setting $k_{1}=j_{1}+1$, $k_{2}=\ldots =k_{j_{2}}=1$, $k_{j_{2}+1}=j_{3}+1$, 
$k_{j_{2}+2}=\ldots =k_{j_{2}+j_{4}}=1$,
$k_{j_{2}+j_{4}+1}=j_{5}+1$, $\ldots $, 
$k_{j_{2}+j_{4}+\ldots +j_{r-2}+1}=j_{r-1}+1$,
$k_{j_{2}+j_{4}+\ldots +j_{r-2}+2}=\ldots = k_{m-1}=1$, we get the desired result.
\hfill $\blacksquare$
\begin{Theorem}
If the $\varphi$--distribution of $a_{2}$ is not concentrated at zero,
then the formal power series corresponding to random variables $a_{1}, a_{2}$ and $Z$ 
satisfy the relation
\begin{equation}\tag{4.4}
\eta_{Z}(z)=\frac{z\eta_{a_1}(\eta_{a_2}(z))}{\eta_{a_2}(z)}.
\end{equation}
If the distribution of $a_{2}$ is concentrated at zero, then $\eta_{Z}(z)=z\mu_1(X)$,
where $\mu_1$ is the $\varphi$--distribution of $a_{1}$.
\end{Theorem}
{\it Proof.}
Suppose the $\varphi$--distribution of $a_2$ is not the Dirac delta.
Then, writing the RHS of equation (4.4) as $D(z)=\sum_{n=1}^{\infty}d_{n}z^{n}$, 
we need to show that $d_{n}=\eta_{Z}(n)$ for $n \in {\mathbb N}$. We clearly have
$$
N_{Z}(1)=\varphi(Z_{3})=\varphi_{1}(a_1)=N_{a_1}(1)=d_{1}.
$$ 
Now, for $n>1$, it holds that
$$
d_{n}=
\sum_{m=1}^{n}
N_{a_1}(m)\sum _{k_{1}+k_{2}+\ldots +k_{m-1}=n-1}
N_{a_1}(k_{1})N_{a_2}(k_{2})\ldots N_{a_2}(k_{m-1}),
$$ 
where all $k_{i}$'s are assumed to be positive integers, 
Using Proposition 2.2 and the explicit form of $Z_{i}$'s, we get
\begin{eqnarray*}
N_{Z}(n)
&=&\varphi (Z(P^{\perp}Z)^{n-1})\\
&=&\varphi(Z_{3}(Z_{1}+Z_{2})^{n-1})\\
&=&\varphi(Z_{3}(Z_{1}+Z_{2})^{n-2}Z_{1})\\
&=&
\sum_{r=1}^{n-2}
\sum_{i_{1}\neq i_{2}\neq \ldots \neq i_{r}}
\sum_{j_{1}+j_{2}+\ldots + j_{r}=n-2}
\varphi(Z_{3}Z_{i_{1}}^{j_{1}}Z_{i_{2}}^{j_{2}}\ldots Z_{i_{r}}^{j_{r}}Z_{1}),
\end{eqnarray*}
where $j_{1},j_{2}, \ldots , j_{r}$ are positive integers.
Note that on the RHS of the formula for $d_{n}$ we have the sum of products 
$$
N_{a_1}(m)N_{a_2}(k_{1})N_{a_2}(k_{2})\ldots N_{a_2}(k_{m-1}),
$$
where
$1\leq m \leq n$ and $(k_{1},k_{2}, \ldots , k_{m-1})$ is a sequence of positive integers which
add up to $n-1$. Now, every product obtained
in Propositions 4.3-4.4 is of exactly this form, with indices $(k_{1},k_{2}, \ldots , k_{m-1})$
satisfying the conditions just mentioned.
Moreover, all these products are different from each other and 
exhaust all possibile values of $(m, k_{1}, \ldots , k_{m-1})$. 
In fact, there are two disjoint cases distinguished: 
Proposition 4.3 covers the case $k_{1}=1$ and Proposition 4.4 -- the case $k_{1}\neq 1$.
In the first case, even $j_{i}$'s determine the subsequence of $(k_{1}, k_{2}, \ldots , k_{m-1})$ 
consisting of numbers which are greater than $1$ 
and odd $j_{i}$'s determine their positions in the sequence by determining all 
$k_{i}$'s which are equal to $1$.   
In the second case, odd $j_{i}$'s determine the subsequence of 
$(k_{1},k_{2}, \ldots , k_{m-1})$ consisting of numbers which are greater than $1$
and even $j_{i}$'s determine their positions in the sequence by determining all $k_{i}$'s 
which are equal to $1$. Therefore, both cases give together every tuple 
$(m, k_{1},k_{2}, \ldots , k_{m-1})$ under consideration exactly once.
This proves (4.4) in the case when the distribution 
of $a_2$ is not concentrated at zero.
Finally, if the $\varphi$--distribution of $a_{2}$ is $\delta_{0}$, then $\varphi_2(a_{2}^{n})=\delta_{n,0}$,
which easily gives $N_Z(n)=N_{a_1}(1)\delta_{n,1}$, where $\delta_{i,j}$ is the Kronecker symbol,
and since $N_{a_1}(1)=\varphi_1 (a_1)=\mu_1(X)$, we have $\eta_{Z}(z)=z\mu_1(X)$. 
This completes the proof.
\hfill $\blacksquare$\\
\indent{\par}
This shows that Eq.(1.7) holds for the orthogonal multiplicative 
convolution of distributions $\mu_1, \mu_2\in \Sigma$, provided $\mu_2\neq \delta_{0}$,
where the functions $\eta$ are treated as formal power series. 
In turn, explicit computations of distributions $\mu_1\angle \mu_2$ based on 
Definition 4.2 give
\begin{equation}\tag{4.5}
\delta_{0}\angle \mu=\delta_{0}\;\;\;{\rm and}\;\;\;\mu \angle \delta_{0}=\delta_{\mu(X)}
\end{equation}
for any $\mu\in \Sigma$. Let us also remark that Theorem 4.5 gives 
another formula for formal power series corresponding to $\varphi$--distributions, namely
\begin{equation}\tag{4.6}
\rho_{\mu_1\angle\mu_2}(z)=\rho_{\mu_1}(\eta_{\mu_2}(z)),
\end{equation}
where $\rho_{\mu_1}$ is given by (3.3), which turns out slightly more
useful in computations involving $\mu_1\angle \mu_2$ than (4.4) (see Section 5).
\begin{Corollary}
If $\mu_{1}$ and $\mu_{2}$ are $\varphi$--distributions of certain random variables,
then the $\varphi$--moment of $\mu_1\angle \mu_2$ of order $n\in {\mathbb N}$ depends on 
$\varphi$--moments of $\mu_1$ of orders $1\leq k \leq n$ and $\varphi$--moments of 
$\mu_2$ of orders $1\leq k \leq n-1$. 
\end{Corollary}
{\it Proof.}
This fact is a consequence of Proposition 4.1 and Theorem 4.5. In fact, the combinatorial formula for $d_{n}$
in the proof of Theorem 4.5 gives the explicit formula for $N_{\mu\angle\nu}(n)$ in terms of
$N_{\mu}(m)$ and $N_{\nu}(k)$ for $1\leq m \leq n$ and $1\leq k \leq n-1$.
Thus, in view of (1.2), a similar property holds for the $\varphi$--moments of $\mu\angle\nu$.
\hfill $\blacksquare$ 
\begin{Example}
{\rm If $\mu_1, \mu_2\in \Sigma$, the lowest order first return moments 
$N_{\mu_1 \angle \mu_2}(n)$ are given by
\begin{eqnarray*}
N_{\mu_1 \angle \mu_2}(1) &=& N_{\mu_1}(1),\\
N_{\mu_1 \angle \mu_2}(2) &=& N_{\mu_1}(2)N_{\mu_2}(1),\\
N_{\mu_1 \angle \mu_2}(3) &=& N_{\mu_1}(3)N_{\mu_2}^{2}(1)+N_{\mu_1}(2)N_{\mu_2}(2),\\
N_{\mu_1 \angle \mu_2}(4) &=& N_{\mu_1}(4)N_{\mu_2}^{3}(1)+2N_{\mu_1}(3)N_{\mu_2}(2)N_{\mu_2}(1)
+
N_{\mu_1}(2)N_{\mu_2}(3).
\end{eqnarray*}
When we return to the $\varphi$--moments $\mu_1\angle \mu_2 (X^{n})$, using (1.2),
we can express them in terms of `universal polynomials' 
$Q_{n}(\mu_1(X), \ldots , \mu_1(X^n), \mu_2(X), \ldots, \mu_2(X^{n-1}))$, but one should note that 
they are not homogenous as in the case of free convolutions [25].}
\end{Example}

\section{Convolutions of measures on ${\mathbb R}_{+}$ and ${\mathbb T}$.}

In this section we show that the orthogonal multiplicative convolution 
can be defined for arbitrary probability measures on ${\mathbb R}_{+}$ 
which are not concentrated at zero, 
as well as for arbitrary probability measures on the unit circle ${\mathbb T}$.

For the following result, we refer the reader to
the works of Belinschi and Bercovici [3,4]. 
\begin{Theorem}
There is a bijection between ${\mathcal M}_{{\mathbb R}_{+}}$ and the class of analytic 
self-maps $\eta$ of ${\mathbb C}\setminus {\mathbb R}_{+}$, such that
\begin{enumerate}
\item $\eta(\overline{z})=\overline{\eta(z)}$ {\it for all} $z\in {\mathbb C}\setminus {\mathbb R}_{+}$.
\item $\lim_{z\rightarrow 0^{-},z<0}\eta (z)=0$
\item ${\rm arg}\eta(z)\in [{\rm arg}z, \pi)$ {\it for all} $z\in {\mathbb C}^{+}$.
\end{enumerate}
Moreover, the map $\eta$ corresponding to $\mu\in {\mathcal M}_{{\mathbb R}_{+}}$ is given by
the transform $\eta_{\mu}$.
\end{Theorem}

We shall use this result to prove that formula (4.4), which was shown to
hold for formal power series corresponding to distributions of random variables,
can also be used to define $\mu_1\angle\mu_2$ for 
$\mu_1,\mu_2\in {\mathcal M}_{{\mathbb R}_{+}}$ and $\mu_2\neq \delta_{0}$,
if formal power series are replaced by analytic functions on ${\mathbb C}\setminus {\mathbb R}_{+}$. 
\begin{Proposition}
If $\mu_1, \mu_2\in {\mathcal M}_{{\mathbb R}_{+}}$ and $\mu_2$ is not concentrated at zero, 
then there exists a unique probability measure $\mu\in {\mathcal M}_{{\mathbb R}_{+}}$, such that
\begin{equation}\tag{5.1}
\eta_{\mu}(z)=\frac{z\eta_{\mu_{1}}(\eta_{\mu_{2}}(z))}{\eta_{\mu_{2}}(z)}
\end{equation}
for $z\in {\mathbb C}\setminus {\mathbb R}_{+}$. 
The measure $\mu$ will be defined to be the orthogonal multiplicative convolution
of $\mu_1$ and $\mu_2$, denoted $\mu_1\angle \mu_2$.
\end{Proposition}
{\it Proof.}
For $z\in {\mathbb C}\setminus {\mathbb R}_{+}$, let us define
$$
\eta(z)=\frac{z\eta_{\mu_1}(\eta_{\mu_2}(z))}{\eta_{\mu_2}(z)}.
$$
Since $\mu_1,\mu_2\in {\mathcal M}_{{\mathbb R}_{+}}$, transforms
$\eta_{\mu_1}$ and $\eta_{\mu_2}$ are analytic self-maps of
${\mathbb C}\setminus {\mathbb R}_{+}$ which satisfy the conditions of Theorem 5.1.
This implies that $\eta(z)$ is an analytic self-map of 
${\mathbb C}\setminus {\mathbb R}_{+}$ since $\mu_2\neq \delta_{0}$.
Condition (1) of Theorem 5.1 holds for $\eta_{\mu_1}$ and $\eta_{\mu_2}$ and therefore it holds
for $\eta$. 
Then, it is known that  
$$
0 \leq \lim_{z\rightarrow 0^{-}, \,z<0}
\frac{\eta_{\mu_2}(z)}{z}\leq \infty,
$$
and this limit is different from zero (it equals zero only if $\mu_2=\delta_{0}$).
Thus, we obtain 
\begin{eqnarray*}
\lim_{z\rightarrow 0^{-},\,z<0}\eta(z)
&=& 
\lim_{z\rightarrow 0^{-},\,z<0}
\frac{z}{\eta_{\mu_2}(z)}
\lim_{z\rightarrow 0^{-},\, z<0}
\eta_{\mu_1}(\eta_{\mu_2}(z))\\
&=&
c\cdot 
\lim_{z\rightarrow 0^{-},\,z<0}
\eta_{\mu_1}(\eta_{\mu_2}(z))\\
&=&
0
\end{eqnarray*}
for some $0\leq c <\infty$, where we used condition 2 of Theorem 5.1 for the function
$\eta_{\mu_1}$ (if $z\rightarrow 0^{-}$, $\psi_{\mu_2}(z)\rightarrow 0^{-}$ and 
thus $\eta_{\mu_2}(z)\rightarrow 0^{-}$ and therefore 
$\eta_{\mu_1}(\eta_{\mu_2}(z))\rightarrow 0$ as $z\rightarrow 0^{-}$).
This proves that condition (2) of Theorem 5.1 holds for the map $\eta$.

Next, for any $z\in {\mathbb C}^{+}$, we have $\psi_{\mu_2}(z)\in {\mathbb C}^{+}$, using the integral
representation of $\psi_{\mu_2}(z)$ given by (1.1). 
This implies that $\eta_{\mu_2}(z)\in {\mathbb C}^{+}$, since
$$
\eta_{\mu_2}(z)=\frac{\psi_{\mu_2}(z)}{1+\psi_{\mu_2}(z)}=\frac{\psi_{\mu_2}(z)+|\psi_{\mu_2}(z)|^{2}}{1+\psi_{\mu_2}(z)|^{2}},
$$
and therefore
$$
\Im \eta_{\mu_2}(z)=\frac{\Im \psi_{\mu_2}(z)}{|1+\psi_{\mu_2}(z)|^{2}}.
$$
Hence, applying Theorem 5.1 to the function $\eta_{\mu_1}$, we get 
$$
{\rm arg}\eta_{\mu_1}(\eta_{\mu_2}(z))\geq {\rm arg}\eta_{\mu_2}(z),
$$ 
and thus
$$
{\rm arg}\eta(z)={\rm arg}z+{\rm arg}\eta_{\mu_1}(\eta_{\mu_2}(z))-{\rm arg}\eta_{\mu_2}(z)\geq {\rm arg}z.
$$
Similarly, 
\begin{eqnarray*}
{\rm arg}\eta(z)&=&{\rm arg}\eta_{\mu_1}(\eta_{\mu_2}(z))+{\rm arg}z-{\rm arg}\eta_{\mu_2}(z)\\
&=&{\rm arg}\eta_{\mu_1}(\eta_{\mu_2}(z))-({\rm arg}\eta_{\mu_2}(z)-{\rm arg}z)\\
&<&\pi
\end{eqnarray*}
since we have 
$$
{\rm arg}\eta_{\mu_2}(z)\geq {\rm arg}z\;\;\; {\rm and}\;\;\;
{\rm arg}\eta_{\mu_1}(\eta_{\mu_2}(z))<\pi,
$$ 
which follows from $\eta_{\mu_2}(z)\in {\mathbb C}^{+}$. This proves that 
$\eta$ satisfies condition (3) of Theorem 5.1, which completes the proof.
\hfill $\blacksquare$\\
\indent{\par}
For computations, it is convenient to use (4.6) and represent the transforms involved as continued fractions. 
The continued fraction representation of $\eta_{\mu}$ for compactly supported
$\mu\in {\mathcal M}_{{\mathbb R}_{+}}$ is obtained from that of the Cauchy transform $G_{\mu}$, or its reciprocal $F_{\mu}$, by using the formula
\begin{equation}\tag{5.2}
\eta_{\mu}(z)=1-zF_{\mu}\left(\frac{1}{z}\right).
\end{equation}
If $\mu\in {\mathcal M}_{{\mathbb R}_{+}}$ is associated with sequences of Jacobi parameters 
$(\alpha, \omega)$, where $\alpha=(\alpha_{n})_{n\geq 0}$ and 
$\omega=(\omega_{n})_{n\geq 0}$, which we write $J(\mu)=(\alpha, \omega)$, 
then, using the continued fraction representation of $G_{\mu}$, we obtain 
\begin{equation}\tag{5.3}
\eta_{\mu}(z)=\alpha_{0}z+\cfrac{\omega_{0}z^2}{1-\alpha_{1}z-\cfrac{\omega_1 z^2}
{1-\alpha_{2}z-\cfrac{\omega_{2}z^2}{1-\alpha_{3}z-\cfrac{\omega_{3}z^2} {\ldots}}}}
\end{equation}
and a related continued fraction for $\rho_{\mu}$. Conversely, this formula, 
together with the first formula of (1.7), enables us to compute the Jacobi 
sequences corresponding to the orthogonal multiplicative convolution.
In the case of compactly supported measures, the Jacobi sequences uniquely determine
the corresponding measure.  
\begin{Example}
{\rm Let $\delta_{a}$ be the Dirac measure with $a>0$ and let $\mu \in {\mathcal M}_{{\mathbb R}_{+}}$
be compactly supported, with the associated Jacobi sequences $J(\mu)=(\alpha, \omega)$.
Then $\rho_{\mu}(z)=a$ and therefore, $\rho_{\delta_{a}\angle \mu}(z)=\rho_{\delta_{a}}(\eta_{\mu}(z))=a$,
which gives $\delta_{a}\angle \mu=\delta_{a}$.
In turn, $\rho_{\mu\angle \delta_{a}}(z)=\rho_{\mu}(az)$.
Let us denote the corresponding 
transformation of compactly supported measures by $S_{a}$.
In terms of reciprocal Cauchy transforms, it can be defined by 
a convex linear combination
$$
F_{S_{a}\mu}(z)=\frac{1}{a}F_{D_{a}\mu}(z)+(1-\frac{1}{a})z,
$$
where $D_{\lambda}\mu$ is the dilation of measure $\mu$ by $\lambda$, defined by
$D_{\lambda}\mu(E)=\mu(\lambda^{-1}E)$. In terms of the transforms $\eta$, we have $\eta_{S_{a}\mu}(z)=\eta_{D_{a}\mu}(z)/a$, using (5.2). We also have $\eta_{D_{a}\mu}(z)=\eta_{\mu}(az)$.
Note that $S_{a}\mu=(T_{a}\circ D_{a})\mu$, 
where $T_{t^{-1}}$ is the so-called $t$-transformation of measures [9].
If $\mu$ is compactly supported, we can write $\eta_{S_{a}\mu}$ 
in the form of a continued fraction
$$
\eta_{S_{a}\mu}(z)=\alpha_{0}z+\cfrac{\omega_{0}az^2}{1-\alpha_{1}az-\cfrac{\omega_1 a^{2}z^2}
{1-\alpha_{2}az-\cfrac{\omega_{2}a^2z^2}{1-\alpha_{3}az-\cfrac{\omega_{3}a^2z^2} {\ldots}}}}.
$$
In particular, $\mu\angle \delta_{1}=\mu$, thus $\delta_1$ is the right unit
w.r.t. the operation $\angle$ on ${\mathcal M}_{{\mathbb R}_{+}}$, and 
$\delta_{1}\angle \mu=\delta_{1}$.}
\end{Example}
\begin{Example}
{\rm If $\mu_1=(1-p)\delta_{0}+p\delta_{1}$ and $\mu_2=(1-q)\delta_{0}+q\delta_{1}$, then
the corresponding reciprocal Cauchy transforms are
$$
F_{\mu_1}(z)=\frac{z-1+p}{z(z-1)}, \;\;\;F_{\mu_2}(z)=\frac{z-1+q}{z(z-1)},
$$
which gives
$$
\rho_{\mu_1}=p+\frac{(p-p^2)z}{1-(1-p)z},\;\;\;
\eta_{\mu_2}=qz+\frac{(q-q^2)z^2}{1-(1-q)z},
$$
and therefore, in view of (4.6), we obtain
$$
\rho_{\mu_1\angle \mu_2}(z)=p+\frac{(p-p^2)qz}{1+(pq-1)z}.
$$
Using (5.3), we obtain $J(\mu_1\angle \mu_2)=(\alpha, \omega)$, with $\alpha=(p,1-pq, 0,0, \ldots )$
and $\omega=((p-p^2)q, 0,0, \ldots )$. The corresponding measure has two atoms (explicit dependence
on $p$ and $q$ is rather complicated and is omitted).}
\end{Example}

A result analogous to Theorem 5.1 also holds for the set ${\mathcal M}_{{\mathbb T}}$
of probability measures on the unit circle ${\mathbb T}=\{z\in {\mathbb C}:\, |z|=1\}$.
Here, the class of self-maps $\eta$ of the open unit disc 
${\mathbb D}=\{z\in {\mathbb C}:\,|z|<1\}$ is used, which is again related
to the transforms $\eta_{\mu}$. For details, we refer the reader to 
the works of Belinschi and Bercovici [3,4].
\begin{Theorem}
There is a bijection between ${\mathcal M}_{{\mathbb T}}$ and 
the class of analytic self-maps $\eta$ of ${\mathbb D}$ such that  $\eta(0)=0$ and 
$|\eta(z)|\leq |z|$ for all $z\in {\mathbb D}$. Moreover, the map $\eta$ corresponding 
to $\mu$ is given by the transform $\eta_{\mu}$.
\end{Theorem}
\begin{Proposition}
If $\mu_1, \mu_2\in {\mathcal M}_{{\mathbb T}}$, then there exists a unique 
probability measure $\mu\in {\mathcal M}_{{\mathbb T}}$ such that
\begin{equation}\tag{5.4}
\eta_{\mu}(z)=
\frac{z\eta_{\mu_{1}}(\eta_{\mu_{2}}(z))}{\eta_{\mu_{2}}(z)}
\end{equation}
for $z\in {\mathbb D}$, where we understand that $\eta_{\mu_1}(0)/0=\eta_{\mu_1}'(0)$.
This measure $\mu$ will be defined to be the orthogonal multiplicative 
convolution $\mu_1\angle \mu_2$.
\end{Proposition}
{\it Proof.}
First of all, observe that $|\eta(z)|\leq |z|$ since $\eta_{\mu_1}$ has the same property by Theorem 5.3.
Moreover, $\eta$ is a quotient of analytic functions on ${\mathbb D}$, and if $\eta_{\mu_2}(z_{0})=0$, then 
$$
\eta(z_{0})=\lim_{z\rightarrow z_{0}}z
\frac{\eta_{\mu_1}(\eta_{\mu_2}(z))}{\eta_{\mu_2}(z)}=z_{0}\eta_{\mu_1}'(0),
$$
thus $\eta$ is a well-defined analytic self-map of ${\mathbb D}$. Finally,
$\eta(0)=0\cdot \eta_{\mu_1}'(0)=0$.
Thus, $\eta$ satisfies the conditions of Theorem 5.3 and therefore
$\eta=\eta_{\mu}$ for some $\mu \in {\mathcal M}_{{\mathbb T}}$, which completes the proof.
\hfill $\blacksquare$
\begin{Example}
{\rm For any $\mu\in {\mathcal M}_{{\mathbb T}}$ and $a\in {\mathbb T}$, we get
$\delta_{a}\angle \mu=\delta_{a}$ and $\mu \angle \delta_{a}=S_{a}\mu$, where
$S_{a}\mu$ is defined by the same equation $\eta_{S_{a}\mu}(z)=\eta_{D_{a}\mu}(z)/a$
as in the case of measures on ${\mathbb R}_{+}$ 
(on ${\mathcal M}_{{\mathbb T}}$, the operation $D_{a}$ should be interpreted as a rotation).}
\end{Example}

\section{Orthogonal loop product of graphs}

In this section we establish a relation between the orthogonal 
multiplicative convolution of compactly supported 
probability measures and a new type of product of graphs 
called the `orthogonal loop product'.

As in the monotone case, we begin with recalling the definition of the orthogonal product 
of graphs [2,16]. Then we will modify this product in an appropriate manner to 
find a relation with the orthogonal multiplicative convolution.
\begin{Definition}
{\rm The {\it orthogonal product} of two rooted graphs $({\mathcal G}_{1},e_{1})$ and $({\mathcal G}_{2},e_{2})$ 
is the rooted graph $({\mathcal G}_{1}\vdash {\mathcal G}_{2},e)$ obtained by attaching a copy of ${\mathcal G}_{2}$ by its root $e_{2}$ to each vertex of ${\mathcal G}_{1}$ but the root $e_{1}$, where $e$ is taken to be equal to $e_{1}$.
If its set of vertices is identified with $V_{1}\vdash V_{2}:=(V_{1}^{0}\times V_{2})\cup \{e_{1}\times e_{2}\}$,
then $e$ is identified with $e_{1}\times e_{2}$.}
\end{Definition}

It is worth noting that the orthogonal product of graphs resembles their comb product. 
The difference is that in the comb product the second graph is glued by its root to {\it all }
vertices of the first graph, whereas in the orthogonal product the second graph is glued 
to {\it all vertices but the root} of the first graph. 
\begin{Definition}
{\rm Suppose the edges of the orthogonal product $({\mathcal G}_{1}\vdash {\mathcal G}_{2},e)$
are $I$--colored. 
The {\it orthogonal loop product} of rooted graphs 
$({\mathcal G}_{1},e_1)$ and $({\mathcal G}_{2},e_2)$
is the rooted graph $({\mathcal G}_{1}\vdash_{\ell} {\mathcal G}_{2},e)$
obtained from $({\mathcal G}_{1}\vdash {\mathcal G}_{2},e)$
by attaching a loop of color $1$ to all vertices but the root 
of each copy of ${\mathcal G}_{2}$, and a loop of color $2$ to the root
of ${\mathcal G}_{1}$.}
\end{Definition}

The justification of the above definition is similar to that in Remark 3.1.
It is enough to deduce the glueing rules from the expressions for $R_{1}$ and $R_{2}$,
given in the proof of Proposition 4.1, interpreting $a_{\iota}$ as the adjacency matrix of
graph ${\mathcal G}_{\iota}$, $\iota \in I$. Also, there is an easy analog
of Remark 3.1. Finally, it follows from [2] that the additive analog of the 
Multiplication Theorem holds for orthogonal independence. An example of the orthogonal loop product 
of graphs is given in Fig.3 (cf. Figs.1-2).
\begin{figure}
\unitlength=1mm
\special{em.linewidth 2pt}
\linethickness{0.5pt}
\begin{picture}(60.00,30.00)(20.00,0.00)
\put(0.00,10.00){\line(1,0){20.00}}
\put(35.00,10.00){\line(0,1){6.00}}

\put(0.00,10.00){\circle*{1.60}}
\put(0.00,8.00){\circle{4.00}}
\put(10.00,10.00){\circle*{1.00}}
\put(20.00,10.00){\circle*{1.00}}

\put(35.00,10.00){\circle*{1.60}}
\put(35.00,6.50){\circle{7.00}}
\put(35.00,16.00){\circle*{1.00}}
\put(30.00,20.00){\circle*{1.00}}
\put(40.00,20.00){\circle*{1.00}}

\put(80.00,22.00){\circle{4.00}}
\put(85.00,18.00){\circle{4.00}}
\put(90.00,22.00){\circle{4.00}}

\put(95.00,22.00){\circle{4.00}}
\put(100.00,18.00){\circle{4.00}}
\put(105.00,22.00){\circle{4.00}}

\put(1.00,11.00){\scriptsize $e_{1}$}
\put(71.00,11.00){\scriptsize $e$}
\put(10.00,11.00){\scriptsize $x$}
\put(20.00,11.00){\scriptsize $x'$}
\put(37.00,10.00){\scriptsize $e_{2}$}
\put(36.00,13.50){\scriptsize $y$}
\put(39.00,22.00){\scriptsize $y''$}
\put(29.00,22.00){\scriptsize $y'$}
\put(9.00,15.00){\scriptsize ${\mathcal G}_{1}$}
\put(27.00,15.00){\scriptsize ${\mathcal G}_{2}$}
\put(55.00,15.00){\scriptsize ${\mathcal G}_{1}\vdash_{\ell} {\mathcal G}_{2}$}





\put(70.00,10.00){\line(1,0){30.00}}
\put(85.00,10.00){\line(0,1){6.00}}
\put(100.00,10.00){\line(0,1){6.00}}

\put(70.00,8.00){\circle{4.00}}
\put(70.00,6.50){\circle{7.00}}
\put(85.00,6.50){\circle{7.00}}
\put(100.00,6.50){\circle{7.00}}

\qbezier(80,20)(80,16)(85,16)
\qbezier(85,16)(90,16)(90,20)
\qbezier(95,20)(95,16)(100,16)
\qbezier(100,16)(105,16)(105,20)

\qbezier(30,20)(30,16)(35,16)
\qbezier(35,16)(40,16)(40,20)

\put(70.00,10.00){\circle*{1.60}}
\put(85.00,10.00){\circle*{1.00}}
\put(100.00,10.00){\circle*{1.00}}

\put(85.00,16.00){\circle*{1.00}}
\put(100.00,16.00){\circle*{1.00}}
\put(90.00,20.00){\circle*{1.00}}
\put(80.00,20.00){\circle*{1.00}}
\put(105.00,20.00){\circle*{1.00}}
\put(95.00,20.00){\circle*{1.00}}

\end{picture}
\caption{An example of the orthogonal loop product}
\end{figure}
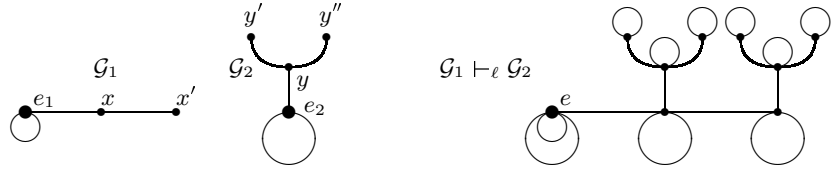
\begin{Example}
{\rm In the orthogonal loop product of graphs given in Fig.3, we have 
$$
(\beta_{1},w_{1},\beta_{2},w_{2}, \beta_{3}, w_{3}, \beta_{4}, w_{4})\in D_{8}(e),
$$
where the sequence of edges of color 1 is the same as in Example 3.1, namely
$$
(\beta_{1}, \beta_{2}, \beta_{3}, \beta_{4})=((e,x),(x,x'), (x',x), (x,e)),
$$
and the (alternating) f-walks attached to the vertices $x,x',x,e$ (in that order)
are: $w_{1}=w_{3}=(x,x)_{2}$, $w_{2}=(x',x')_{2}$ and $w_{4}=(e,e)$
Note that this d-walk has lenght $8$.
Again, as in the case of d-walks on the comb loop product, 
if any of the loops at $x$ or $x'$ is replaced by an alternating 
f-walk, we shall still get a rooted alternating d-walk.  
Of course, the simplest rooted alternating d-walk is given by
$(\beta_{1},w_{1},\beta_{4},w_{4})\in D_{4}(e)$ and $(\beta_{0},w_{4})\in D_{2}(e)$,
where $\beta_{0}=(e,e)_{1}$. }
\end{Example}
\begin{Theorem}
The Multiplication Theorem holds for orthogonal independence, the associated
loop product ${\mathcal G}_{1}\vdash_{\ell} {\mathcal G}_{2}$, and the
multiplicative convolution $\mu_1\angle \mu_2$.
\end{Theorem}
{\it Proof.}
The proof is similar to that of Theorem 3.1, but the combinatorics is based on 
a different formula for formal power series, namely 
$$
\eta_{\mu_1 \angle \mu_2}(z)=\frac{z\eta_{\mu_1}(\eta_{\mu_2}(z))}{\eta_{\mu_2}(z)},
$$ 
which leads to the combinatorial formula
$$
N_{\mu_1 \angle \mu_2}(n)=
\sum_{r=1}^{n}
N_{\mu_1}(r)\sum _{k_{1}+k_{2}+\ldots +k_{r-1}=n-1}
N_{\mu_2}(k_{1})N_{\mu_2}(k_{2})\ldots
N_{\mu_2}(k_{r-1}),
$$
where it is assumed that $k_{1}, k_{2}, \ldots , k_{r-1}$ are positive integers.
Recall that in the orthogonal
product ${\mathcal G}_{1}\vdash {\mathcal G}_{2}$ there is one copy of ${\mathcal G}_{1}$ 
(with $e_{1}$ identified with the root
$e$ of the product graph) with a copy of ${\mathcal G}_{2}$ attached by its root to every vertex $x$ of ${\mathcal G}_{1}$
but the root $e_{1}$. 
Therefore, each d-walk $w\in D_{2n}(e)$ is a sequence of $r$ edges of
${\mathcal G}_{1}$, which themselves must form an f-walk $c=(e,v_{1}, v_{2}, \ldots , v_{r-1},e)$ of color 1, 
interlaced with alternating f-walks $w_{i}\in F(v_{i})$, $1\leq i \leq r-1$ and a loop 
of color $2$, $w_{r}=(e,e)_{2}$. Note that this is the only way to produce a rooted 
alternating d-walk since the only edge of color $2$ incident on $e$ is the loop and therefore, in
order to get an alternating walk, the first f-walk which begins and ends with an edge of color $1$
must be followed by the loop at $e$ of color $2$ to make a double return to $e$.  
The contribution from each product of type 
$$
N_{\mu_1}(r)N_{\mu_2}(k_{1})N_{\mu_2}(k_{2})\ldots
N_{\mu_2}(k_{r-1})
$$ 
to the RHS of the above formula 
is equal to the number of all such f-walks $c$ of color 1 which have $r$ edges and are interlaced with $r-1$ f-walks 
of lenghts $2k_{1}-1, 2k_{2}-1, \ldots , 2k_{r-1}-1$ attached to all vertices of $c$ but one (we choose this vertex
to be the root since we want the considered walk to be an f-walk). 
The summation over all $1\leq r \leq n$ and $2k_{1}+ 2k_{2}+\ldots + 2k_{r-1}-r+1=2n-1$ indicates that
the first f-walk $u_1$ in the d-walk $w=(u_1,u_2)$ is of lenght $2n-1$, and $u_2$ is the loop of lenght $1$.
The summation over $1\leq r \leq n$ gives exactly the cardinality of $D_{2n}(e)$, which finishes the proof.\hfill $\blacksquare$
\begin{Example}
{\rm Let us apply the formulas of Example 4.1 to the enumeration of rooted alternating d-walks
on the graph ${\mathcal G}_{1}\vdash_{\ell} {\mathcal G}_{2}$ in Fig.3 (we keep the notation 
of Example 3.2 for spectral distributions $\mu_1$ and $\mu_2$ and thus we get the same values
of the $N_{\mu_1}(k)$'s and the $N_{\mu_2}(j)$'s). Using Example 4.1 and Theorem 6.1, we get
$D_{2}(e)=D_{4}(e)=D_{6}(e)=D_{8}(e)=1$, which can be verified directly by counting 
rooted alternating d-walks on ${\mathcal G}_{1}\vdash_{\ell} {\mathcal G}_{2}$.}
\end{Example}

\section{Subordination operators}

In analogy to the additive case [16], where we introduced and studied operators  
related to the subordination property for the free additive convolution, 
we shall now present an analogous approach to
the subordination property for the free multiplicative convolution.
We will mainly refer to sets ${\mathcal M}_{{\mathbb R}_{+}}$ and ${\mathcal M}_{{\mathbb T}}$,
which correspond to positive and unitary operators, 
but the operatorial subordination will also hold for all bounded operators.

Let $({\mathcal H}_{\iota}, \xi_{\iota})$, where $\iota \in I$, 
be Hilbert spaces with distinguished unit vectors.
Then their Hilbert space free product $({\mathcal H}_{1}, \xi_{1})*({\mathcal H}_{2},\xi_{2})$
is $({\mathcal H},\xi)$ where
\begin{equation}\tag{7.1}
{\mathcal H}={\mathbb C}\;\xi\oplus\bigoplus_{n=1}^{\infty}
\bigoplus_{\iota_1 \ne \iota_2 \ne ... \ne \iota_n}
{\mathcal H}_{\iota_1}^0\otimes {\mathcal H}_{\iota_2}^0\otimes \ldots \otimes {\mathcal H}_{\iota_n}^0,
\end{equation}
with ${\mathcal H}_\iota^0={\mathcal H}_{\iota} \ominus {\mathbb C} \xi_{\iota}$ and $\xi$ denoting a unit vector
(canonical scalar product is used).
For any $h\in{\mathcal H}_{\iota}$, denote by $h^0$ the orthogonal projection of $h$ onto ${\mathcal H}_{\iota}^0$.
Moreover, let
\begin{equation}\tag{7.2}
{\mathcal H}^{(n)}(\iota)=
\bigoplus_{\stackrel {\iota_{1}\neq \iota_{2}\neq \ldots \neq \iota_{n}}
{\scriptscriptstyle \iota_{1}\neq j}}
{\mathcal H}_{\iota_1}^0\otimes {\mathcal H}_{\iota_2}^0\otimes \ldots \otimes {\mathcal H}_{\iota_n}^0,
\end{equation}
\begin{equation}\tag{7.3}
{\mathcal K}^{(n)}(\iota)=
\bigoplus_{\stackrel {\iota_{1}\neq \iota_{2}\neq \ldots \neq \iota_{n}}
{\scriptscriptstyle \iota_{n}\neq \iota}}
{\mathcal H}_{\iota_1}^0\otimes {\mathcal H}_{\iota_2}^0\otimes \ldots \otimes {\mathcal H}_{\iota_n}^0,
\end{equation}
for $\iota\in I$ and $n\in {\mathbb N}$, and, for convenience,
we set ${\mathcal H}^{(0)}(\iota)={\mathcal K}^{(0)}(\iota)={\mathbb C}\xi$ with
the canonical projection $P_{\xi}:{\mathcal H}\rightarrow {\mathbb C}\xi$.

Since our index set $I$ consists of two elements, the above 
notation gives identifications
\begin{equation}\tag{7.4}
{\mathcal H}^{(n)}(\iota)=
\left\{
\begin{array}{ll}
{\mathcal K}^{(n)}(\iota) & n \;{\rm odd},\\
{\mathcal K}^{(n)}(\overline{\iota}) & n \; {\rm even}.
\end{array}
\right.
\end{equation} 
Nevertheless, it is convenient to use both sequences, (7.2) and (7.3), as well
as direct sums
\begin{equation}\tag{7.5}
{\mathcal H}(\iota)=\bigoplus_{n=1}^{\infty}{\mathcal H}^{(n)}(\iota) \;\;\;{\rm and}\;\;\;
{\mathcal K}(\iota)=\bigoplus_{n=1}^{\infty}{\mathcal K}^{(n)}(\iota),
\end{equation}
where $\iota\in I$. 

As in our previous work, we decompose the free product of Hilbert spaces 
as (two different) orthogonal direct sums
\begin{equation}\tag{7.6}
{\mathcal H}=\bigoplus_{n=1}^{\infty}{\mathcal H}^{(n-1)}(\iota)\oplus {\mathcal H}^{(n)}(\bar{\iota}),
\end{equation}
where $\iota\in I$, and denote by 
\begin{equation}\tag{7.7}
P_{\iota}(n):\;{\mathcal H}\rightarrow  {\mathcal H}^{(n-1)}(\iota)\oplus {\mathcal H}^{(n)}(\bar{\iota})
\end{equation}
the associated canonical projections.
Finally, we define the so-called {\it vacuum state} $\varphi(\cdot)=\langle \cdot\xi,\xi\rangle$
on $\mathcal{B}(\mathcal{H})$.

Of particular interest will be {\it s-free products of Hilbert spaces}, $({\mathcal H}_{\iota},\xi_{\iota})$
and $({\mathcal H}_{\overline{\iota}}, \xi_{\iota})$, denoted $({\mathcal H}_{\iota}\oright {\mathcal H}_{\overline{\iota}}, \xi)$
and defined [16] as the pair $({\mathcal K}_{\iota}, \xi)$
where
\begin{equation}\tag{7.8}
{\mathcal K}_{\iota}={\mathbb C}\xi\oplus {\mathcal K}(\overline{\iota}),
\end{equation}
and $\iota\in I$. The direct sum decompositions
\begin{equation}\tag{7.9}
{\mathcal K}_{\iota}=
\bigoplus_{n\;{\rm odd}}{\mathcal H}^{(n-1)}(\iota)\oplus {\mathcal H}^{(n)}(\bar{\iota})
\end{equation}
hold for each $\iota\in I$.

Let $x\in B({\mathcal H}_{1})$ and $y\in B({\mathcal H}_{2})$
be fixed random variables.
The corresponding free random variables $X_{1}:=\lambda(x)$ and 
$X_2:=\lambda(y)$ are elements of $B({\mathcal H})$, where 
$\lambda$ is the free product representation on ${\mathcal H}$, and can be 
decomposed according to (7.6) as 
\begin{equation}\tag{7.10}
X_{\iota}=\sum_{j=1}^{\infty} X_{\iota}(j),
\end{equation}
where $X_{\iota}(j) = P_{\iota}(j)X_{\iota}P_{\iota}(j)$, $j\in {\mathbb N}$, 
can be viewed as replicas of $x$ and $y$, respectively, where $\iota\in I$.
Using (7.6) and (7.9), we also have the decompositions
\begin{equation}\tag{7.11}
1=\sum_{j=1}^{\infty}P_{\iota}(j)\;\;\;{\rm and}\;\;\;
1_{{\mathcal K}_{\iota}}=\sum_{j \;{\rm odd}}P_{\iota}(j)
\end{equation}
of the units in $B({\mathcal H})$ 
and $B({\mathcal K}_{\iota})$, where $\iota\in I$, respectively.
Here, and also in the sequel, we denote by the same symbols, $X_{\iota}(j)$ and
$P_{\iota}(j)$, the corresponding operators on ${\mathcal K}_{\iota}$.
 
We shall use representations of free random variables
as `orthogonal series' with the unit singled out, namely
\begin{equation}\tag{7.12}
X_{\iota}=1+\sum_{j=1}^{\infty}x_{\iota}(j)
\end{equation}
where $x_{\iota}(j)=X_{\iota}(j)-P_{\iota}(j)$ and $\iota \in I$, $j\in {\mathbb N}$.
If no confusion arises, we will distinguish the first term in the above series
by a special notation and write $x_{\iota}=x_{\iota}(1)$. Let us remark that the above form is 
suitable for the study of multiplicative convolutions, where units play a special role
(roughly speaking, they have to be subtracted when it comes to proving 
some kind of independence).

We begin with the decomposition of the product $X_{\overline{\iota}}X_{\iota}$
which corresponds to Eq.(1.3) for $\mu_1,\mu_2\in {\mathcal M}_{{\mathbb R}_{+}}$,
and therefore involves positive operators. 
\begin{Theorem}
Let $X_{\iota}$, $\iota \in I$, be positive random variables which are free
with respect to $\varphi$, and denote $t_{\iota}=1+x_{\iota}$,
$T_{\overline{\iota}}=\sqrt{X_{\iota}-x_{\iota}}X_{\overline{\iota}}\sqrt{X_{\iota}-x_{\iota}}$. 
Then $X_{\overline{\iota}}X_{\iota}$ has the same $\varphi$-distribution as 
$\sqrt{t_{\iota}}T_{\overline{\iota}}\sqrt{t_{\iota}}$ and the pair $(t_{\iota}-1,T_{\overline{\iota}}-1)$
is monotone independent with respect to $\varphi$.
\end{Theorem}
{\it Proof.}
Note that $T_{\overline{\iota}}$ introduced above is well-defined since 
$$
X_{\iota}-x_{\iota}=P_{\iota}(1)+\sum_{j=2}^{\infty}X_{\iota}(j)
$$
is positive as an orthogonal direct sum of positive random variables. 
Using a similar argument, we get positivity of 
$$
t_{\iota}=X_{\iota}(1)+\sum_{j=2}^{\infty}P_{\iota}(j).
$$
Computations of square roots give
$$
\sqrt{t_{\iota}}=\sqrt{X_{\iota}(1)}+\sum_{j=2}^{\infty}P_{\iota}(j) \;\;\;{\rm and}\;\;\;
\sqrt{X_{\iota}-x_{\iota}}=P_{{\iota}}(1)+\sum_{j=2}^{\infty}\sqrt{X_{\iota}(j)},
$$
where $\sqrt{X_{\iota}(j)}=P_{\iota}(j)\lambda(\sqrt{z})P_{\iota}(j)$ for $j\in {\mathbb N}$,
with $z=x$ or $z=y$ if $\iota=1$ or $\iota =2$, respectively, which
leads to
$$
\sqrt{t_{\iota}}\sqrt{X_{\iota}-x_{\iota}}=\sqrt{X_{\iota}}.
$$
Now, we use the fact that $X_{\overline{\iota}}X_{\iota}$ has the same $\varphi$-distribution
as $\sqrt{X_{\iota}}X_{\overline\iota}\sqrt{X_{\iota}}$, which completes the proof of the first part of the 
theorem.

We will show now that the pair $(t_{\iota}-1,T_{\overline{\iota}}-1)$ is 
monotone independent w.r.t. $\varphi$. For simplicity, we denote
$x=t_{\iota}-1$ and $y=T_{\overline{\iota}}-1$.
\indent{\par}
{\it Case 1.}
If $y^{n}$ is in the middle of the moment, we compute
$\varphi(w_{1}xy^{n}xw_{2})$, 
where $n\in {\mathbb N}$ and $w_{1},w_{2}\in {\rm alg}(x,y)$. Note that 
the range of $x$ is ${\mathbb C}\oplus {\mathcal H}_{\iota}^{0}$, therefore, we only need to find
the action of $y^{n}$ onto $\xi$ and $h\in {\mathcal H}_{{\iota}}^{0}$.
Using the explicit form of $\sqrt{X_{\iota}-x_{\iota}}$, we get
\begin{eqnarray*}
y^{n}\xi &=& (
(1+\sum_{k=2}^{\infty}r_k)(1+\sum_{j=1}^{\infty}z_j)(1+\sum_{k=2}^{\infty}r_k)-1
)^{n}\xi\\
&=&
\varphi(y^{n})\xi \;\;\;\;\;\; 
{\rm mod}\;{\mathcal K}(\iota),
\end{eqnarray*}
where $r_{k}=\sqrt{X_{\iota}(k)}-P_{\iota}(k)$ and $z_j=X_{\overline{\iota}}(j)-P_{\overline{\iota}}(j)$ 
and therefore $r_{k}\xi=0$ for any $k\geq 2$. Similarly, we get
$$ 
y^{n}h = 
\langle y^{n}h,h\rangle h \;\;\;\;{\rm mod}\;{\mathcal K}(\overline{\iota})\ominus {\mathcal K}^{(1)}(\overline{\iota})
$$
for any $h\in {\mathcal H}_{{\iota}}^{0}$ of norm $\parallel h \parallel =1$.
Note that $h$ plays the role of a cyclic vector for the restriction of the free product representation $\pi_{1}*\pi_{2}$ to the unital algebra 
generated by $y$. Thus $\langle y^nh,h\rangle 
=\varphi(y^n)h$, which gives
$$
\varphi(w_{1}xy^{n}xw_{2})=
\varphi(y^{n})\varphi(w_{1}x^{2}w_{2}),
$$
i.e. the required condition for monotone independence.
\indent{\par}
{\it Case 2.}
If $y^{n}$ is at the end of the moment, we get
$$
\varphi(w_{1}xy^{n})=\varphi(y^{n})\varphi(w_{1}x),
$$
using the relation for $y^n \xi$ obtained above and the fact that 
${\mathcal K}(\iota)\subset {\rm Ker}x$,
which finishes the proof. \hfill $\blacksquare$
\begin{Corollary}
If, in addition, $X_{\iota}$ has $\varphi$-distribution $\mu_{\iota}\neq \delta_{0}$ for $\iota\in I$,
then the $\varphi$--distribution of $T_{\iota}$ is given by the s-free multiplicative convolution 
$\mu_{\iota}\boxslash \mu_{\overline{\iota}}$.
\end{Corollary}
{\it Proof.}
First, note that the variable $t_{\iota}$ has $\varphi$--distribution $\mu_{\iota}$. 
In view of (1.5) and Theorem 7.1, we have
$\eta_{\mu_{1}\boxtimes\, \mu_{2}}=\eta_{\mu_{1}}\circ \eta_{\sigma_2}=
\eta_{\mu_{2}}\circ \eta_{\sigma_1}$, 
where $\sigma_{1}$ and $\sigma_{2}$ are the $\varphi$--distributions of 
$T_{1}$ and $T_{2}$, respectively. From the definition of the s-free multiplicative convolution, 
it follows that $\sigma_{1}=\mu_{1}\boxslash \mu_{2}$ and $\sigma_{2}=\mu_2\boxslash \mu_1$,
which proves our assertion.
\hfill $\blacksquare$
\begin{Remark}
{\rm In analogy to the terminology used in analytic subordination [8,23], we can say 
that $X_{\overline{\iota}}X_{\iota}$ is `subordinate to $t_{\iota}$ and $t_{\overline{\iota}}$, with
$T_{\overline{\iota}}$ and $T_{\iota}$, respectively, being the corresponding `subordination operators'.
Let us remark, however, that operators $T_{\iota}$, $\iota\in I$, are not 
good candidates for `subordination branches' of the product of free random variables
since they are not suitable for further decomposition. In that sense,
they are not multiplicative analogues of the 
`additive subordination branches' [16], or the branches
of the free product of graphs introduced by Quenell [20].
In the next section we will introduce and study such analogues, 
whose decompositions will correspond to decompositions of s-free multiplicative 
convolutions.}
\end{Remark}

A similar operatorial subordination result can be established for bounded operators
which includes unitary operators (the latter are related to ${\mathcal M}_{{\mathbb T}}$).
\begin{Theorem}
Let $X_{\iota}$, $\iota \in I$, be bounded random variables which are free
with respect to $\varphi$. Then, for each $\iota \in I$, the operator
$X_{\overline{\iota}}X_{\iota}$ has the same $\varphi$--distribution as 
$Q_{\overline{\iota}}q_{\iota}$, where $q_{\iota}=1+x_{\iota}$ and  
$Q_{\overline{\iota}}=X_{\overline{\iota}}(X_{\iota}-x_{\iota})$. 
Moreover, the pair $(q_{\iota}-1,Q_{\overline{\iota}}-1)$ is monotone 
independent with respect to $\varphi$. Finally, if $X_{\iota}$, $\iota\in I$, are unitary, then
the operators $q_{\iota},Q_{\iota}$, $\iota \in I$, are unitary.
\end{Theorem}
{\it Proof.}
The proof of the statements concerning bounded operators is similar
to that of Theorem 7.1. Therefore, we shall just prove unitarity of
$Q_{\iota}$ and $q_{\iota}$.  
Writing operators $q_{\iota}$ and $X_{\iota}-x_{\iota}$
in the form of orthogonal series
$$
q_{\iota}=X_{\iota}(1)+\sum_{j=2}^{\infty}P_{\iota}(j), \;\;\; 
X_{\iota}-x_{\iota}=P_{\iota}(1)+\sum_{j=2}^{\infty}X_{\iota}(j),
$$
we obtain
$$
q_{\iota}q_{\iota}^{*}=X_{\iota}(1)X_{\iota}^{*}(1)+\sum_{j=2}^{\infty}P_{\iota}(j)=1
$$
and 
$$
(X_{\iota}-x_{\iota})(X_{\iota}-x_{\iota})^{*}=P_{\iota}(1)+\sum_{j=2}^{\infty}X_{\iota}(j)X_{\iota}^{*}(j)=1.
$$
Similarly, $q_{\iota}^{*}q_{\iota}=1$ and $(X_{\iota}-x_{\iota})^{*}(X_{\iota}-x_{\iota})=1$.
This proves unitarity of $q_{\iota}$ and $X_{\iota}-x_{\iota}$,
from which we obtain unitarity of $Q_{\iota}$. 
The remaining part of the proof is similar to that of Theorem 7.1.\hfill $\blacksquare$
\begin{Corollary}
Let $\mu_{\iota}$ and $\sigma_{\iota}$ be the $\varphi$--distributions of $X_{\iota}$ and
$Q_{\iota}$, respectively, for $\iota \in I$. Then the $\varphi$--distribution of 
$X_{\overline{\iota}}X_{\iota}$ is given by 
$\mu_\iota\circlearrowright \sigma_{\overline{\iota}}$ for each $\iota \in I$.
In particular, if the operators $X_{\iota}$ are unitary, then
$\sigma_{\iota}=\mu_{\iota}\boxslash \mu_{\overline{\iota}}$ for each
$\iota \in I$.
\end{Corollary}
{\it Proof.}
The first statement follows from Theorem 7.3 and the definition of 
the monotone multiplicative convolution of distributions.
The proof of the second statement is similar to that of Corollary 7.2
(the s-free convolution of measures from ${\mathcal M}_{*}$ is used).
\hfill $\blacksquare$

\section{Subordination branches}

In order to introduce operator-valued `multiplicative subordination branches' 
associated with the product of bounded free random variables, let us introduce 
random variables 
\begin{equation}\tag{8.1}
R_{\iota}(m)=
1+\sum_{j\in J(m)}x_{\iota}(j)
\end{equation}
where the notation $J(m)=\{m,m+2,m+4,\ldots \}$ is used for index sets
and, by abuse of notation, $1$ denotes $1_{{\mathcal K}_{\iota}}$ if $m$ is odd,  
or $1_{{\mathcal K}_{\overline{\iota}}}$ if $m$ is even.
Note that operators $R_{\iota}(m)$ are elements of $B({\mathcal K}_{\iota})$ or $B({\mathcal K}_{\overline{\iota}})$,
depending on whether we have $m$ odd or even, respectively.

Therefore, on each of the `s-free Fock spaces', 
${\mathcal K}_{\iota}$, $\iota \in I$, we get an `interlaced'
sequence of operators
\begin{equation}\tag{8.2}
R_{\iota}(1),R_{\overline{\iota}}(2), R_{\iota}(3), R_{\overline{\iota}}(4), \ldots \in B({\mathcal K}_{\iota}).
\end{equation}
Each of these two sequences is used in the definition of 
{\it one} sequence of operator-valued `multiplicative subordination branches'.
\begin{Proposition}
If $X_1,X_2$ are positive (unitary), then operators $R_{\iota}(m)$, where $\iota \in I$ and $m\in {\mathbb N}$, 
are positive (unitary).
\end{Proposition}
{\it Proof.}
This fact can be easily checked by decomposing the units according to (7.11).
We choose to show unitarity of 
$$
R_{\iota}^{\scriptscriptstyle {\rm odd}}=
\sum_{\scriptscriptstyle j \; {\rm odd}}X_{\iota}(j)\;\;\; {\rm and}\;\;\; 
R_{\iota}^{\scriptscriptstyle {\rm even}}=P_{\xi}+\sum_{\scriptscriptstyle j \; {\rm even}}X_{\iota}(j)
$$ 
for unitary $X_{\iota}$. Since $X_{\iota}(j)X_{\iota}^{*}(j)=P_{\iota}(j)X_{\iota}X_{\iota}^{*}P_{\iota}(j)=P_{\iota}(j)$
for each $\iota \in I$ and $j\in {\mathbb N}$, we get
$$
R_{\iota}^{\scriptscriptstyle {\rm odd}}(R_{\iota}^{\scriptscriptstyle {\rm odd}})^{*}=
\sum_{\scriptscriptstyle j \;{\rm odd}}X_{\iota}(j)X_{\iota}^{*}(j)
=
\sum_{\scriptscriptstyle j \;{\rm odd}}P_{\iota}(j)
=
1,
$$
$$
R_{\iota}^{\scriptscriptstyle {\rm even}}(R_{\iota}^{\scriptscriptstyle {\rm even}})^{*}=
P_{\xi}+\sum_{\scriptscriptstyle j \;{\rm even}}X_{\iota}(j)X_{\iota}^{*}(j)
=
P_{\xi}+\sum_{\scriptscriptstyle j \;{\rm even}}P_{\iota}(j)
=
1.
$$
Unitarity of $R_{\iota}(m)$ for arbitrary $m$ can be proved in similar way.
Positivity of $R_{\iota}(m)$ follows easily from the appropriate 
decomposition of the unit.
\hfill $\blacksquare$
\begin{Definition}
{\rm Let $X_{\iota}$, $\iota \in I$, be random variables which are
free w.r.t. $\varphi$. Random variables
\begin{equation}\tag{8.3}
U_{\iota}(m)=R_{\iota}(m)R_{\overline{\iota}}(m+1),
\end{equation}
where $\iota\in I$ and $m\in {\mathbb N}$, will be called 
{\it multiplicative bounded subordination branches}.
If $X_{\iota}$, $\iota\in I$ are unitary, then they 
will be called {\it multiplicative unitary subordination branches}.} 
\end{Definition}

In addition to branches $U_{\iota}(m)$, we shall need operators obtained 
from $U_{\iota}(m)$ by `order reversal'. Namely, let
\begin{equation}\tag{8.4}
V_{\iota}(m)= R_{\overline{\iota}}(m+1)R_{\iota}(m)
\end{equation}
where $\iota \in I$ and $m\in {\mathbb N}$. Again, we use a simpler notation for $m=1$, 
namely $V_{\iota}=V_{\iota}(m)$.
\begin{Definition}
{\rm Let $X_{\iota}$, $\iota \in I$ be positive random variables
which are free w.r.t. $\varphi$. Random variables
\begin{equation}\tag{8.5}
Y_{\iota}(m) = 
\sqrt{R_{\overline{\iota}}(m+1)}
R_{\iota}(m)
\sqrt{R_{\overline{\iota}}(m+1)},
\end{equation}
where $\iota \in I$ and $m\in {\mathbb N}$, will be called 
{\it multiplicative positive subordination branches}.}
\end{Definition}

For simplicity, we will also use the terms: 
`bounded branches', `unitary branches' and `positive branches'. 
As in the case of operators $R_{\iota}(m)$, where $m\in {\mathbb N}$ and  $\iota\in I$, 
we get two sequences of alternating (bounded, unitary) branches
\begin{equation}\tag{8.6}
U_{\iota}(1),U_{\overline{\iota}}(2), U_{\iota}(3), U_{\overline{\iota}}(4), \ldots \; \in B({\mathcal K}_{\iota})
\end{equation}
for each $\iota\in I$, and similar sequences of positive branches. 
Recall that existence of `interlaced branches' was also observed for the 
additive (self-adjoint) subordination branches [16].

Of particular interest are branches of 1st order, 
and it is therefore of advantage to use a special notation for the first 
two operators in each `interlaced' sequence given by (8.1). 
For $\iota\in I$, we shall use
\begin{equation}\tag{8.7}
R_{\iota}^{\scriptscriptstyle {\rm odd}}  =
1+\sum_{\scriptscriptstyle j\;{\rm odd}}x_{\iota}(j)\;\;\;{\rm and}\;\;\;
R_{\iota}^{\scriptscriptstyle {\rm even}}  =
1+\sum_{\scriptscriptstyle j\;{\rm even}}x_{\iota}(j)
\end{equation}
where again, the unit $1$ has to be interpreted as $1_{{\mathcal K}_{\iota}}$ and
$1_{{\mathcal K}_{\overline{\iota}}}$, respectively.
Note that although we have two different units here, 
the operators which appear in the same `interlaced' sequence
contain the same unit.  
\begin{Remark}
{\rm Let us recall that in the additive case [15, Theorem 8.4], the branches $B_1$ and $B_2$ of the 
sum $X_1+X_2$ can be decomposed as 
\begin{equation}\tag{8.8}
B_{1}= S_{1}^{\scriptscriptstyle {\rm odd}}+S_{2}^{\scriptscriptstyle {\rm even}},
\end{equation}
where 
$S_{1}^{\scriptscriptstyle {\rm odd}}=\sum_{j \;{\rm odd}}X_{1}(j)$, 
$S_{2}^{\scriptscriptstyle {\rm even}}=\sum_{j\;{\rm even}}X_{2}(j)$,
and that the pair $(S_{1}^{\scriptscriptstyle {\rm odd}},S_{2}^{\scriptscriptstyle {\rm even}})$ 
is s-free w.r.t. the pair of states $(\varphi, \psi)$,
where $\psi$ is associated with any unit vector $\zeta\in {\mathcal H}_{1}^{0}$.
Note that we have 
$R_{1}^{\scriptscriptstyle {\rm odd}}=S_{1}^{\scriptscriptstyle {\rm odd}}$ 
and $R_{2}^{\scriptscriptstyle {\rm even}}=P_{\xi}+S_{2}^{\scriptscriptstyle {\rm even}}$,
but we prefer to have a new notation to maintain a uniform style 
for `subtracting units'.}
\end{Remark}

This, in turn, leads to a special notation for branches of 1st order, 
$U_{\iota}=U_{\iota}(1)$, $Y_{\iota}=Y_{\iota}(1)$, $\iota\in I$.
For instance,
\begin{eqnarray*}
Y_{1}&=&Y_{1}(1)=\sqrt{R_{2}^{{\scriptscriptstyle {\rm even}}  }}
R_{1}^{{\scriptscriptstyle {\rm odd}}  }\sqrt{R_{2}^{{\scriptscriptstyle {\rm even}}  }},\\
Y_{2}&=&Y_{2}(1)=\sqrt{R_{1}^{{\scriptscriptstyle {\rm even}}  }}
R_{2}^{{\scriptscriptstyle {\rm odd}}  }\sqrt{R_{1}^{{\scriptscriptstyle {\rm even}}  }}.
\end{eqnarray*}
Many computations can be reduced to branches of 1st order since their $\varphi$-distributions 
agree with the distributions of branches of higher orders with respect to 
suitably chosen states.

Let us examine the $\varphi$--distributions of branches of first order.
In fact, we will show that $Y_{\iota}$ and $U_{\iota}$ have the same $\varphi$-distributions
as $T_{\iota}$ for given $\iota \in I$.
\begin{Lemma}
The variables $Q_{\iota}$, $U_{\iota}$ and 
$V_{\iota}$ have the same $\varphi$--distributions for any given $\iota \in I$.
If $X_{\iota}$, $\iota\in I$,
are positive, then the variables $T_{\iota}, Y_{\iota}, U_{\iota}$
and $V_{\iota}$ have the same $\varphi$-distributions for any given $\iota \in I$.
\end{Lemma}
{\it Proof.}
We shall give the proof in the more difficult case of 
positive operators (the proof for the case of bounded operators is similar).
For $n\geq 1$, we have
\begin{eqnarray*}
\varphi(Y_{\iota}^{n})&=&
\left(
\sqrt{R_{\overline{\iota}}^{{\scriptscriptstyle {\rm even}}  }}
(R_{\iota}^{{\scriptscriptstyle {\rm odd}}  }R_{\overline{\iota}}^{{\scriptscriptstyle {\rm even}}  })^{n-1}
R_{\iota}^{{\scriptscriptstyle {\rm odd}}  }
\sqrt{R_{\overline{\iota}}^{{\scriptscriptstyle {\rm even}}  }}
\right)
\\
&=&
\varphi
\left(
(R_{\iota}^{{\scriptscriptstyle {\rm odd}}  }R_{\overline{\iota}}^{{\scriptscriptstyle {\rm even}}  })^{n-1}R_{\iota}^{{\scriptscriptstyle {\rm odd}}  }
\right)\\
&=&
\varphi
\left(
U_{\iota}^{n}
\right)
\end{eqnarray*}
since $\sqrt{R_{\overline{\iota}}^{{\scriptscriptstyle {\rm even}}  }}\xi=\xi$ (the only summand in 
the direct sum decomposition of this square root which gives a non-zero contribution 
is $1_{\overline{\iota}}(1)$). This shows that $Y_{\iota}$ and 
$U_{\iota}$ have the same $\varphi$-distributions.

Let us now compare the moments of $T_{\iota}$ with those of $Y_{\iota}$. We have
\begin{eqnarray*}
\varphi(T_{\iota}^{n})&=&
\varphi
\left(
\sqrt{X_{\overline{\iota}}-x_{\overline{\iota}}}(X_{\iota}(X_{\overline{\iota}}-x_{\overline{\iota}}))^{n-1}X_{\iota}\sqrt{X_{\overline{\iota}}-x_{\overline{\iota}}}
\right)\\
&=&
\varphi
\left(
(X_{\iota}(X_{\overline{\iota}}-x_{\overline{\iota}}))^{n}
\right)
\end{eqnarray*} 
since $\sqrt{X_{\overline{\iota}}-x_{\overline{\iota}}}\xi=\xi$.
To fix attention, suppose that $\iota=2$ and denote, for convenience,
$z_{k}=x_{1}(k)$ and $y_{k}=x_2(k)$.
Now, let us observe that in the finite set of
variables of type $z_k$ and $y_j$, where $k\geq 2$ and $j\geq 1$, 
which are used when calculating 
$X_{2}(X_{1}-x_1)^{n-1}X_{2}\xi$ 
there are no $z_{k}$'s for odd $k$'s, or 
$y_{j}$'s for even $j$'s. 
This is because the only variable which may give a non-zero contribution 
when acting on $\xi$ is $y_{1}$, then only $y_1$ and $z_2$ may follow, 
of which the first operator can be followed by $y_{1}$ or $z_2$, whereas the second --
by $z_2,y_1$ or $y_3$, etc. This leads to 
$$
\varphi(P(X_{2},X_{1}-x_{1}))=\varphi 
(P(R_{2}^{{\scriptscriptstyle {\rm odd}}  }, R_{1}^{{\scriptscriptstyle {\rm even}}  }))
$$
for any polynomial $P$ in two noncommuting variables. This shows that $T_{2}$ and $Y_{2}$ have the same 
$\varphi$-distributions.

Finally, we will show that $U_{2}$ and $V_{2}$ have the same 
$\varphi$-distributions.
Using the notation used in the previous paragraph, we have 
\begin{eqnarray*}
\varphi \left((U_{2}-1)^{n}\right)
&=&
\varphi
(
(\sum_{{\scriptscriptstyle k\;{\rm even}}}z_{k}
+
\sum_{{\scriptscriptstyle j\; {\rm odd}}}y_{j}
+
\sum_{\stackrel{k\;{\rm even},\; j\;{\rm odd}}
{\scriptscriptstyle |j-k|=1}}
y_{j}z_{k}
)^{n}
)
\\
\varphi \left((V_{2}-1)^{n}\right)
&=&
\varphi
(
(\sum_{{\scriptscriptstyle k\;{\rm even}}}z_{k}
+
\sum_{{\scriptscriptstyle j\; {\rm odd}}}y_{j}
+
\sum_{\stackrel{k\;{\rm even},\; j\;{\rm odd}}
{\scriptscriptstyle |j-k|=1}}
z_{k}y_{j}
)^{n}
)
\end{eqnarray*}
for any natural $n$, and therefore, we need to show that the right-hand-sides of the above
equations are equal to each other.
Note that the products of variables $z_{k}$ and $y_{j}$ in words $w$ which appear in 
mixed moments $\varphi(w)$ giving non-zero contributions to the above moments satisfy the following
three conditions:
\indent{\par}
(i) whenever $z_{k}$ stands next to $z_{m}$ or $y_{k}$ next to $y_{m}$, 
it holds that $k=m$, 
\indent{\par}
(ii) whenever $y_{j}$ stands next to $z_{k}$, it holds that  $|j-k|=1$,
\indent{\par}
(iii) the word $w$ begins and ends with the letter $y_{1}$.\\ 
Let us denote by $W$ the set of words in letters  
$$
L=\{z_{k},y_{j}:\,k\;{\rm even},\;j\; {\rm odd}\}
$$
subject to conditions (i)-(iii).
In turn, let $W'$ and $W''$, respectively, denote the sets of all words in letters 
\begin{eqnarray*}
L'&=&\{z_{k},y_{j}, z_{k}y_{k-1}, z_{k}y_{k+1}:\; k\;{\rm even}, \; j\;{\rm odd}\},\\
L''&=&\{z_{k},y_{j}, y_{k-1}z_{k}, y_{k+1}z_{k}:\; k\;{\rm even}, \; j\;{\rm odd}\}
\end{eqnarray*}
subject to conditions (i) and (iii), where products of type $z_{k}y_{j}$ and $y_{j}z_{k}$
are treated as letters. When they are used to form words, they are denoted by 
$[z_{k}y_{j}]$ and $[y_{j}z_{k}]$, respectively. Note that $W\subset W'$ and $W\subset W''$.
More importantly, there is a bijection 
$$
\tau: W'\rightarrow W''
$$
defined by `shifting brackets to the left'. Namely, $\tau$ is uniquely defined
by $\tau(w_{1})=w_{1}$ and the recursion
$$
\tau(w_{1}y_{m}z_{k}^{p}[z_{k}y_{j}]w_{2})=w_{1}[y_{m}z_{k}]z_{k}^{p}\tau(y_{j}w_{2})
$$
for any $w_{1}\in W$, $w_{2}\in W'$, $p\geq 0$ and $|j-k|=1$, $|k-m|=1$, i.e.
we assume that $[z_ky_j]$ is the first such pair (counting from the left).
In this recursion, we shift the brackets to the left from this pair
to the closest possible pair $y_{m},z_{k}$ lying to its left 
(in view of (iii), such a pair must exist). It is not hard to see that
its inverse is given by $\tau^{-1}(u_{2})=u_{2}$ and the recursion
$$
\tau^{-1}(w_{1}'[y_{j}z_{k}]z_{k}^{p}y_{m}w_{2}')=\tau^{-1}(w_{1}'y_{j})z_{k}^{p}[z_{k}y_{m}]w_{2}'
$$
where $w_{1}'\in W''$, $w_{2}'\in W$, $p\geq 0$ and $|j-k|=1$, $|k-m|=1$. 
Again, note that for every such pair $[y_{j}z_{k}]$ there must exist $y_{m}$ 
standing to its right in view of (iii).
Clearly, for every $w\in W'$ it holds that $\varphi(\tau(w))=\varphi(w)$. 
This completes the proof. \hfill $\blacksquare$\\
\indent{\par}
We are ready to prove a theorem, which will lead to a relation
between the distributions of two consecutive branches. 
For that purpose, introduce operators
\begin{equation}\tag{8.9}
Z_{\iota}(m)=
1+x_{\iota}(m)
\end{equation}
with the same convention concerning the units as before ($1=1_{{\mathcal K}_{\iota}}$ if $m$ is odd
and $1=1_{{\mathcal K}_{\overline{\iota}}}$ if $m$ is even), 
as well as operators obtained from the subordination branches by subtracting the (appropriate) unit:
\begin{equation}\tag{8.10}
y_{\iota}(m)=Y_{\iota}(m)-1 
\end{equation}
for any $m\in {\mathbb N}$ and $\iota \in I$.
\begin{Theorem}
For given $m\in {\mathbb N}$ and $\iota\in I$, let $\varphi, \psi$ be states associated with any unit 
vectors $\zeta\in {\mathcal H}^{(m-1)}(\iota)$ and $\zeta'\in {\mathcal H}^{(m)}(\overline{\iota})$.
Then, operators $U_{\iota}(m)$ and $Z_{\iota}(m)U_{\overline{\iota}}(m+1)$
have the same $\varphi$--distributions. 
Similarly, if $X_{\iota}$, $\iota\in I$, are positive, 
then operators $Y_{\iota}(m)$ and $Z_{\iota}(m)Y_{\overline{\iota}}(m+1)$ have the same $\varphi$--distributions. 
Moreover, the pairs $(x_{\iota}(m), u_{\overline{\iota}}(m+1))$ and 
$(x_{\iota}(m),y_{\overline{\iota}}(m+1))$ 
are orthogonal with respect to $(\varphi, \psi)$ for each $\iota \in I$ and $m\in {\mathbb N}$.
\end{Theorem}
{\it Proof.}
We shall prove the more difficult case of positive random variables. 
Without loss of generality, let $\iota=1$ and, 
for notational simplicity, consider the case $m=1$ (the proof of the 
general case is similar). Denote
$$
Z_{1}=Z_{1}(1), \;R_{1}=R_{1}(1), \; R_{3}=R_{1}(3), \; R_{2}=R_{2}(2), \; Y_{2}=Y_{2}(2).
$$
In the first part of the theorem we need to prove that
$Y_{1}$ has the same $\varphi$--distribution as $Z_{1}Y_{2}$.
Using orthogonal decompositions and `square root calculus', we obtain
$$
Z_{1}=X_{1}(1)+\sum_{j\in J(3)}P_{1}(j),\;\;\;{\rm and}\;\;\;
\sqrt{R_{3}}=P_{1}(1)+\sum_{j\in J(3)}\sqrt{X_{1}(j)}
$$
which gives
$$
\sqrt{R_{3}}Z_{1}\sqrt{R_{3}}=
\sum_{\scriptscriptstyle j\;{\rm odd}}X_{1}(j)=R_{1},
$$
and therefore
$$
Y_{1}=\sqrt{R_{2}}R_{1}\sqrt{R_{2}}
=
\sqrt{R_{2}}
\sqrt{R_{3}}Z_{1}\sqrt{R_{3}}
\sqrt{R_{2}}.
$$
Therefore,
\begin{eqnarray*}
\varphi(Y_{1}^{n})&=&\varphi\left(
(\sqrt{R_{2}}
\sqrt{R_{3}}Z_{1}\sqrt{R_{3}}
\sqrt{R_{2}})^{n}\right)\\
&=&
\varphi
\left(
\sqrt{R_2}\sqrt{R_{3}}(Z_{1}\sqrt{R_{3}}R_{2}\sqrt{R_{3}})^{n-1}Z_{1}\sqrt{R_{3}}\sqrt{R_{2}}
\right)\\
&=&\varphi
\left(
\sqrt{R_2}\sqrt{R_{3}}(Z_{1}Y_{2})^{n-1}Z_{1}\sqrt{R_{3}}\sqrt{R_{2}}
\right)\\
&=&
\varphi
\left(
(Z_{1}Y_{2})^{n-1}Z_{1}
\right)\\
&=&
\varphi
\left(
(Z_{1}Y_{2})^{n}
\right)
\end{eqnarray*}
since $Y_2\xi=\xi$ and also $\sqrt{R_3}\sqrt{R_2}\xi=\xi$.
Therefore, the $\varphi$--distribution of $Y_{1}$ agrees with the $\varphi$--distribution of
$Z_{1}Y_{2}$.

In the proof of orthogonality, 
denote, as before, $x_1=x_1(1)$, and $y_{2}=Y_{2}-1_{{\mathcal K}_{1}}$.
Using orthogonal decompositions and `square root calculus', we can write
\begin{eqnarray*}
Y_2
&=&
\sqrt{R_{3}}R_{2}\sqrt{R_{3}}\\
&=&
(P_{1}(1)+\sum_{j\in J(3)}\sqrt{X_{1}(j)})
(P_{\xi}+\sum_{k\in J(2)}X_{2}(k))
(P_{1}(1)+\sum_{j\in J(3)}\sqrt{X_{1}(j)})
\end{eqnarray*}
and therefore $Y_{2}\xi=\xi$, which implies that $y_{2}\xi=0$. Therefore,
we get $\varphi(wy_{2})=0$ and therefore, by taking the adjoints, $\varphi(y_{2}w)=0$
for any $w\in {\rm alg}(x_1,y_2)$. This gives the first orthogonality condition.
Now, let us prove the second orthogonality condition, i.e. 
$$
\varphi(w_{1}x_1y_{2}^{k}x_1w_{2})
=
\psi(y_{2}^{k})(\varphi(w_{1}x_1^{2}w_{2})-\varphi(w_{1}x_1)\varphi(x_1w_{2}))
$$
for any $k\in {\mathbb N}$ and $w_{1},w_{2}\in {\rm alg}(x_1,y_{2})$.
Note that the `lowest order term' in the expression for $y_2$, which turns out to be of special importance
when $x_1$ stands next to $y_2$, is of the form 
$$
P_{1}(1)X_{2}(2)P_{1}(1)=P_{1}(1)P_{2}(2)X_{2}P_{2}(2)P_{1}(1).
$$ 
Recall that $P_{2}(2)$ is the projection onto 
${\mathcal H}_{1}^{0}\oplus ({\mathcal H}_{2}\otimes {\mathcal H}_{1}^{0})$ and 
$P_{1}(1)$ is the projection onto ${\mathbb C}\xi \oplus {\mathcal H}_{1}^{0}$, thus
$P_{1}(1)P_{2}(2)=P_{{\mathcal H}_{1}^{0}}$, and therefore we get
$$
x_{1}y_{2}^{m}x_{1}=
x_{1}P_{{\mathcal H}_{1}^{0}}y_{2}^{m}P_{{\mathcal H}_{1}^{0}}x_{1},
$$
for $m\geq 1$. This gives
$$
\varphi(w_{1}x_{1}y_{2}^{m}x_{1}w_{2})
=\varphi(w_{1}x_{1}P_{{\mathcal H}_{1}^{0}}x_{1}w_{2})
\psi(y_{2}^{m}),
$$
since, for any $\zeta\in {\mathcal H}_{1}^{0}$ of norm one, we have
$$
y_{2}^{m}\zeta=\psi(y_{2}^{m})\zeta \;\;\;{\rm mod}\;{\mathcal K}(2)\ominus {\mathcal K}^{(1)}(2),
$$
where $\psi$ is the state associated with $\zeta$. Note that $\psi$ gives the same moments of the 
variable $y_{2}$ irrespective which $\zeta$ is taken and 
they agree with the corresponding moments in the state $\varphi_{2}$. Finally, since 
$P_{{\mathcal H}_{1}^{0}}=P_{{\mathbb C}\xi\oplus {\mathcal H}_{1}^{0}}-P_{{\mathbb C}\xi}$, we have
$$
\varphi(w_{1}x_{1}P_{{\mathcal H}_{1}^{0}}x_{1}w_{2})=\varphi(w_{1}x_{1}^{2}w_{2})-\varphi(w_{1}x_{1})\varphi(x_{1}w_{2}),
$$
which completes the proof of the second orthogonality condition.
\hfill $\blacksquare$
\begin{Corollary}
Let $\mu_{\iota}$ and $\sigma_{\iota}$ be $\varphi$--distributions
of (bounded operators) $X_{\iota}$ and $U_{\iota}$, respectively, for each $\iota \in I$,
and let $\psi_{m-1}$ be the state associated with any unit vector
$\zeta\in {\mathcal H}^{(m-1)}(\iota)$, where $m\in {\mathbb N}$. Then
\begin{enumerate}
\item
the $\psi_{m-1}$--distribution of $U_{\iota}(m)$ agrees with 
$\sigma_{\iota}$ for every $m\geq 2$ and $\iota \in I$,
\item 
it holds that $\sigma_{\iota}=\mu_{\iota}\angle \sigma_{\overline{\iota}}$
for each $\iota \in I$,
\item
if $X_1,X_2$ are positive and $\mu_1,\mu_2 \neq \delta_{0}$, then 
$\sigma_{\iota}$ agrees with the s-free convolution $\mu_{\iota}\boxslash \mu_{\overline{\iota}}$ 
of probability measures on ${\mathbb R}_{+}$ for each $\iota\in I$,
\item
if $X_1,X_2$ are unitary, then $\sigma_{\iota}$ agrees with the s-free convolution 
$\mu_{\iota}\boxslash\,\mu_{\overline{\iota}}$ of probability measures on ${\mathbb T}$
for each $\iota\in I$.
\end{enumerate}
\end{Corollary}
{\it Proof.}
Assertion (1) follows from the definition of operators $U_{\iota}(m)$.
Then, (2) is a consequence of (1) and Theorem 8.3.
Finally, it follows from Theorems 7.3 and 8.3 that the $\varphi$--distributions of 
$U_{\iota}$ satisfy the subordination equations (2.1).
Uniqueness of the subordination functions for the cases considered in (3) and (4) 
give the assertions.
\hfill $\blacksquare$\\  
\indent{\par}
Corollary 8.1, together with Corollary 4.6, naturally lead to the definition of a 
sequence of iterations of the s-free multiplicative convolution. Thus, for 
given $\mu_1, \mu_2\in \Sigma$ let
\begin{equation}\tag{8.11}
\mu_1 \angle_1 \mu_2 =\mu_1\angle \mu_2, \;\;\;\mu_1\angle_n\mu_2=\mu_1\angle (\mu_2\angle_{n-1}\mu_1),
\end{equation}
for $n\geq 2$ (of course, this sequence can also be defined by for 
any measures, for which the operation $\angle$ has been defined). 
\begin{Corollary}
If $\mu_1,\mu_2\in {\mathcal M}_{{\mathbb R}_{+}}\setminus \{\delta_{0}\}$ 
are compactly supported, then 
\begin{equation}\tag{8.12}
w-\lim_{n\rightarrow \infty}(\mu_1\angle_n \mu_2)= \mu_1\boxslash \mu_2
\end{equation}
\begin{equation}\tag{8.13}
w-\lim_{n\rightarrow \infty}(\mu_1\circlearrowright (\mu_2\angle_n \mu_1))=
\mu_1\boxtimes \mu_2,
\end{equation}
and the corresponding $\eta$--transforms converge uniformly on the compact subsets
of ${\mathbb C}\setminus {\mathbb R}_{+}$ to the $\eta$--transforms
of the limit measures. An analogous result holds for $\mu_1,\mu_2\in {\mathcal M}_{*}$.
\end{Corollary}
{\it Proof.}
In view of Corollary 4.6, for fixed $m\in {\mathbb N}$, 
$(\mu_1\angle_n\mu_2)(m)=(\mu_1\angle_m\mu_2)(m)$ for all $n\geq m$, from which
we obtain convergence of moments as $\rightarrow \infty$ for any distributions
$\mu_1, \mu_2\in\Sigma$. 
Since the measures are compactly supported, 
this implies weak convergence of the corresponding measures.
Therefore, the $\eta$--transforms converge
uniformly on compact subsets of ${\mathbb C}\setminus {\mathbb R}_{+}$.
Moreover, by Corollary 8.1, the limit distributions agree
with the distributions of the positive subordination branch $Y_{1}$, which 
is $\mu_1\boxslash \mu_2$. The arguments for the weak convergence
of the second sequence are similar. An analogous proof holds for
$\mu_1, \mu_2\in {\mathcal M}_{*}$.
\hfill $\blacksquare$
\begin{Remark}
{\rm Informally, the weak limits of Corollary 8.2 can be written
in the following form:
\begin{eqnarray*}
\mu_1 \,\boxslash \,\mu_2&=&\mu_1 \angle (\mu_2\angle (\mu_1\angle (\mu_2\angle (\ldots)))),\\
\mu_1 \,\boxtimes \,\mu_2&=&\mu_1 \circlearrowright (\mu_2\angle (\mu_1\angle (\mu_2\angle (\ldots)))),
\end{eqnarray*}
whereas their transforms in the `continued composition form':
\begin{eqnarray*}
\rho_{\mu_1 \,\boxslash\,\mu_2}(z)&=&\rho_{\mu_1}(z\rho_{\mu_2}(z\rho_{\mu_1}(z\rho_{\mu_2}(\ldots )))),\\
\eta_{\mu_1 \,\boxtimes \,\mu_2}(z)&=&\eta_{\mu_1}(z\rho_{\mu_2}(z\rho_{\mu_1}(z\rho_{\mu_2}(\ldots )))),
\end{eqnarray*}
where the right-hand sides are understood as the uniform limits on compact subsets of 
${\mathbb C}\setminus {\mathbb R}_{+}$ or ${\mathbb D}$.
Actually, these formulas can be used to compute some examples,
including simple examples of $\mu_1\boxtimes \mu_2$ (without using S-transforms).}
\end{Remark}
\begin{Remark}
{\rm By Corollary 8.2, one can {\it define} the s-free and 
free convolutions of any compactly supported measures
$\mu_1, \mu_2$ on ${\mathbb C}$ as the weak limits of the form (8.12)-(8.13). 
This allows us to use operatorial subordination results of Sections 7-8, including 
those which concern distributions of products of free bounded random variables 
as well as distributions of bounded branches and express them in terms of s-free convolutions, 
denoted with the same symbol $\boxslash$ and understood as weak limits 
of type (8.12)-(8.13). 
Using the weak limits, we also get
\begin{equation}\tag{8.14}
\delta_{0}\boxslash \mu=\delta_{0}\;\;\; {\rm and}\;\;\;
\mu \boxslash \delta_{0}=\delta_{\mu(X)}
\end{equation}
and
\begin{equation}\tag{8.15}
\mu\boxtimes \delta_{0}=\delta_{0}=\delta_{0}\boxtimes \mu,
\end{equation}
for any compactly supported probability measure on ${\mathbb C}$,
where we used the well-known relations
$\mu \circlearrowright \delta_{0}=\delta_{0}= 
\delta_{0}\circlearrowright \delta_{\mu(X)}$.}
\end{Remark}
\begin{Example}
{\rm Using Corollary 8.2 and the results of Example 5.1, we obtain 
$\mu\boxslash \delta_a=\mu\angle \delta_a=
S_{a}(\mu)$ and $\delta_{a}\boxslash \mu=\delta_{a}\angle\mu=\delta_{a}$ 
for compactly supported $\mu\in {\mathcal M}_{{\mathbb R}_{+}}$ and $a>0$. This gives 
$$
\delta_{a}\boxtimes \mu=\delta_{a}\circlearrowright S_{a}\mu= D_{a}\mu,
$$
in view of (1.5), since the corresponding $\eta$--transform is of the form
$$
\eta_{\delta_{a}}(\eta_{S_{a}\mu}(z))
=\eta_{D_{a}\mu}(z),
$$
where we used $\eta_{\delta_{a}}(z)=az$ and $\eta_{S_{a}\mu}(z)=\eta_{D_{a}\mu}(z)/a$.
On the other hand,
$$
\mu\boxtimes \delta_{a}=\mu\circlearrowright \delta_{a}=D_{a}\mu,
$$
since $\eta_{\mu\circlearrowright \delta_{a}}(z)=\eta_{\mu}(az)=\eta_{D_{a}\mu}(z)$.
In particular, we have $\delta_{a}\boxslash \delta_{b}=\delta_{a}$ for $a,b>0$.}
\end{Example}
\begin{Example}
{\rm Using the results of Example 5.3, we obtain
$$
\mu \boxslash \delta_{a}=S_{a}\mu\;\;\; {\rm and}\;\;\; \delta_{a}\boxslash \mu = \delta_{a},
$$
for compactly supported $\mu\in {\mathcal M}_{{\mathbb T}}$ and $a\in{\mathbb T}$,
which gives
$$
\mu\boxtimes \delta_{a}=D_{a}\mu= \delta_{a}\boxtimes \mu,
$$
by the same arguments as those used in Example 8.1.}
\end{Example}

Let us finally show that the multiplicative subordination branches are 
related to the notion of s-free independence, as it was the case
of their additive counterparts. For that purpose, let us recall
this concept [16].
\begin{Definition}
{\rm Let $({\mathcal A},\varphi, \psi)$ be a unital algebra with a pair
of linear normalized functionals.
Let ${\mathcal A}_{1}$ be a unital subalgebra of ${\mathcal A}$ and let
${\mathcal A}_{2}$ be a non-unital subalgebra with an `internal' unit $1_{2}$, i.e.
$1_{2}b=b=b1_{2}$ for every $b\in {\mathcal A}_{2}$.
We say that the pair $({\mathcal A}_{1},{\mathcal A}_{2})$ is {\it free with subordination}, or 
simply {\it s-free},
with respect to $(\varphi, \psi)$ if $\psi(1_{2})=1$ and it holds that\\
\indent{\par}
(i)
$\varphi(a_{1}a_{2}\ldots a_{n})=0$ whenever $a_{j}\in {\mathcal A}_{i_{j}}^{0}$ and
$i_{1}\ne i_{2}\neq \ldots \neq i_{n}$
\indent{\par}
(ii)
$\varphi(w_{1}1_{2}w_{2})=\varphi(w_{1}w_{2})-\varphi(w_{1})\varphi(w_{2})$
for any $w_{1},w_{2}\in {\rm alg}({\mathcal A}_{1}, {\mathcal A}_{2})$,\\[5pt]
where ${\mathcal A}_{1}^{0}={\mathcal A}_{1}\cap {\rm ker}\varphi$ and ${\mathcal A}_{2}^{0}={\mathcal A}_{2}\cap {\rm ker}\psi$.
We say that the pair $(a,b)$ of random variables from ${\mathcal A}$
is {\it s-free} with respect to $(\varphi, \psi)$ if there exists $1_2\in {\mathcal A}$
such that $({\mathcal A}_{1},{\mathcal A}_{2})$, where ${\mathcal A}_{1}:={\rm alg}(1,a)$ and  
${\mathcal A}_{2}:={\rm alg}(1_{2},b)$, is s-free w.r.t. $(\varphi, \psi)$.}
\end{Definition}
\begin{Proposition}
The pair $(R_{\iota}^{\scriptscriptstyle {\rm odd}}-1, R_{\overline{\iota}}^{\scriptscriptstyle {\rm even}}-1)$ 
is s-free with respect to the pair $(\varphi, \psi)$, where 
$\psi$ is the state associated with any unit vector $\zeta \in {\mathcal H}_{\iota}^{0}$,
and the projection associated with the second variable is $1_{2}=1-P_{\xi}$.
\end{Proposition}
{\it Proof.}
For simplicity, choose $\iota=1$ and denote 
$r_{1}=R_{1}^{\scriptscriptstyle {\rm odd}}-1$, $r_{2}=R_{2}^{\scriptscriptstyle {\rm even}}-1$.
We have
$$
r_1=\sum_{j \;{\rm odd}}x_1(j)\;\;\;{\rm and}\;\;\;
r_2=\sum_{j \;{\rm even}}x_2(j).
$$
Clearly, the definition of $1_2$ immediately gives $\psi(1_2)=1$ as well as 
condition (ii) of s-freeness. Let us show that condition (i) also holds.
Let ${\mathcal A}_{1}={\rm alg}(r_1,1)$ and ${\mathcal A}_{2}={\rm alg}(r_2,1_2)$.
If in (i) we take $a_{n}\in {\mathcal A}_{2}\cap {\rm ker}\psi$, then $\lambda(a_{n})\xi=0$,
where $\lambda$ is the free product representation, and hence (i)
holds. If $a_{n},a_{n-2},\ldots\in {\mathcal A}_{1}\cap {\rm ker}\varphi$ and
$a_{n-1}, a_{n-3}, \ldots \in {\mathcal A}_{2}\cap {\rm ker}\psi$, then
the GNS representation space for such a mixed moment is ${\mathcal K}_{1}$
(none of $a_{n-1},a_{n-3}, \ldots $ ever gets to act on $\xi$). Therefore,
in the computation of this kind of moment, one can replace $r_1,r_{2}$ and $1_2$
by $X_1, X_2$ and $1$, respectively. But then (i) follows from the usual
freeness condition.  \hfill $\blacksquare$

\section{Free product and s-free loop product of graphs}

In this section we prove the Multiplication Theorem for 
s-free independence and free independence.

Let us recall the definition of the $s$-free product of rooted graphs [2,16].
\begin{Definition}
{\rm By the $s$-{\it free} product of rooted graphs $({\mathcal G}_{1},e_{1})$
and $({\mathcal G}_{2}, e_{2})$, denoted  $({\mathcal G}_{1}\oright {\mathcal G}_{1}, e_{1})$,
we understand the inductive limit of the sequence $({\mathcal B}_{1}(m),e_{1})_{m\in {\mathbb N}}$ 
of rooted graphs, where ${\mathcal B}_{1}(1)={\mathcal G}_{1}\vdash {\mathcal G}_{2}$
and ${\mathcal B}_{\iota}(m)$ is obtained from ${\mathcal B}_{\iota}(m-1)$ by attaching
by its root a copy of ${\mathcal G}_{1}$ (if $m$ is even), or a copy of ${\mathcal G}_{2}$ (if $m$ is odd) 
to every vertex of the difference ${\mathcal B}_{\iota}(m-1)\setminus {\mathcal B}_{\iota}(m-2)$.}
\end{Definition}
\begin{Remark}
{\rm It can be seen that $({\mathcal B}_{1}(m),e_{1})_{m\in {\mathbb N}}$ is
a sequence of growing graphs with the root kept to be $e_{\iota}$, and therefore, 
the inductive limit exists. In a similar way we define  $({\mathcal G}_{2}\oright {\mathcal G}_{1}, e_{2})$.
It is easy to see that 
\begin{equation}\tag{9.1}
{\mathcal B}_{1}(m)={\mathcal G}_{1}\vdash {\mathcal B}_{2}(m-1)\;\;\;{\rm and}\;\;\;
{\mathcal B}_{2}(m)={\mathcal G}_{2}\vdash {\mathcal B}_{1}(m-1),
\end{equation}
which can also serve an inductive definition of the considered sequences. 
These formulas are important since they correspond to the decomposition
of the s-free additive convolution of the form
\begin{equation}\tag{9.2}
\mu_1 \boxright \mu_2=\mu_1\vdash(\mu_2\vdash(\mu_1\vdash(\mu_2\vdash\ldots ))),
\end{equation}
where $\mu_{\iota}$ is the spectral distribution of ${\mathcal G}_{\iota}$ [2].}
\end{Remark}

Using the decomposition (9.2), we proved the Addition Theorem 
for s-free independence in [2].
Our first goal in this Section is to prove the corresponding Multiplication Theorem.
For that purpose, we introduce a new notion of graph product, called 
the `s-free loop product'. We follow it up with a proposition which
justifies the definition.
\begin{Definition}
{\rm Suppose that the s-free product of graphs $({\mathcal G}_{1}\oright {\mathcal G}_{2}, 
e_{1})$ is naturally colored. 
The {\it s-free loop product} of $({\mathcal G}_{1}, e_{1})$
and $({\mathcal G}_{2}, e_{2})$ is 
the graph $({\mathcal G}_{1}\oright_{\ell} {\mathcal G}_{2}, e_{1})$
obtained from $({\mathcal G}_{1}\oright {\mathcal G}_{2}, e_{1})$ 
by attaching a loop of color $2$ to the root $e_{1}$.}
\end{Definition}
\begin{figure}
\unitlength=1mm
\special{em.linewidth 2pt}
\linethickness{0.5pt}
\begin{picture}(120.00,35.00)(0.00,5.00)
\put(10.00,11.50){\circle{3.00}}
\put(3.00,13.00){\circle*{1.00}}
\put(10.00,10.00){\circle*{1.00}}
\put(17.00,13.00){\circle*{1.00}}
\put(-2.00,13.00){\scriptsize ${\mathcal G}_{1}$}
\put(9.00,7.00){\scriptsize $e_{1}$}
\put(2.00,15.00){\scriptsize $x$}
\put(16.00,15.00){\scriptsize $x'$}

\qbezier(3,13)(3,10)(10,10)
\qbezier(10,10)(17,10)(17,13)
\put(10.00,20.00){\circle*{1.00}}
\put(3.00,23.00){\circle*{1.00}}
\put(17.00,23.00){\circle*{1.00}}
\put(-2.00,23.00){\scriptsize ${\mathcal G}_{2}$}
\put(9.00,17.00){\scriptsize $e_{2}$}
\put(2.00,25.00){\scriptsize $y$}
\put(16.00,25.00){\scriptsize $y'$}

\qbezier(3,23)(3,20)(10,20)
\qbezier(10,20)(17,20)(17,23)
\put(30.00,15.00){\scriptsize ${\mathcal G}_{1}\oright_{\ell} {\mathcal G}_{2}$}
\put(80.00,11.50){\circle{3.00}}
\put(80.00,13.00){\circle{6.00}}
\put(80.00,10.00){\circle*{1.00}}
\put(60.00,18.00){\circle*{1.00}}
\put(100.00,18.00){\circle*{1.00}}
\put(79.00,7.00){\scriptsize $e$}

\qbezier(60,18)(60,10)(80,10)
\qbezier(80,10)(100,10)(100,18)
\put(50.00,23.00){\circle*{1.00}}
\put(70.00,23.00){\circle*{1.00}}
\qbezier(50,23)(50,18)(60,18)
\qbezier(60,18)(70,18)(70,23)

\put(90.00,23.00){\circle*{1.00}}
\put(110.00,23.00){\circle*{1.00}}
\qbezier(90,23)(90,18)(100,18)
\qbezier(100,18)(110,18)(110,23)
\put(45.00,26.00){\circle*{1.00}}
\put(55.00,26.00){\circle*{1.00}}
\qbezier(45,26)(45,23)(50,23)
\qbezier(50,23)(55,23)(55,26)

\put(65.00,26.00){\circle*{1.00}}
\put(75.00,26.00){\circle*{1.00}}
\qbezier(65,26)(65,23)(70,23)
\qbezier(70,23)(75,23)(75,26)

\put(85.00,26.00){\circle*{1.00}}
\put(95.00,26.00){\circle*{1.00}}
\qbezier(85,26)(85,23)(90,23)
\qbezier(90,23)(95,23)(95,26)

\put(105.00,26.00){\circle*{1.00}}
\put(115.00,26.00){\circle*{1.00}}
\qbezier(105,26)(105,23)(110,23)
\qbezier(110,23)(115,23)(115,26)

\put(50.00,24.50){\circle{3.00}}
\put(70.00,24.50){\circle{3.00}}
\put(90.00,24.50){\circle{3.00}}
\put(110.00,24.50){\circle{3.00}}

\put(42.00,28.00){\circle*{1.00}}
\put(48.00,28.00){\circle*{1.00}}
\qbezier(42,28)(42,26)(45,26)
\qbezier(45,26)(48,26)(48,28)

\put(52.00,28.00){\circle*{1.00}}
\put(58.00,28.00){\circle*{1.00}}
\qbezier(52,28)(52,26)(55,26)
\qbezier(55,26)(58,26)(58,28)

\put(62.00,28.00){\circle*{1.00}}
\put(68.00,28.00){\circle*{1.00}}
\qbezier(62,28)(62,26)(65,26)
\qbezier(65,26)(68,26)(68,28)

\put(72.00,28.00){\circle*{1.00}}
\put(78.00,28.00){\circle*{1.00}}
\qbezier(72,28)(72,26)(75,26)
\qbezier(75,26)(78,26)(78,28)

\put(82.00,28.00){\circle*{1.00}}
\put(88.00,28.00){\circle*{1.00}}
\qbezier(82,28)(82,26)(85,26)
\qbezier(85,26)(88,26)(88,28)

\put(92.00,28.00){\circle*{1.00}}
\put(98.00,28.00){\circle*{1.00}}
\qbezier(92,28)(92,26)(95,26)
\qbezier(95,26)(98,26)(98,28)

\put(102.00,28.00){\circle*{1.00}}
\put(108.00,28.00){\circle*{1.00}}
\qbezier(102,28)(102,26)(105,26)
\qbezier(105,26)(108,26)(108,28)

\put(112.00,28.00){\circle*{1.00}}
\put(118.00,28.00){\circle*{1.00}}
\qbezier(112,28)(112,26)(115,26)
\qbezier(115,26)(118,26)(118,28)

\put(42.00,29.00){\circle{2.00}}
\put(48.00,29.00){\circle{2.00}}
\put(52.00,29.00){\circle{2.00}}
\put(58.00,29.00){\circle{2.00}}

\put(62.00,29.00){\circle{2.00}}
\put(68.00,29.00){\circle{2.00}}
\put(72.00,29.00){\circle{2.00}}
\put(78.00,29.00){\circle{2.00}}
\put(82.00,29.00){\circle{2.00}}
\put(88.00,29.00){\circle{2.00}}
\put(92.00,29.00){\circle{2.00}}
\put(98.00,29.00){\circle{2.00}}

\put(102.00,29.00){\circle{2.00}}
\put(108.00,29.00){\circle{2.00}}
\put(112.00,29.00){\circle{2.00}}
\put(118.00,29.00){\circle{2.00}}

\qbezier(40.25,30)(40.25,28)(42,28)
\qbezier(42,28)(43.75,28)(43.75,30)
\qbezier(46.25,30)(46.25,28)(48,28)
\qbezier(48,28)(49.75,28)(49.75,30)

\qbezier(50.25,30)(50.25,28)(52,28)
\qbezier(52,28)(53.75,28)(53.75,30)
\qbezier(56.25,30)(56.25,28)(58,28)
\qbezier(58,28)(59.75,28)(59.75,30)

\qbezier(60.25,30)(60.25,28)(62,28)
\qbezier(62,28)(63.75,28)(63.75,30)
\qbezier(66.25,30)(66.25,28)(68,28)
\qbezier(68,28)(69.75,28)(69.75,30)

\qbezier(70.25,30)(70.25,28)(72,28)
\qbezier(72,28)(73.75,28)(73.75,30)
\qbezier(76.25,30)(76.25,28)(78,28)
\qbezier(78,28)(79.75,28)(79.75,30)

\qbezier(80.25,30)(80.25,28)(82,28)
\qbezier(82,28)(83.75,28)(83.75,30)
\qbezier(86.25,30)(86.25,28)(88,28)
\qbezier(88,28)(89.75,28)(89.75,30)

\qbezier(90.25,30)(90.25,28)(92,28)
\qbezier(92,28)(93.75,28)(93.75,30)
\qbezier(96.25,30)(96.25,28)(98,28)
\qbezier(98,28)(99.75,28)(99.75,30)

\qbezier(100.25,30)(100.25,28)(102,28)
\qbezier(102,28)(103.75,28)(103.75,30)
\qbezier(106.25,30)(106.25,28)(108,28)
\qbezier(108,28)(109.75,28)(109.75,30)

\qbezier(110.25,30)(110.25,28)(112,28)
\qbezier(112,28)(113.75,28)(113.75,30)
\qbezier(116.25,30)(116.25,28)(118,28)
\qbezier(118,28)(119.75,28)(119.75,30)

\end{picture}
\caption{An example of the s-free loop product}
\end{figure}
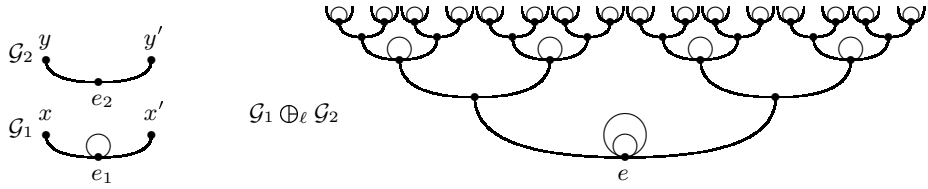
\begin{Example}
{\rm Consider the s-free loop product of graphs given in Fig.4, which 
can be viewed as a `binary tree with loops' 
(the product graph is of course infinite, only our picture 
is truncated at distance 5 from the root). 
Note that graph ${\mathcal G}_{1}$ has a loop at the root, 
and so do copies of ${\mathcal G}_{1}$ in the product graph,
whereas ${\mathcal G}_{2}$ does not have any loops. 
We produce ${\mathcal G}_{1}\oright {\mathcal G}_{2}$ in the 
inductive way (Definition 9.1 is used).
Now, in order to obtain the s-free loop product, we attach
one additional loop of color $2$ to the root (we draw this loop 
larger than loops of color $1$). On the level
of the adjacency matrices, this additional loop
corresponds to the projection $P_{\xi}$ as the proof of Proposition 9.1
shows. Note that the main purpose of attaching this extra loop 
is to fit the s-free multiplicative convolution into the general scheme
of the Multiplication Theorem, where counting rooted alternating d-walks 
leads to the moments of multiplicative convolutions. 
Direct computations give, for instance, 
$D_{2}(e_1)=1$, $D_{4}(e_1)=0$, $D_{6}(e_{1})=4$, $D_{8}(e_1)=0$.}
\end{Example}
\begin{Proposition}
Let $A_{\iota}$ be the adjacency matrix of ${\mathcal G}_{\iota}$, where
$\iota \in I$. The adjacency matrix of the s-free loop product of ${\mathcal G}_{1}$
and ${\mathcal G}_{2}$ takes the form
\begin{equation}\tag{9.3}
A({\mathcal G}_{1}\oright_{\ell}{\mathcal G}_{2})=R_{1}+R_{2}
\end{equation}
where
\begin{equation}\tag{9.4}
R_{1}=\sum_{\scriptscriptstyle j\; {\rm odd}}A_{1}(j)\;\;\;and
\;\;\;
R_{2}=P_{\xi}+\sum_{\scriptscriptstyle k\; {\rm even}}A_{2}(k).
\end{equation}
Moreover, the $\varphi_{e}$-distribution of $R_{1}$ is $\mu_{1}$, the $\psi$-distribution 
of $R_2$ is $\mu_2$ and the pair $(R_{1}-1,R_{2}-1)$ is s-free independent w.r.t. $(\varphi_{e}, \psi)$, 
where $\psi$ is the state associated with any unit vector $\zeta\in {\mathcal H}_{1}^{0}=l_{2}(V_{1}^{0})$. 
\end{Proposition}
{\it Proof.}
It follows from [2,16] that the adjacency matrix of the s-free product of graphs 
is of the form
$$
A({\mathcal G}_{1}\oright{\mathcal G}_{2})=S_{1}+S_{2},
$$
where 
$$
S_{1}=\sum_{\scriptscriptstyle j\; {\rm odd}}A_{1}(j)\;\;\;{\rm and}
\;\;\;
S_{2}=\sum_{\scriptscriptstyle k\; {\rm even}}A_{2}(k).
$$
By Definition 9.1, the adjacency matrix of the corresponding 
s-free loop product is obtained by adding the projection $P_{0}$ 
to the adjacency matrix of the subgraph
of color $2$, namely $S_{2}$, which corresponds to 
attaching the loop of color $2$ to the root. 
This gives
$$
A({\mathcal G}_{1}\oright_{\ell}{\mathcal G}_{2})=R_{1}+R_{2},
$$
where 
$$
R_{1}=S_{1}\;\;\;{\rm and}
\;\;\;
R_{2}=P_{\xi}+ S_{2}.
$$
In view of Proposition 8.2, the pair $(R_{1}-1,R_{2}-1)$ is s-free independent 
w.r.t. the pair $(\varphi, \psi)$, where $\psi$ is the state associated with 
any unit vector $\zeta\in {\mathcal H}_{1}^{0}$ (here, $1$ denotes the unit on 
$l_{2}(V)$, where $V$ is the set of vertices of ${\mathcal G}_{1}\oright_{\ell}{\mathcal G}_{2}$,
which corresponds to ${\mathcal K}_{1}$).
Moreover, it is easy to check that 
the $\varphi_{e}$-distribution of $R_{1}$ is $\mu_{1}$ and the $\psi$-distribution
of $R_2$ is $\mu_2$. This completes the proof.\hfill $\blacksquare$\\
\indent{\par}
Thus, the first part of the Multiplication Theorem for s-free independence
is proved. We now need to find a connection between the `first return moments'
of $R_{2}R_{1}$ and cardinalities of the sets $D_{2n}$. In that context, 
note that in order that the set $D$ of all rooted alternating d-walks on 
${\mathcal G}_{1}\oright_{\ell} {\mathcal G}_{2}$ be non-empty, either ${\mathcal G}_{1}$ or 
${\mathcal G}_{2}$ must have a loop at $e_{1}$ or $e_{2}$, respectively. 
In fact, otherwise, starting an alternating 
walk from the root of ${\mathcal G}_{1}\oright {\mathcal G}_{2}$, one cannot return 
to it since the $m$-th edge of that walk must be of different color than the $m-1$-th edge
and thus it must belong to ${\mathcal B}_{j}(m)\setminus {\mathcal B}_{j}(m-1)$, and therefore 
its distance from the root (understood as the distance of the vertex closer to the root) 
equals $m-1$ and thus tends to $\infty$ as $m\rightarrow \infty$. 
\begin{Theorem}
The Multiplication Theorem holds for s-free independence, the associated loop 
product ${\mathcal G}_{1}\oright_{\ell}{\mathcal G}_{2}$, and the multiplicative 
convolution $\mu_1\boxslash \mu_2$.
\end{Theorem}
{\it Proof.}
For notational simplicity, denote $\sigma_1=\mu_1\boxslash \mu_2$,
$\sigma_2=\mu_2\boxslash \mu_1$ and 
$$
{\mathcal B}_{\iota}={\mathcal G}_{\iota}\oright {\mathcal G}_{\overline{\iota}}, \;\;\;
{\mathcal B}_{\iota}^{\ell}={\mathcal G}_{\iota}\oright_{\ell} {\mathcal G}_{\overline{\iota}}
$$
From Proposition 9.1 and the results of Section 8, it follows that
\begin{equation}\tag{9.5}
N_{R_2R_1}(n)=N_{\sigma_1}(n)
\end{equation}
for $n\in {\mathbb N}$. Now, we need to prove that
\begin{equation}\tag{9.6}
N_{\sigma_1}(n)=|D_{2n}(e)|
\end{equation}
for $k\in {\mathbb N}$. Note first that any f-walk 
$w\in F(e_{1})$ on ${\mathcal B}_{1}$ 
must begin and terminate with an edge of color 1 since no edge of 
color $2$ is incident on $e_{1}$. Since ${\mathcal B}_{1}^{\ell}$
is obtained from ${\mathcal B}_{1}$ by attaching a loop of color $2$
to the root $e_1$, the set of rooted alternating d-walks on 
${\mathcal B}_{1}^{\ell}$ 
of lenght $2n$ is in 1-1 correspondence with the set of rooted alternating 
f-walks on ${\mathcal B}_{1}$ of lenght
$2n-1$. In fact, each rooted alternating d-walk of lenght $2n$ 
on ${\mathcal B}_{1}^{\ell}$ is of the form $w=(u_1,u_2)$, 
where $u_1$ is a rooted alternating f-walk of lenght $2n-1$ 
and $u_2$ is the rooted loop of color 2. Therefore, 
\begin{equation}\tag{9.7}
|D_{2n}(e_1)|=|F_{2n-1}(e_1)|,
\end{equation}
where $D_{2n}(e_1)$ refers to the loop product and $F_{2n-1}(e_1)$ -- to the usual
product. 
This fact will be used in the induction proof given below.
If $n=1$, we use this correspondence to get
$N_{\sigma_1}(1)=N_{\mu_1}(1)=|F_{1}(e_{1})|=|D_{2}(e_1)|$.
Suppose now that (9.6) holds if $n$ is replaced by $1\leq k \leq n-1$, where
$\sigma_1$ is understood to be the s-free convolution of any
$\sigma'$ and $\sigma''$ (thus, in particular, $\mu_2$ and $\mu_1$).
Using the proof of Theorem 6.1, we have
$$
N_{\sigma_1}(n)=
\sum_{r=1}^{n}N_{\mu_1}(r)\sum _{k_{1}+k_{2}+\ldots +k_{r-1}=n-1}
N_{\sigma_2}(k_{1})N_{\sigma_2}(k_{2}) \ldots N_{\sigma_{2}}(k_{r-1}).
$$
In view of (9.7), it is enough to show that the expression on the RHS is 
equal to $|F_{2n-1}(e_1)|$. By the inductive assumption, we replace in the above 
formula each $N_{\sigma_2}(k_{i})$ by $|D_{2k_{i}}(e_{1})|$, and then, using (9.7), 
by $|F_{2k_{i}-1}(e_{1})|$, $1\leq i \leq r-1$. 
Therefore, we need to justify the above formula, viewing 
$N_{\sigma_{1}}(n)$ as the number of rooted alternating f-walks on ${\mathcal B}_{1}$
and $N_{\sigma_{2}}(k_{i})$, $1\leq i \leq r-1$, as numbers of rooted alternating 
f-walks on ${\mathcal B}_{2}$. Since
$$
{\mathcal B}_{1}={\mathcal G}_{1}\vdash {\mathcal B}_{2}
$$
we can observe that each rooted alternating f-walk of lenght $2n-1$ 
on ${\mathcal B}_{1}$ consists of a rooted f-walk 
$c=(v_{0}, v_{1}, \ldots , v_{r})\in F_{r}(e_{1})$ 
of color $1$, for some $1\leq r \leq n$, with $r-1$ 
alternating f-walks $w_{i}\in F(v_{i})$ on copies of
${\mathcal B}_{2}$ attached to vertices $v_{i}$, where $1\leq i \leq r-1$.
The latter must begin and terminate with edges of color $2$ since 
no edge of color $1$ is incident on the roots of these copies. 
Thus, the expression given on the RHS of the above formula
gives the numbers of all alternating rooted f-walks of total lenght 
$$
r+(2k_{1}-1)+(2k_{2}-1)+\ldots + (2k_{r-1}-1)=r+2(n-1)-(r-1)=2n-1
$$
summed over $1\leq r \leq n$, which completes the proof.\hfill $\blacksquare$
\begin{Example}
{\rm The lowest order first return moments $N_{\mu_1 \boxslash \,\mu_2}(n)$ are given by
\begin{eqnarray*}
N_{\mu_1 \boxslash \,\mu_2}(1) &=& N_{\mu_1}(1),\\
N_{\mu_1 \boxslash \,\mu_2}(2) &=& N_{\mu_1}(2)N_{\mu_2}(1),\\
N_{\mu_1 \boxslash \,\mu_2}(3) &=& N_{\mu_1}(3)N_{\mu_2}^{2}(1)+N_{\mu_1}(2)N_{\mu_1}(1)N_{\mu_2}(2),\\
N_{\mu_1 \boxslash \,\mu_2}(4) &=& N_{\mu_1}(4)N_{\mu_2}^{3}(1)+2N_{\mu_1}(3)N_{\mu_1}(1)N_{\mu_2}(2)N_{\mu_2}(1)\\
&&
+N_{\mu_1}(2)N_{\mu_1}^{2}(1)N_{\mu_2}(3) + N_{\mu_1}^{2}(2)N_{\mu_2}(2)N_{\mu_2}(1)
\end{eqnarray*}
where we used the results of Example 6.2. The enumeration of d-walks 
given in Example 9.1 can be easily verified by substituting to
the above formulas only the non-vanishing moments: $N_{\mu_1}(1)=1$, $N_{\mu_1}(2)=2$,
$N_{\mu_2}(2)=2$. We have $N_{\sigma_1}(1)=1=D_{2}(e_1)$, $N_{\sigma_1}(2)=0=D_{4}(e_1)$,
$N_{\sigma_1}(3)=4=D_{6}(e_1)$ and $N_{\sigma_1}(4)=0=D_{8}(e_1)$, where
$\sigma_1=\mu_1\boxslash \mu_2$. }
\end{Example}
Finally, we will prove the Multiplication Theorem for free independence.
In this case, one does not need to introduce a new type of graph product
since no `unitization' of adjacency matrices is necessary to make them
freely independent (basically, this is because units are identified in the case of the free
product of algebras). Nevertheless, if neither of the graphs, ${\mathcal G}_{1}$ or
${\mathcal G}_{2}$, have loops at their roots, then the set of rooted alternating 
d-walks on ${\mathcal G}_{1}*{\mathcal G}_{2}$ is empty.
\begin{Theorem}
The Multiplication Theorem holds for free independence, the associated product graph
${\mathcal G}_{1}*{\mathcal G}_{2}$, and the multiplicative convolution $\mu_1\boxtimes\mu_2$.
\end{Theorem}
{\it Proof.}
First of all, we know [2] that
$$
A({\mathcal G}_{1}*{\mathcal G}_{2})=S_{1}+S_{2}
$$
where $S_{\iota}=\sum_{k=1}^{\infty}A_{\iota}(k)$, $\iota \in I$.
Moreover, the pairs $(S_{1},S_{2})$ and $(S_{1}-1,S_{2}-1)$
are free w.r.t. $\varphi$. Also, it is clear that the moments of
$S_2S_1$ agree with the moments of $\mu_1\boxtimes \mu_2$.
It remains to be shown that the latter coincide with the 
corresponding numbers of rooted alternating d-walks.
Using Proposition 2.1, we have $\mu_1\boxtimes\,\mu_2=\mu_1 \circlearrowright (\mu_2\boxslash \mu_1)$ 
and therefore, we can employ the combinatorial formula used in the proof
of Theorem 3.1 to get
$$
N_{\sigma_1}(n)=\sum_{r=1}^{n}
N_{\mu_1}(r)
\sum_{k_{1}+k_{2}+\ldots + k_{r}=n}
N_{\sigma_2}(k_{1})
N_{\sigma_2}(k_{2})
\ldots
N_{\sigma_2}(k_{r}).
$$
Now, we can decompose the free product of graphs as
$$
{\mathcal  G}_{1}*{\mathcal  G}_{2}={\mathcal  G}_{1}\vartriangleright {\mathcal  B}_{2}
$$  
with the root $e$ of ${\mathcal G}_{1}*{\mathcal G}_{2}$ obtained by identifying 
the root $e_{1}$ of ${\mathcal G}_{1}$ with 
the root $e_{2}$ of ${\mathcal B}_{2}$ (see [2,16,20]).
In view of Proposition 9.1, 
$$
N_{\sigma_2}(k_{i})=|F_{2k_{i}-1}(e_{2})|,
$$
where $e_{2}$ is the root in ${\mathcal B}_{2}$.
Of course, $N_{\mu_1}(r)=|F_{r}(e_{1})|$, where $e_{1}$ is treated as the root of ${\mathcal G}_{1}$.
Therefore, the RHS of the above combinatorial formula gives the number of rooted walks 
$w$ on ${\mathcal G}_{1}*{\mathcal G}_{2}$ consisting of an f-walk $c=(v_{0},v_{1}, \ldots , v_{r})\in F_{r}(e)$,
interlaced with rooted alternating f-walks $w_{1},w_{2}, \ldots , w_{r}$ on copies of
${\mathcal  B}_{2}$ attached to vertices $v_{1}, v_{2}, \ldots , v_{r}$, respectively.
In other words, we have
$$
w=(\beta_{1}, w_{1}, \beta_{2}, w_{2}, \ldots , \beta_{r},w_{r}),
$$  
and therefore, it consists of two rooted alternating f-walks, namely 
$$
u_{1}=(\beta_{1}, w_{1}, \beta_{2}, w_{2}, \ldots , \beta_{r})\;\;\; {\rm and}\;\;\;u_{2}=w_{r},
$$
where $u_1$ begins and terminates with an edge of color $1$, and 
$u_{2}$ begins and terminates with an edge of color $2$. 
Therefore, $w$ is a rooted alternating
d-walk of lenght $2n$ which is `subordinate' to an f-walk $c$ of lenght $r$.
The summation over $1\leq r \leq n$ gives all rooted d-walks of lenght $2n$, which
completes the proof. \hfill $\blacksquare$\\
\indent{\par}
In view of this result, the free product of graphs is `complete'
in our `category' of product graphs since it is naturally related
to the rooted alternating d-walks and no additional loops are necessary to
fit it into the scheme of the Multiplication Theorem.
\begin{figure}
\indent{\par}
\unitlength=1mm
\special{em.linewidth 2pt}
\linethickness{0.5pt}
\begin{picture}(120.00,55.00)(10.00,-15.00)
\put(15.00,2.00){\circle{4.00}}
\put(8.00,3.00){\circle*{1.00}}
\put(15.00,0.00){\circle*{1.00}}
\put(22.00,3.00){\circle*{1.00}}
\put(0.00,3.00){\scriptsize ${\mathcal G}_{1}$}
\put(14.00,-3.00){\scriptsize $e_{1}$}
\put(7.00,5.00){\scriptsize $x$}
\put(21.00,5.00){\scriptsize $x'$}

\qbezier(8,3)(8,0)(15,0)
\qbezier(15,0)(22,0)(22,3)
\put(15.00,10.00){\circle*{1.00}}
\put(8.00,13.00){\circle*{1.00}}
\put(22.00,13.00){\circle*{1.00}}
\put(0.00,13.00){\scriptsize ${\mathcal G}_{2}$}
\put(14.00,7.00){\scriptsize $e_{2}$}
\put(7.00,15.00){\scriptsize $y$}
\put(21.00,15.00){\scriptsize $y'$}

\qbezier(8,13)(8,10)(15,10)
\qbezier(15,10)(22,10)(22,13)
\put(36.00,10.00){\scriptsize ${\mathcal G}_{1}* {\mathcal G}_{2}$}
\put(80.00,13.00){\circle{6.00}}
\put(80.00,10.00){\circle*{1.00}}
\put(60.00,18.00){\circle*{1.00}}
\put(100.00,18.00){\circle*{1.00}}
\put(79.00,7.00){\scriptsize $e$}

\qbezier(60,18)(60,10)(80,10)
\qbezier(80,10)(100,10)(100,18)
\put(50.00,23.00){\circle*{1.00}}
\put(70.00,23.00){\circle*{1.00}}
\qbezier(50,23)(50,18)(60,18)
\qbezier(60,18)(70,18)(70,23)

\put(90.00,23.00){\circle*{1.00}}
\put(110.00,23.00){\circle*{1.00}}
\qbezier(90,23)(90,18)(100,18)
\qbezier(100,18)(110,18)(110,23)
\linethickness{0.4pt}
\put(45.00,26.00){\circle*{1.00}}
\put(55.00,26.00){\circle*{1.00}}
\qbezier(45,26)(45,23)(50,23)
\qbezier(50,23)(55,23)(55,26)

\put(65.00,26.00){\circle*{1.00}}
\put(75.00,26.00){\circle*{1.00}}
\qbezier(65,26)(65,23)(70,23)
\qbezier(70,23)(75,23)(75,26)

\put(85.00,26.00){\circle*{1.00}}
\put(95.00,26.00){\circle*{1.00}}
\qbezier(85,26)(85,23)(90,23)
\qbezier(90,23)(95,23)(95,26)

\put(105.00,26.00){\circle*{1.00}}
\put(115.00,26.00){\circle*{1.00}}
\qbezier(105,26)(105,23)(110,23)
\qbezier(110,23)(115,23)(115,26)

\put(50.00,24.50){\circle{3.00}}
\put(70.00,24.50){\circle{3.00}}
\put(90.00,24.50){\circle{3.00}}
\put(110.00,24.50){\circle{3.00}}
\linethickness{0.3pt}
\put(42.00,28.00){\circle*{1.00}}
\put(48.00,28.00){\circle*{1.00}}
\qbezier(42,28)(42,26)(45,26)
\qbezier(45,26)(48,26)(48,28)

\put(52.00,28.00){\circle*{1.00}}
\put(58.00,28.00){\circle*{1.00}}
\qbezier(52,28)(52,26)(55,26)
\qbezier(55,26)(58,26)(58,28)

\put(62.00,28.00){\circle*{1.00}}
\put(68.00,28.00){\circle*{1.00}}
\qbezier(62,28)(62,26)(65,26)
\qbezier(65,26)(68,26)(68,28)

\put(72.00,28.00){\circle*{1.00}}
\put(78.00,28.00){\circle*{1.00}}
\qbezier(72,28)(72,26)(75,26)
\qbezier(75,26)(78,26)(78,28)

\put(82.00,28.00){\circle*{1.00}}
\put(88.00,28.00){\circle*{1.00}}
\qbezier(82,28)(82,26)(85,26)
\qbezier(85,26)(88,26)(88,28)

\put(92.00,28.00){\circle*{1.00}}
\put(98.00,28.00){\circle*{1.00}}
\qbezier(92,28)(92,26)(95,26)
\qbezier(95,26)(98,26)(98,28)

\put(102.00,28.00){\circle*{1.00}}
\put(108.00,28.00){\circle*{1.00}}
\qbezier(102,28)(102,26)(105,26)
\qbezier(105,26)(108,26)(108,28)

\put(112.00,28.00){\circle*{1.00}}
\put(118.00,28.00){\circle*{1.00}}
\qbezier(112,28)(112,26)(115,26)
\qbezier(115,26)(118,26)(118,28)
\linethickness{0.2pt}
\put(42.00,28.75){\circle{1.50}}
\put(48.00,28.75){\circle{1.50}}
\put(52.00,28.75){\circle{1.50}}
\put(58.00,28.75){\circle{1.50}}

\put(62.00,28.75){\circle{1.50}}
\put(68.00,28.75){\circle{1.50}}
\put(72.00,28.75){\circle{1.50}}
\put(78.00,28.75){\circle{1.50}}
\put(82.00,28.75){\circle{1.50}}
\put(88.00,28.75){\circle{1.50}}
\put(92.00,28.75){\circle{1.50}}
\put(98.00,28.75){\circle{1.50}}

\put(102.00,28.75){\circle{1.50}}
\put(108.00,28.75){\circle{1.50}}
\put(112.00,28.75){\circle{1.50}}
\put(118.00,28.75){\circle{1.50}}

\qbezier(40.25,30)(40.25,28)(42,28)
\qbezier(42,28)(43.75,28)(43.75,30)
\qbezier(46.25,30)(46.25,28)(48,28)
\qbezier(48,28)(49.75,28)(49.75,30)

\qbezier(50.25,30)(50.25,28)(52,28)
\qbezier(52,28)(53.75,28)(53.75,30)
\qbezier(56.25,30)(56.25,28)(58,28)
\qbezier(58,28)(59.75,28)(59.75,30)

\qbezier(60.25,30)(60.25,28)(62,28)
\qbezier(62,28)(63.75,28)(63.75,30)
\qbezier(66.25,30)(66.25,28)(68,28)
\qbezier(68,28)(69.75,28)(69.75,30)

\qbezier(70.25,30)(70.25,28)(72,28)
\qbezier(72,28)(73.75,28)(73.75,30)
\qbezier(76.25,30)(76.25,28)(78,28)
\qbezier(78,28)(79.75,28)(79.75,30)

\qbezier(80.25,30)(80.25,28)(82,28)
\qbezier(82,28)(83.75,28)(83.75,30)
\qbezier(86.25,30)(86.25,28)(88,28)
\qbezier(88,28)(89.75,28)(89.75,30)

\qbezier(90.25,30)(90.25,28)(92,28)
\qbezier(92,28)(93.75,28)(93.75,30)
\qbezier(96.25,30)(96.25,28)(98,28)
\qbezier(98,28)(99.75,28)(99.75,30)

\qbezier(100.25,30)(100.25,28)(102,28)
\qbezier(102,28)(103.75,28)(103.75,30)
\qbezier(106.25,30)(106.25,28)(108,28)
\qbezier(108,28)(109.75,28)(109.75,30)

\qbezier(110.25,30)(110.25,28)(112,28)
\qbezier(112,28)(113.75,28)(113.75,30)
\qbezier(116.25,30)(116.25,28)(118,28)
\qbezier(118,28)(119.75,28)(119.75,30)

\linethickness{0.5pt}
\put(80.00,10.00){\circle*{1.00}}
\put(60.00,2.00){\circle*{1.00}}
\put(100.00,2.00){\circle*{1.00}}

\qbezier(60,2)(60,10)(80,10)
\qbezier(80,10)(100,10)(100,2)
\put(50.00,-3.00){\circle*{1.00}}
\put(70.00,-3.00){\circle*{1.00}}
\qbezier(50,-3)(50,2)(60,2)
\qbezier(60,2)(70,2)(70,-3)

\put(60.00,-0.50){\circle{5.00}}
\put(100.00,-0.50){\circle{5.00}}

\put(90.00,-3.00){\circle*{1.00}}
\put(110.00,-3.00){\circle*{1.00}}
\qbezier(90,-3)(90,2)(100,2)
\qbezier(100,2)(110,2)(110,-3)
\linethickness{0.4pt}
\put(45.00,-6.00){\circle*{1.00}}
\put(55.00,-6.00){\circle*{1.00}}
\qbezier(45,-6)(45,-3)(50,-3)
\qbezier(50,-3)(55,-3)(55,-6)

\put(65.00,-6.00){\circle*{1.00}}
\put(75.00,-6.00){\circle*{1.00}}
\qbezier(65,-6)(65,-3)(70,-3)
\qbezier(70,-3)(75,-3)(75,-6)

\put(85.00,-6.00){\circle*{1.00}}
\put(95.00,-6.00){\circle*{1.00}}
\qbezier(85,-6)(85,-3)(90,-3)
\qbezier(90,-3)(95,-3)(95,-6)

\put(105.00,-6.00){\circle*{1.00}}
\put(115.00,-6.00){\circle*{1.00}}
\qbezier(105,-6)(105,-3)(110,-3)
\qbezier(110,-3)(115,-3)(115,-6)

\linethickness{0.3pt}
\put(45.00,-7.25){\circle{2.50}}
\put(55.00,-7.25){\circle{2.50}}
\put(65.00,-7.25){\circle{2.50}}
\put(75.00,-7.25){\circle{2.50}}
\put(85.00,-7.25){\circle{2.50}}
\put(95.00,-7.25){\circle{2.50}}
\put(105.00,-7.25){\circle{2.50}}
\put(115.00,-7.25){\circle{2.50}}

\put(42.00,-8.00){\circle*{1.00}}
\put(48.00,-8.00){\circle*{1.00}}
\qbezier(42,-8)(42,-6)(45,-6)
\qbezier(45,-6)(48,-6)(48,-8)

\put(52.00,-8.00){\circle*{1.00}}
\put(58.00,-8.00){\circle*{1.00}}
\qbezier(52,-8)(52,-6)(55,-6)
\qbezier(55,-6)(58,-6)(58,-8)

\put(62.00,-8.00){\circle*{1.00}}
\put(68.00,-8.00){\circle*{1.00}}
\qbezier(62,-8)(62,-6)(65,-6)
\qbezier(65,-6)(68,-6)(68,-8)

\put(72.00,-8.00){\circle*{1.00}}
\put(78.00,-8.00){\circle*{1.00}}
\qbezier(72,-8)(72,-6)(75,-6)
\qbezier(75,-6)(78,-6)(78,-8)

\put(82.00,-8.00){\circle*{1.00}}
\put(88.00,-8.00){\circle*{1.00}}
\qbezier(82,-8)(82,-6)(85,-6)
\qbezier(85,-6)(88,-6)(88,-8)

\put(92.00,-8.00){\circle*{1.00}}
\put(98.00,-8.00){\circle*{1.00}}
\qbezier(92,-8)(92,-6)(95,-6)
\qbezier(95,-6)(98,-6)(98,-8)

\put(102.00,-8.00){\circle*{1.00}}
\put(108.00,-8.00){\circle*{1.00}}
\qbezier(102,-8)(102,-6)(105,-6)
\qbezier(105,-6)(108,-6)(108,-8)

\put(112.00,-8.00){\circle*{1.00}}
\put(118.00,-8.00){\circle*{1.00}}
\qbezier(112,-8)(112,-6)(115,-6)
\qbezier(115,-6)(118,-6)(118,-8)

\linethickness{0.2pt}
\qbezier(40.25,-10)(40.25,-8)(42,-8)
\qbezier(42,-8)(43.75,-8)(43.75,-10)
\qbezier(46.25,-10)(46.25,-8)(48,-8)
\qbezier(48,-8)(49.75,-8)(49.75,-10)

\qbezier(50.25,-10)(50.25,-8)(52,-8)
\qbezier(52,-8)(53.75,-8)(53.75,-10)
\qbezier(56.25,-10)(56.25,-8)(58,-8)
\qbezier(58,-8)(59.75,-8)(59.75,-10)

\qbezier(60.25,-10)(60.25,-8)(62,-8)
\qbezier(62,-8)(63.75,-8)(63.75,-10)
\qbezier(66.25,-10)(66.25,-8)(68,-8)
\qbezier(68,-8)(69.75,-8)(69.75,-10)

\qbezier(70.25,-10)(70.25,-8)(72,-8)
\qbezier(72,-8)(73.75,-8)(73.75,-10)
\qbezier(76.25,-10)(76.25,-8)(78,-8)
\qbezier(78,-8)(79.75,-8)(79.75,-10)

\qbezier(80.25,-10)(80.25,-8)(82,-8)
\qbezier(82,-8)(83.75,-8)(83.75,-10)
\qbezier(86.25,-10)(86.25,-8)(88,-8)
\qbezier(88,-8)(89.75,-8)(89.75,-10)

\qbezier(90.25,-10)(90.25,-8)(92,-8)
\qbezier(92,-8)(93.75,-8)(93.75,-10)
\qbezier(96.25,-10)(96.25,-8)(98,-8)
\qbezier(98,-8)(99.75,-8)(99.75,-10)

\qbezier(100.25,-10)(100.25,-8)(102,-8)
\qbezier(102,-8)(103.75,-8)(103.75,-10)
\qbezier(106.25,-10)(106.25,-8)(108,-8)
\qbezier(108,-8)(109.75,-8)(109.75,-10)

\qbezier(110.25,-10)(110.25,-8)(112,-8)
\qbezier(112,-8)(113.75,-8)(113.75,-10)
\qbezier(116.25,-10)(116.25,-8)(118,-8)
\qbezier(118,-8)(119.75,-8)(119.75,-10)
\end{picture}
\caption{An example of the free product}
\end{figure}
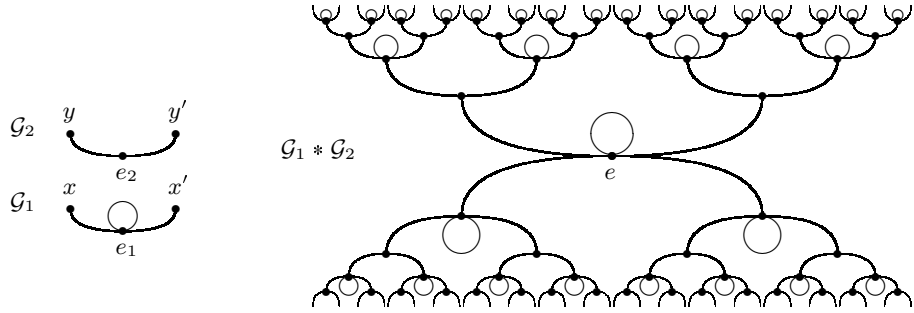
\indent{\par}
\begin{Example}
{\rm The first return moments $N_{\mu_1 \boxtimes \,\mu_2}(n)=N_{\mu_2\boxtimes\, \mu_1}(n)$ 
of lowest orders can be computed using (2.1) and the results of Example 3.2, 
with the $N_{\sigma_1}(n)$, where $\sigma_1=\mu_1\boxslash \mu_2$,
taken from Example 9.2. For the free product in Fig.5, we obtain
$D_{2}(e)=0$, $D_{4}(e)=2$, $D_{6}(e)=0$ and $D_{8}(e)=16$.}
\end{Example}

\end{document}